\documentclass{siamart1116}
\usepackage{amsfonts}
\usepackage{graphicx}
\usepackage{booktabs}
\definecolor{tan}{rgb}{1,0.85,.75}

\newcommand\R{\mathbb{R}}

\DeclareMathOperator{\grad}{grad}
\renewcommand{\div}{\divv}
\DeclareMathOperator{\divv}{div}
\newcommand\x{\times}
\newcommand\dist{\operatorname{dist}}
\newcommand\ct[1]{\multicolumn{1}{c}{#1}}

\numberwithin{theorem}{section}

\headers{Computing spectra without solving eigenvalue problems}{D. Arnold, G. David, M. Filoche, D. Jerison, and S. Mayboroda}

\title{Computing spectra without solving eigenvalue problems\thanks{Submitted to the editors November 12, 2017. This work was supported by grants to each of the authors from the Simons Foundation
(601937, DNA; 601941, GD; 601944, MF; 601948, DJ; 563916, SM)}}

\author{
  Douglas N. Arnold\thanks{School of Mathematics, University of Minnesota, Minneapolis, MN
    (\email{arnold@umn.edu}, \email{svitlana@math.umn.edu}). Arnold was supported by NSF grant DMS-1719694, Mayboroda by NSF INSPIRE grant DMS-1344235 and a Simons Foundation Fellowship.}
  \and
  Guy David\thanks{Univ Paris-Sud, Laboratoire de Math\'ematiques, CNRS, UMR 8658 Orsay, France
    (\email{guy.david@math.u-psud.fr}). Supported by
    an ANR grant, programme blanc GEOMETRYA, ANR-12-BS01-0014.}
  \and
  Marcel Filoche\thanks{Physique de la Mati\`ere Condens\'ee, Ecole Polytechnique, CNRS, Palaiseau, France
    (\email{marcel.filoche@polytechnique.edu}).} 
  \and
  David Jerison\thanks{Mathematics Department, Massachusetts Institute of Technology, Cambridge, MA
    (\email{jerison@math.mit.edu}).  Supported by NSF grant DMS-1500771 and a Simons Foundation Fellowship.}
  \and
  Svitlana Mayboroda\footnotemark[2]
}

\ifpdf
  \DeclareGraphicsExtensions{.eps,.pdf,.png,.jpg}
\else
  \DeclareGraphicsExtensions{.eps}
\fi

\ifpdf
\hypersetup{
  pdftitle={Computing spectra without solving eigenvalue problems},
  pdfauthor={D. Arnold, G. David, M. Filoche, D. Jerison, and S. Mayboroda}
}
\fi

\begin{document}

\maketitle
\begin{abstract}
The approximation of the eigenvalues and eigenfunctions of an elliptic operator is a key computational task in many areas of applied mathematics and computational physics. An important case, especially in quantum physics, is the computation of the spectrum of a Schr\"odinger operator with a disordered potential.  Unlike plane waves or Bloch waves that arise as Schr\"odinger eigenfunctions for periodic and other ordered potentials, for many forms of disordered potentials the eigenfunctions remain essentially localized in a very small subset of the initial domain. A celebrated example is Anderson localization, for which, in a continuous version, the potential is a piecewise constant function on a uniform grid whose values are sampled independently from a uniform random distribution.  We present here a new method for approximating the eigenvalues and the subregions which support such localized eigenfunctions. This approach is based on the recent theoretical tools of the localization landscape and effective potential. The approach is deterministic in the sense that the approximations are calculated based on the examination of a particular realization of a random potential, and predict quantities that depend sensitively on the particular realization, rather than furnishing statistical or probabilistic results about the spectrum associated to a family of potentials with a certain distribution.  These methods, which have only been partially justified theoretically, enable the calculation of the locations and shapes of the approximate supports of the eigenfunctions, the approximate values of many of the eigenvalues, and of the eigenvalue counting function and density of states, all at the cost of solving a single source problem for the same elliptic operator. We study the effectiveness and limitations of the approach through extensive computations in one and two dimensions, using a variety of piecewise constant potentials with values sampled from various different correlated or uncorrelated random distributions.
\end{abstract}



\section{Introduction}

\begin{figure}[htbp]
  \hbox{%
  \vbox{\hsize=.5\textwidth\centerline{\includegraphics[width=.4\textwidth]{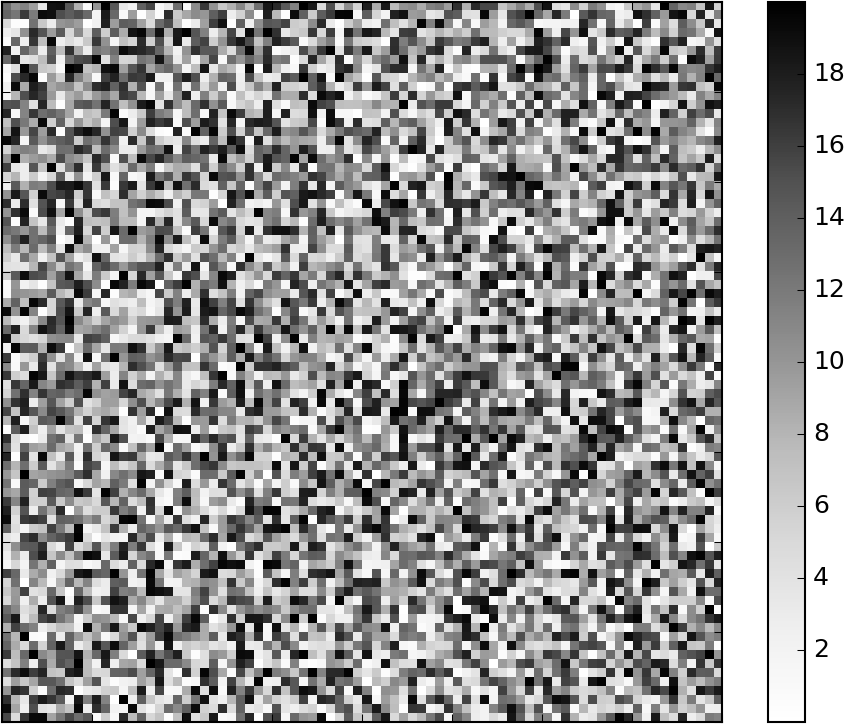}}}%
  \hspace{-.2in}\vbox{\hsize=.5\textwidth
    \centerline{\includegraphics[width=.4\textwidth]{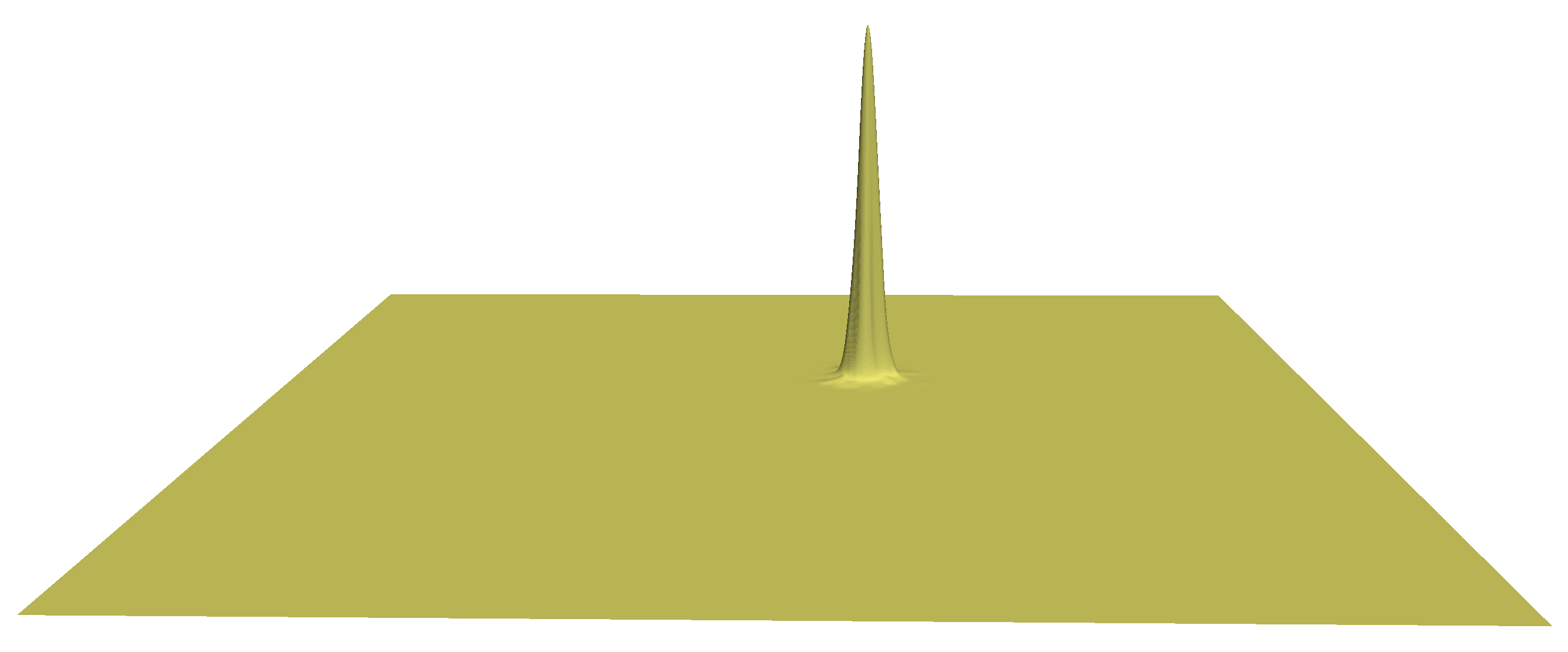}}
    \vspace{0in}
    \centerline{\includegraphics[width=.4\textwidth]{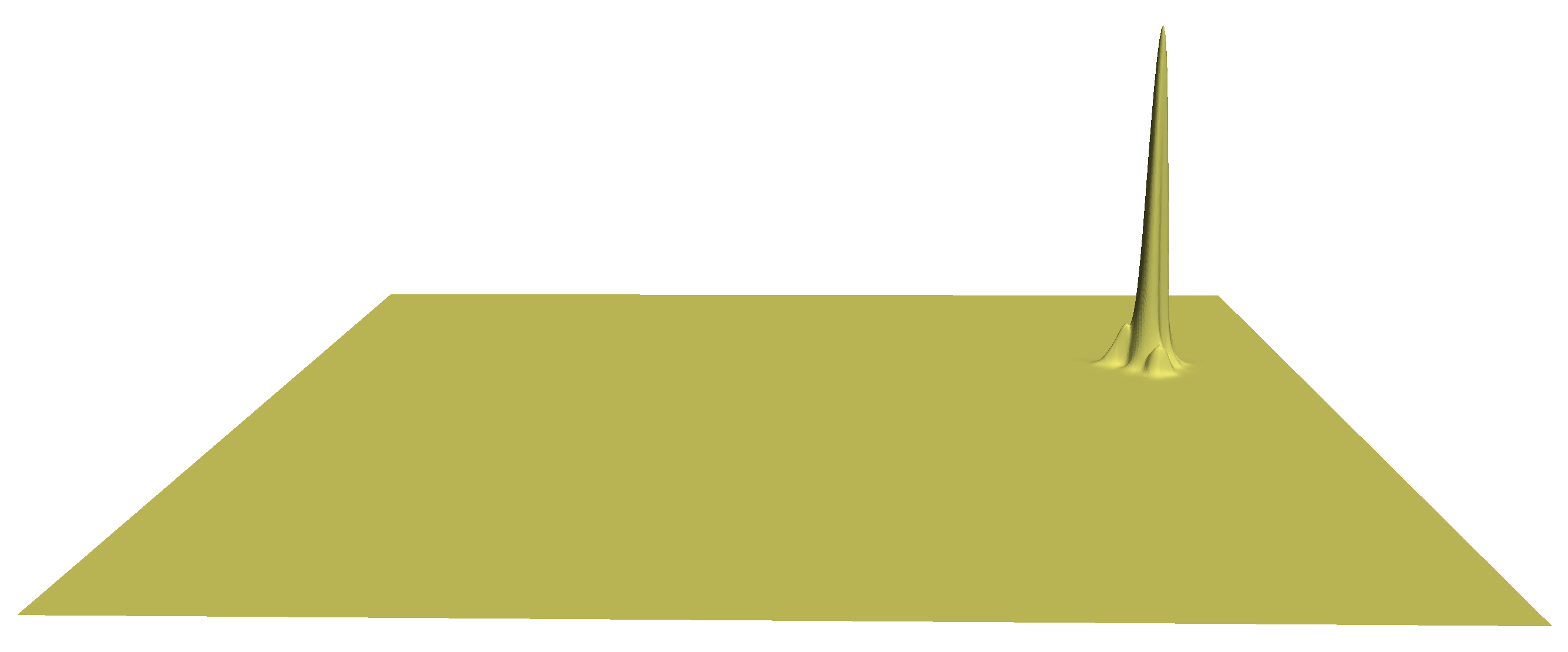}}}
 }%
 \caption{A piecewise constant potential with randomly chosen iid values, and
 the first two eigenfunctions of the corresponding Schr\"odinger operator.}
\label{fg:localization} 
\end{figure}

Eigenfunctions of elliptic operators are often widely dispersed throughout the domain. For example, the eigenfunctions of the Laplacian on a rectangle are tensor products of trigonometric functions, while on a disk they vary trigonometrically in the angular variable and as Bessel functions in the radial argument.  By contrast, in some situations eigenfunctions of an elliptic operator \emph{localize}, in the sense that they are practically zero in much of the domain (after normalizing the $L^2$ or $L^\infty$ norm, say, to $1$), like the two functions pictured on the right of Figure~\ref{fg:localization} (which will be explained shortly). Localization may be brought about by different mechanisms including irregular coefficients of the elliptic operator, certain complexities of the geometry of the domain such as thin necks or fractal boundaries, confining potential wells, and disordered potentials. A celebrated example is Anderson localization~\cite{Anderson1958}, which refers to localization of the eigenfunctions of the Schr\"odinger operator on $\R^n$ 
induced by a potential
with random values. 
The eigenfunctions of the Anderson system model the quantum states of an electron in a disordered alloy, and localization can even trigger a transition of the system from metallic to insulating behavior.  Over the past 60~years, 
analogous phenomena have been observed in
many other fields, and found numerous applications to
the design of optical~\cite{Riboli2011}, acoustic~\cite{Sapoval1997, Felix2007}, electromagnetic~\cite{Laurent2007, Sapienza2010}, and photonic devices~\cite{Filoche2017, Li2017}.

The localized functions shown on the right of Figure~\ref{fg:localization} are the first and second eigenfunctions of the Schr\"odinger operator $H=-\Delta + V$ on the square $[0, 80]\times[0, 80]$ with periodic boundary conditions. The potential $V$ is of Anderson type:
it is a piecewise constant potential obtained by dividing the domain into $80^2$ unit subsquares and assigning to each a constant value chosen randomly from the interval $[0, 20]$.

Although we will consider numerous disordered potentials generated by random distributions in this paper, we emphasize that
our approach is deterministic.  We aim to efficiently predict aspects of the spectrum that depend
on a particular realization of the potential, rather than to give statistical or probabilistic
results about the spectrum associated to a family of potentials with a certain distribution.  For example, the location of the eigenfunctions shown in Figure 1---where they achieve their maxima, what is the shape of their effective supports---depends sensitively on the precise configuration of the disordered potential. The determination of characteristics of the spectrum such as these is an example of the issues addressed here.

In this paper we shall focus on the Schr\"odinger operator $H=-\Delta + V$. The domain $\Omega$ will always be either an interval in one dimension or a square in two, the boundary conditions will be periodic, and the potential $V$ will be piecewise constant with respect to a uniform mesh of $\Omega$ having positive values chosen from some probability distribution, either independently or with correlation. Much of the work can be extended, e.g., to more general domains, potentials, PDEs, and boundary conditions. In particular, we remark that the choice of periodic boundary conditions is for simplicity, and that similar localization occurs with Neumann or Dirichlet boundary conditions. We generally choose unit-sized subsquares for the constant regions of the potential, but this is merely a convenient normalization. For example, instead of considering an $80\times 80$ square broken into unit subsquares as the domain in Figure~\ref{fg:localization}, we could have chosen instead to take the unit square as the domain, with subsquares of side length $1/80$ for the potential. Had we scaled the potential to take values in the range from $0$ to $128,000$, we would have obtained the same localization ($128,000$ being $20\times 80^2$).


There is a large literature concerning the localization of eigenfunctions, approaching the phenomenon from various viewpoints: spectral theory, probability, and quantum mechanics.  But there is nothing like a complete explanatory and predictive theory which can quantitatively and deterministically answer such basic questions as:
\begin{itemize}
 \item How are the eigenvalues and eigenfunctions determined by a particular potential?
 \item Given a potential, do the eigenfunctions localize, and, if so, how many?
 \item What are the size, shape, and location of the approximate supports and how do the eigenfunctions decay away from them?
 \item What are the associated eigenvalues?
\end{itemize}
The present work aims at providing answers to these questions, based on a theory of localization recently conceived and under active development~\cite{Filoche2012, Filoche2013, ADJMF2016, Lefebvre2016, ADFJM2017}.

In the next section of the paper, we briefly survey some classical 
tools used to understand localization. Then, in Section~\ref{sec:lfep},
we introduce two more recent tools: the 
localization landscape function and its associated effective potential, 
introduced in~\cite{Filoche2012, ADJMF2016}. These are easily defined.  The 
landscape function $u$ is the solution to $Hu=1$ subject to periodic boundary conditions, and 
the effective potential $W$ is its reciprocal. Our approach is guided by estimates and relations 
between these objects and the spectrum of the
Schr\"odinger operator from \cite{Filoche2012}, \cite{ADJMF2016}, and the
recent theoretical work of \cite{ADFJM2017}.  In Section~\ref{sec:efun} we show how the structure of wells 
and barriers of the effective potential can be incorporated into numerical algorithms 
to predict the locations and approximate supports of localized 
eigenfunctions.  Then, in Section~\ref{sec:eval}, we show how the values of the minima of the effective potential can 
be used to predict the corresponding eigenvalues and density of states.
Throughout, the  performance of our algorithms is demonstrated for various types of 1D and 2D 
random piecewise constant potentials (uniform, Bernoulli, Gaussian, uncorrelated 
and correlated).

Note that the computation of the effective potential involves the solution of
a single source problem, and so is far less demanding than the computation of a significant portion
of the spectrum by traditional methods. However, a remarkable conclusion of our results is that, for the localized problems
studied here, a great deal of information about the spectrum can be extracted easily from
the effective potential.

\section{Classical confinement}\label{sec:cc}

A simple and well-understood example of 
eigenfunction localization for the Schr\"odinger equation occurs with a classically confining potential. Such a potential is decisively different and simpler than the highly disordered potentials we consider, but we shall draw a connection in the forthcoming discussion, and for that reason we now briefly review localization by confinement. The basic example is
the finite square well potential in one dimension, for which the analytic solution is derived in first courses on quantum mechanics~\cite[vol.~1, p.~78]{Messiah1967}. This problem is posed on the whole real line with the potential $V$ equal to some positive number $\nu$ for $|x|>1$ and vanishing otherwise. The fundamental eigenfunction is then
\begin{equation}
\psi(x) = \begin{cases}
        \cos (\sqrt\lambda x), & |x|\le 1,\\
        c_\nu \exp(-\sqrt{\nu-\lambda}|x|), & |x|>1,
       \end{cases}
\end{equation}
where the eigenvalue $\lambda$ is uniquely determined as the solution to the equation $\cos \sqrt\lambda = \sqrt{\lambda/\nu}$ in the interval $(0,\pi^2/4)$
and $c_\nu=\sqrt{\lambda/\nu}\,\exp(\sqrt{\nu-\lambda})$. Since $\lambda<\nu$, the solution decays exponentially as $|x|\to\infty$, which captures localization in this context. A similar calculation can be made in higher dimensions, for example for spherical wells~\cite[vol.~1, p.~359-361]{Messiah1967}, in which case the well height $\nu$ must be sufficiently large to ensure that there exists an eigenvalue smaller than $\nu$.

A fundamental result for the localization of eigenfunctions of the Schr\"odinger equation with a general confining potential is Agmon's theory \cite{Agmon1982}, \cite[\S~3.3]{Helffer1988}, \cite[Ch.~3]{HislopSigal}, which demonstrates a similar exponential decay for a much larger class of potentials. In this case the domain is all of $\R^n$ and the Schr\"odinger operator is understood as an unbounded operator on $L^2$. One requires that the potential be sufficiently regular and bounded below and that there is an eigenfunction $\psi$ with eigenvalue $\lambda$ such that $V>\lambda$ outside a bounded set. In other words, outside a compact potential well where $V$ dips below the energy level $\lambda$, it remains above it (thus creating confinement). Agmon defined an inner product on the tangent vectors at a point $x\in\R^n$ by
\begin{equation}\label{agmonmetric}
\langle \xi,\eta\rangle_x = \sqrt{[V(x)-\lambda]_+}\,\xi\cdot\eta,
\end{equation}
where the subscript $+$ denotes the positive part.
This defines a Riemannian metric, except that it degenerates to zero at points $x$ where $V(x)\le \lambda$. Its geodesics define a (degenerate) distance $\dist^V_\lambda(x, y)$ between points $x,y\in\R^n$, and, in particular, we may define $\rho(x)=\dist^V_\lambda(x,0)$ to be the distance from the origin to $x$ computed using the Agmon metric. Agmon's theorem states, with some mild restrictions on the regularity on $V$, that for any $\epsilon>0$,
\begin{equation}\label{agmon}
\int_{\R^n} |e^{(1-\epsilon)\rho(x)}\psi|^2\,dx <\infty.
\end{equation}
This result describes exponential decay of the eigenfunction in an $L^2$ sense in regions where $\rho(x)$ grows, which expresses localization in this context.

\begin{figure}[htbp]
  \centerline{%
   \includegraphics[width=35mm]{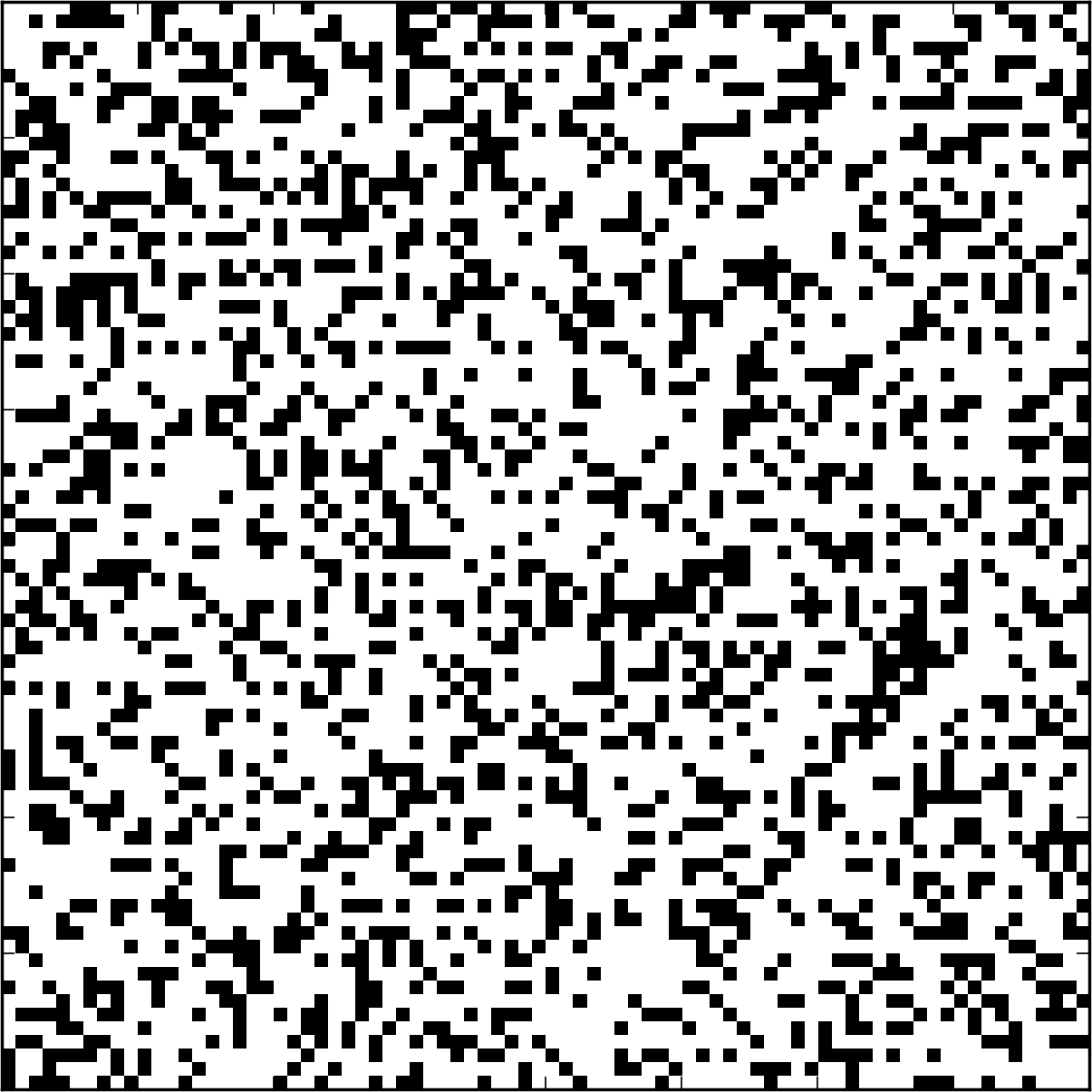}%
  \qquad
   \includegraphics[width=35mm]{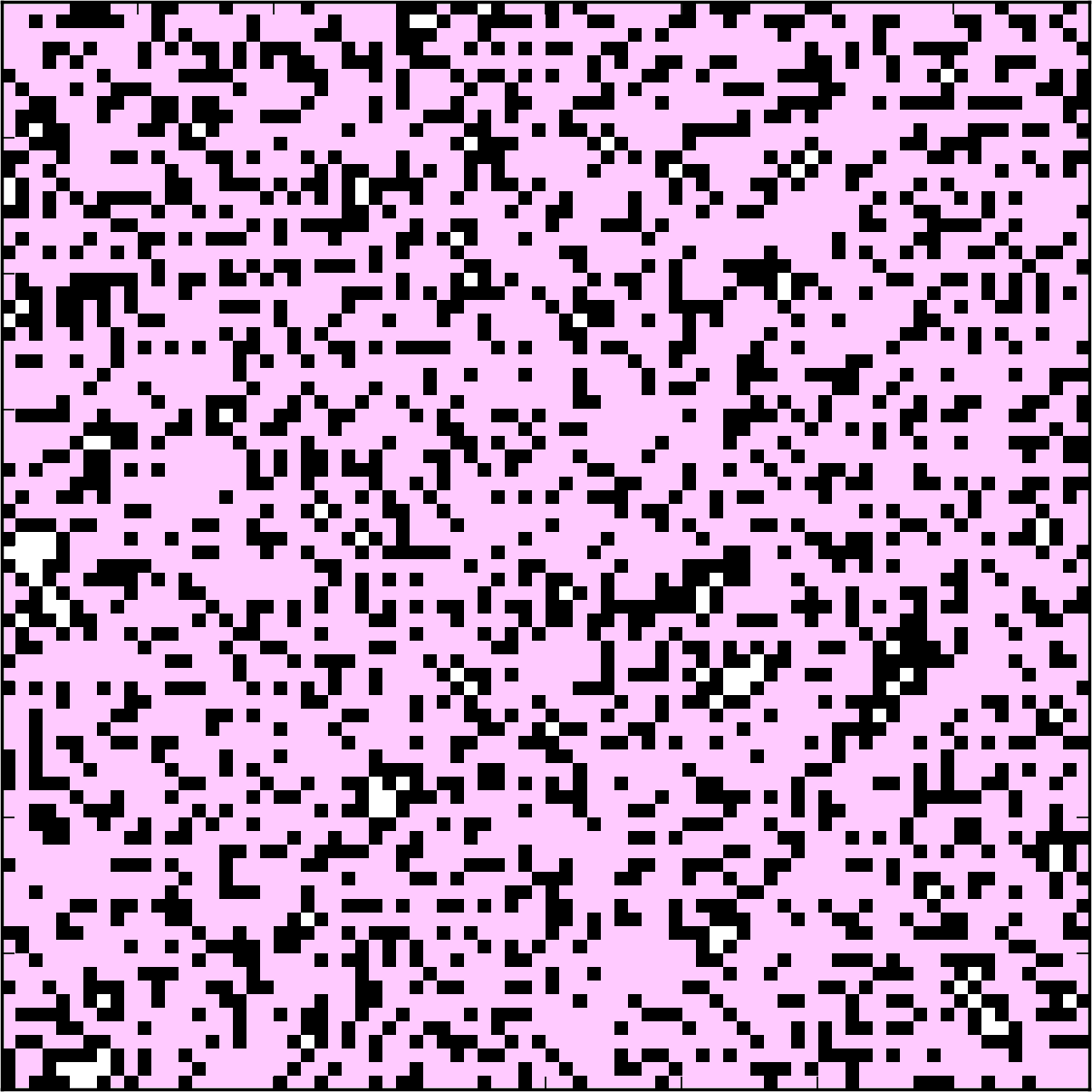}}
   \medskip
  \centerline{%
   \includegraphics[width=45mm]{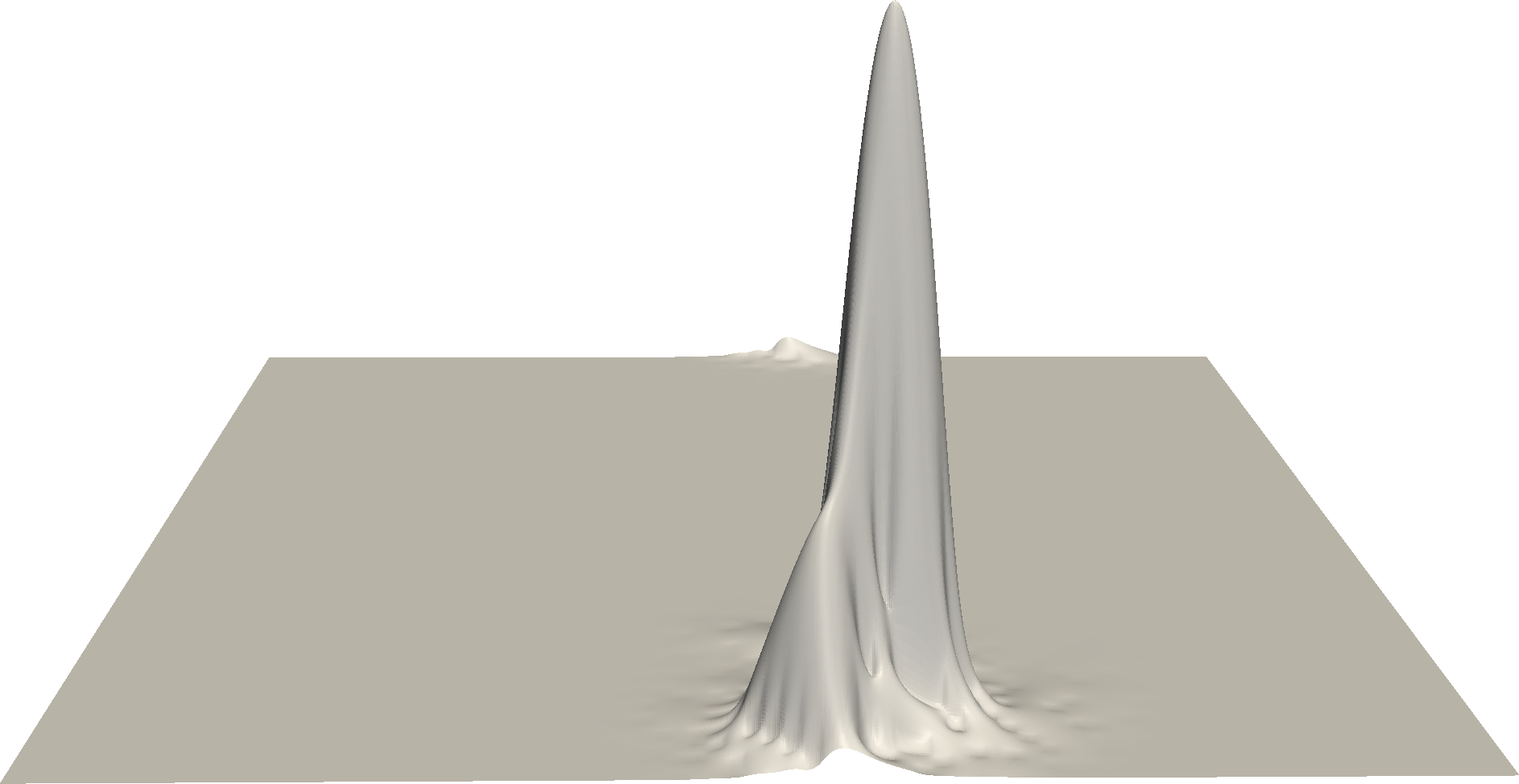}%
  \qquad
   \includegraphics[width=45mm]{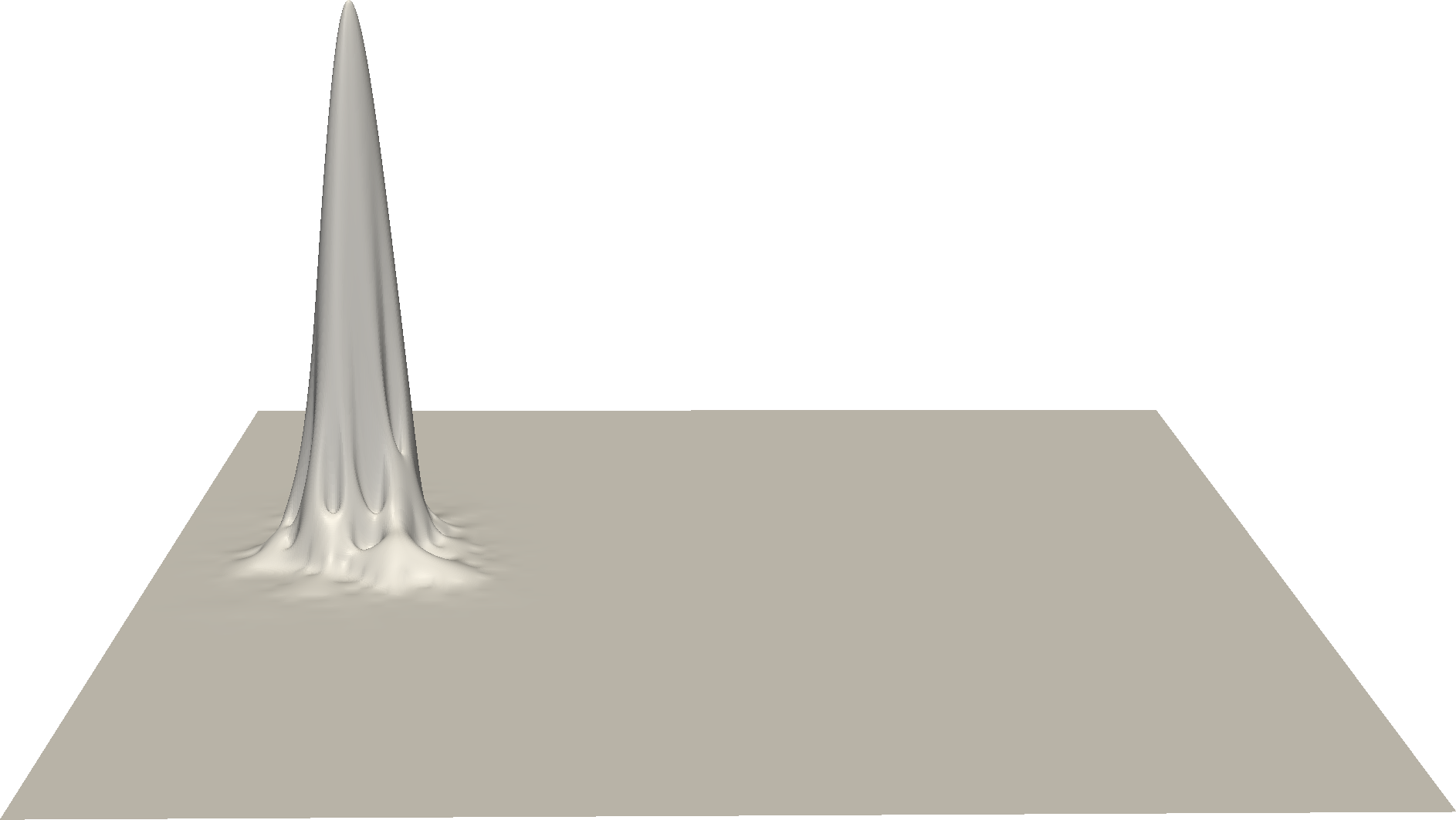}}
 \caption{For this Bernoulli potential there are no wells surrounded
 by thick walls.  Nonetheless, the eigenfunctions localize.}
\label{fg:bernoulli} 
\end{figure}

For the random potentials which we investigate in this paper, the Agmon distance is
highly degenerate, and, consequently, an estimate like \eqref{agmon} is not generally useful. Consider, as a clear example, the Bernoulli potential shown in Figure~\ref{fg:bernoulli}, in which the values $0$ and $4$ are assigned randomly and independently to each of the $80\x80$ subsquares with probabilities $70\%$ and $30\%$, respectively. As shown in color on the right of Figure~\ref{fg:bernoulli}, the region where the potential is zero has a massive connected component which nearly exhausts it. For any positive $\lambda$, the Agmon distance between any two points in this connected component is zero, and hence an estimate like \eqref{agmon} tells us nothing. Nonetheless, as exemplified by the first two eigenfunctions shown in the figure, the eigenfunctions do localize, a phenomenon for which we must seek a different justification.

\section{The landscape function and the effective potential}\label{sec:lfep}

An important step forward was made in \cite{Filoche2012} with the introduction of the \emph{landscape function}, which is simply the solution to the PDE $Hu=1$ together with, for us, periodic boundary conditions. Note that, as long as the potential $V$ is a bounded nonnegative function, nonzero on a set of positive measure, then $u$ is a strictly positive periodic $C^1$ function (indeed, it belongs to $W^2_p$ for any $p<\infty$). The following estimate, taken from \cite{Filoche2012}, relates the landscape function to the eigenvalues.
\begin{proposition}\label{th:landscape-ineq}
 If $\psi:\Omega\to\R$ is an eigenfunction of $H$ with eigenvalue $\lambda$, then
\begin{equation}\label{le}
 \psi(x) \le \lambda u(x) \|\psi\|_{L^\infty}, \quad x\in\Omega.
\end{equation}
\end{proposition}

\begin{figure}[htbp]
  \centerline{%
  \includegraphics[height=2in]{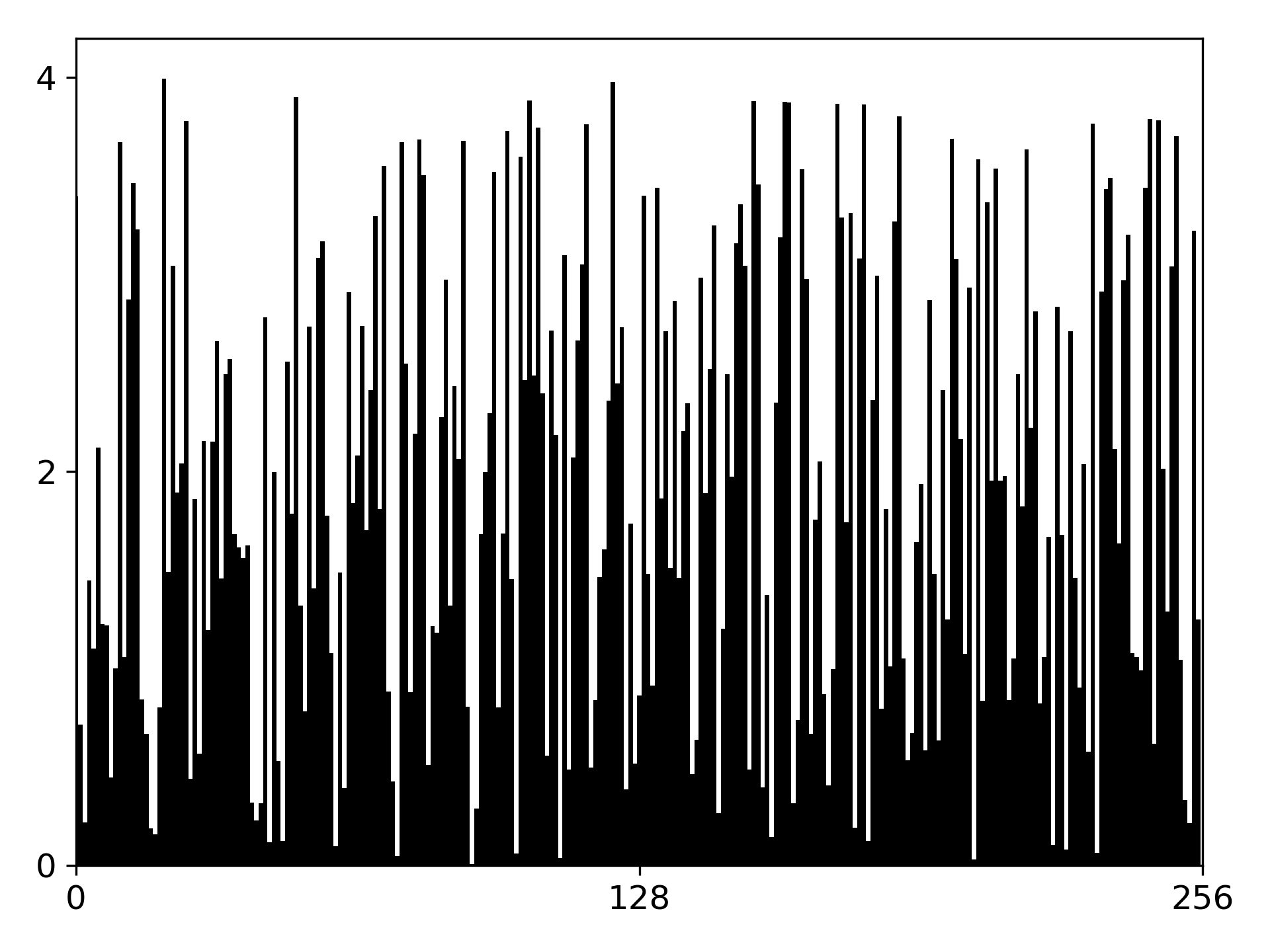}%
  \quad
  \includegraphics[height=2in]{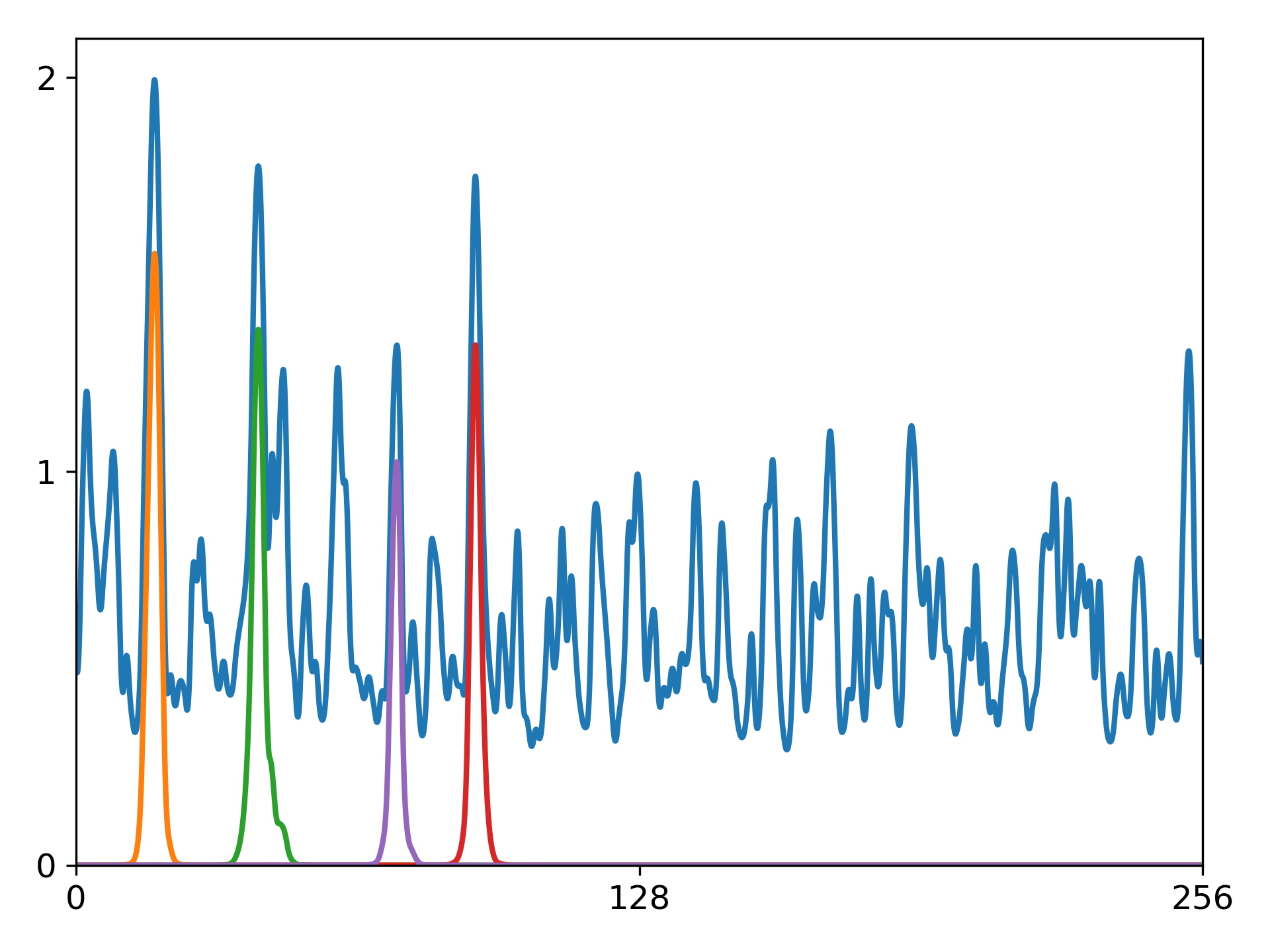}}
 \caption{The potential on the left gives rise to the landscape function, plotted in blue on the right. The first four eigenfunctions are also plotted, scaled so that their maximum value is the reciprocal of their eigenvalue, illustrating the inequality of \eqref{le}.}
\label{fg:le} 
\end{figure}

If we normalize the eigenfunction $\psi$ so that $\|\psi\|_{L^\infty}=1/\lambda$, then the theorem asserts that $\psi\le u$ pointwise, a fact illustrated in Figure~\ref{fg:le}, which shows the one-dimensional case where the potential has 256~values randomly chosen uniformly iid from $[0, 4]$. The argument was made in~\cite{Filoche2012} that if the landscape function $u$ nearly vanishes on the boundary of a subregion $\Omega_0$ of the domain $\Omega$, then \eqref{le} implies that any eigenfunction $\psi$ must nearly vanish there as well, and so $\psi|_{\Omega_0}$ is nearly a Dirichlet eigenfunction for $\Omega_0$ (with the same eigenvalue $\lambda$).  Similarly $\psi$ restricts to a near Dirichlet eigenfunction on the subdomain complementary to $\Omega_0$. Except for the unlikely case in which these two subdomains share the eigenvalue $\lambda$, this suggests that $\psi$ must nearly vanish in one of them, and so be nearly localized to the other. This viewpoint gives an initial insight
into the situation illustrated in Figure~\ref{fg:le}, where it is seen that the eigenfunctions are essentially localized to the subdomains between two consecutive local minima of $u$. However, it must be remarked, in the case shown in Figure~\ref{fg:le} and many other typical cases, the landscape function merely dips, but in no sense vanishes, on the boundary of localization regions, 
and so a new viewpoint is needed in order to satisfactorily explain the localization which is observed.

Such a new viewpoint was developed in the paper \cite{ADJMF2016}, where the emphasis was placed
on the \emph{effective potential} $W$, defined as the reciprocal $1/u$ of the landscape function. A key property of $W$ and the explanation for its name, is that it is the potential of an elliptic operator which is conjugate to the Schr\"odinger operator $H$ and so has the same spectrum.  The following essential identity was derived in \cite{ADJMF2016}.
\begin{proposition}\label{th:ep}
 Suppose that the potential $V\in L^\infty(\Omega)$ is nonnegative and positive on a set
 of positive measure. Let $u>0$ be
 the landscape function and $W=1/u$ the effective potential.
 Define $L:H^1(\Omega)\to H^{-1}(\Omega)$ by
\begin{equation}\label{defL}
 L\phi = -\frac1{u^2}\div(u^2 \grad\phi).
\end{equation}
Then
$$
 (-\Delta + V)(u\phi) = u\,(L + W)\phi,\quad \phi\in H^1(\Omega).
$$
\end{proposition}
In this result, $H^1$ denotes the periodic Sobolev space on $\Omega$ and $H^{-1}$ its dual.  The equation holds in $H^{-1}$, making sense because $u\in C^1$.

%

\begin{corollary}\label{cr:eigfns} Let $V$, $u$, and $W$ be as in the theorem and $\lambda\in\R$.
Then $\psi\in H^1$  satisfies $(-\Delta+V)\psi=\lambda\psi$
if and only if $\phi := \psi/u$ satisfies $(L + W)\phi = \lambda\phi$.
\end{corollary}

\begin{figure}[htbp]
  \centerline{%
   \includegraphics[height=65mm]{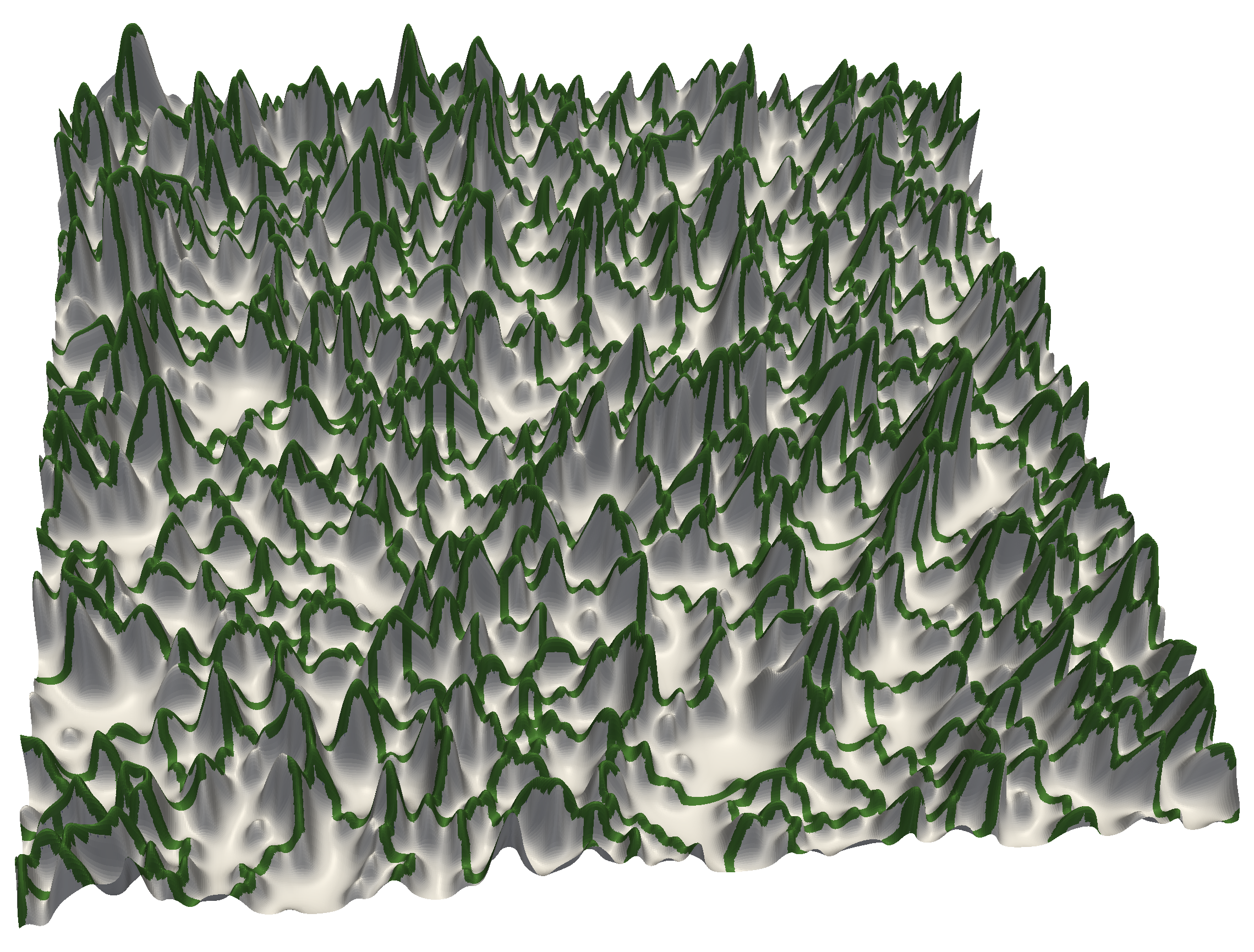}}%
 \caption{The effective potential associated to the Bernoulli potential of Figure~\ref{fg:bernoulli}. It is shown
 with its crestlines which partition the domain into a few hundred basins of attraction surrounding wells.}
\label{fg:berneffpot} 
\end{figure}

Thus the eigenvalues of the operator $L+W$ (with periodicity) are the same as those of the original Schr\"odinger operator, and the eigenfunctions are closely related. However, the corresponding potentials $W$ and $V$ are very different. The effective potential $W$ is often much more regular than the physical potential $V$. More importantly, it has a clear structure of wells and walls.  As we shall see, these induce a sort of localization by confinement, which is not evident from the physical potential. For example, for the Bernoulli potential shown in Figure~\ref{fg:bernoulli}, the effective potential is shown in Figure~\ref{fg:berneffpot}. Note that the effective potential contains many wells: small regions where its value is low, but surrounded, or nearly surrounded, by crestlines where its values are relatively high. If we think of a gradient flow starting at a generic point in the domain and ending at a local minimum, thereby associating the point to one of the local minima of $W$, the crestlines, displayed in green in Figure~\ref{fg:berneffpot} are the boundaries of the basins of attraction of the local minima. There are several algorithms to compute the precise location of these crestlines. The one used in this case is the watershed transform as described in \cite{BeareLehmann}.

In recent work~\cite{ADFJM2017}, we have established a rigorous connection between the well and wall structure of $W$ and the exponential decay of eigenfunctions.   We define the $W$-distance
$\dist^W_\lambda$, to be the Agmon distance, as defined after \eqref{agmonmetric}
in Section~\ref{sec:cc}, but with the potential $V$ replaced by the effective potential $W$.  Then we show, roughly speaking,
that whenever a well of the effective potential exists with sufficient separation of the well depth from the height of the surrounding barriers, then eigenfunctions $\psi$ of the operator $H$  with eigenvalue $\lambda$ are localized to $\{W\le \lambda\}$ in the sense that they decay exponentially in the $W$-distance associated to the eigenvalue. More precisely, in \cite{ADFJM2017}, we prove:

\begin{theorem}
 Suppose $(\psi,\lambda)$ is an eigenpair of the Schr\"odinger operator $H=-\Delta + V$.
 Let $W$ be the effective potential and $\dist^W_\lambda(x,y)$ the associated $W$-distance.
 Let $\delta>0$,
 $$
 S = \{\,x\in\Omega\,|\, W(x)\le \lambda+\delta\,\},
 $$
 a sublevel set of $W$, and $h(x)$ the $W$-distance from $x$ to $S$.  Then there exists
 a constant $C$ depending only on $\|V\|_{L^\infty}$ and $\delta$ (but not the domain $\Omega$) such that
 $$
 \int_\Omega e^{h(x)}\psi^2\,dx \le C \int_\Omega\psi^2\,dx.
 $$
\end{theorem}

This result can be found, stated in considerably more generality and in sharper form in \cite[Corollary~3.5]{ADFJM2017}.  The dependence of the constant is given there explicitly.  It grows at most linearly
in $\|V\|_{L^\infty}$. The same paper also proves exponential decay for the gradient of the eigenfunctions.  On the other hand these theoretical results do not capture fully the accuracy with which $W$ predicts the behavior of the eigenfunctions. Our numerical results show that  the eigenfunctions typically occupy a single connected component of the set  $W\le \lambda$ and decay exponentially across each green crestline of Figure~\ref{fg:berneffpot}, whereas the theorems do not rule out that a resonance occurs resulting in eigenfunctions that have significant mass in several different components of $W\le \lambda$.

\section{Eigenfunction prediction}\label{sec:efun}

We may apply the theory described in the previous section to predict the location and extent of the supports of localized eigenfunctions. For this we proceed in four steps.

\begin{enumerate}
 \item Compute the landscape function $u$ from the PDE $Hu=1$, and define the effective potential $W=1/u$.
 \item To approximate a desired number of localized eigenfunctions, identify the same number of local minima of $W$, selected in order of increasing minimal value. The location of these minima will be our prediction of the place where the eigenfunctions localize, in the sense that the maximum of the localized eigenfunction will occur nearby there.
 \item To each of the selected local minimum we associate an energy level $E$ given by the minimum value of $W$ in the well times a constant greater than $1$.  The constant is chosen so that $E$ is close to fundamental eigenvalue of the well, or perhaps somewhat larger. (We show in the next section how this may be achieved.)
 \item From the corresponding sublevel set, consisting of all $x$ for which $W(x) \le E$, we compute the connected component which contains the
 selected local minimum.  This is the region we predict to be occupied by the eigenfunction.
\end{enumerate}

%

\begin{figure}[htbp]
  \centerline{%
  \includegraphics[height=2in]{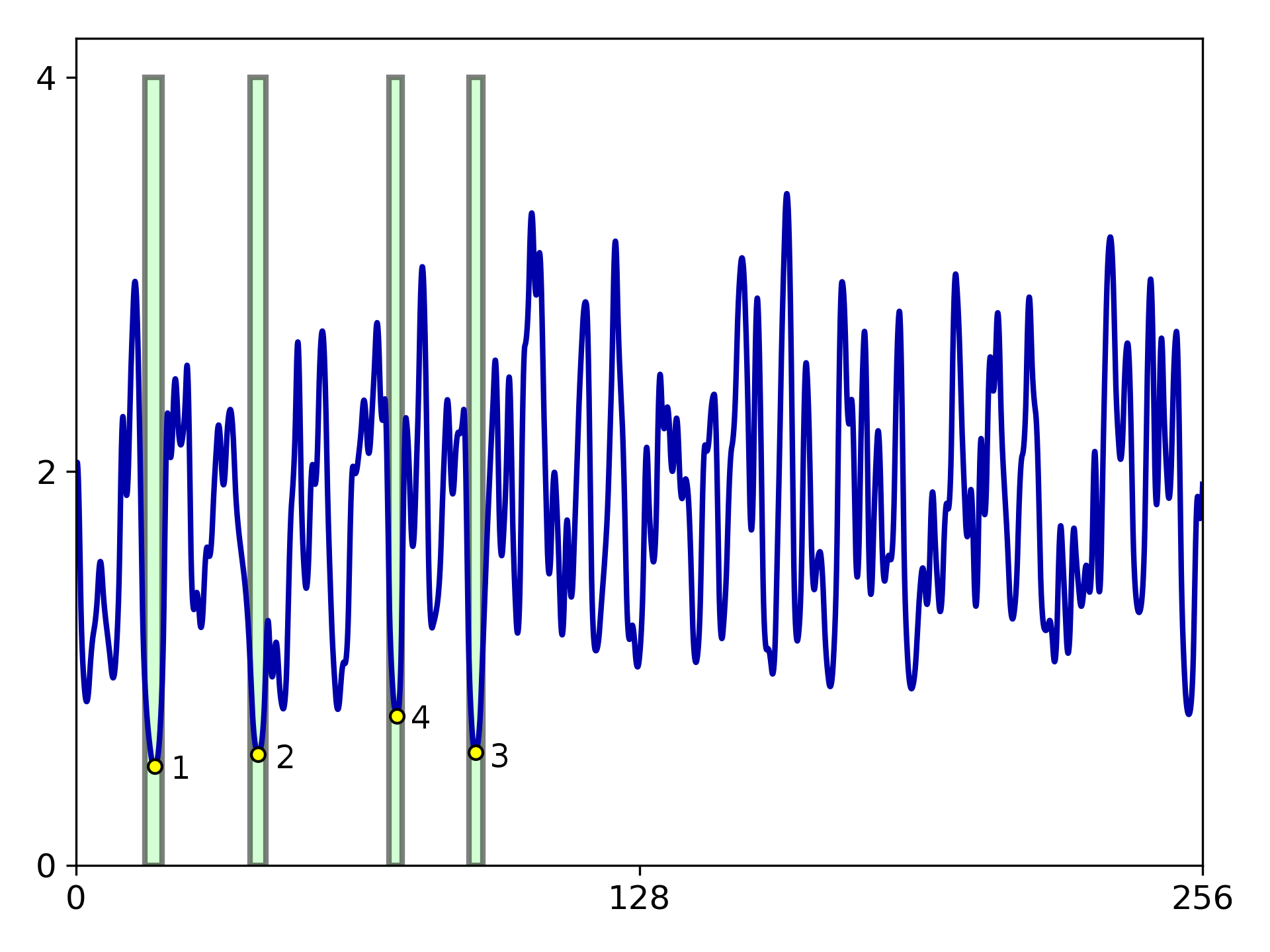}%
  \quad
  \includegraphics[height=2in]{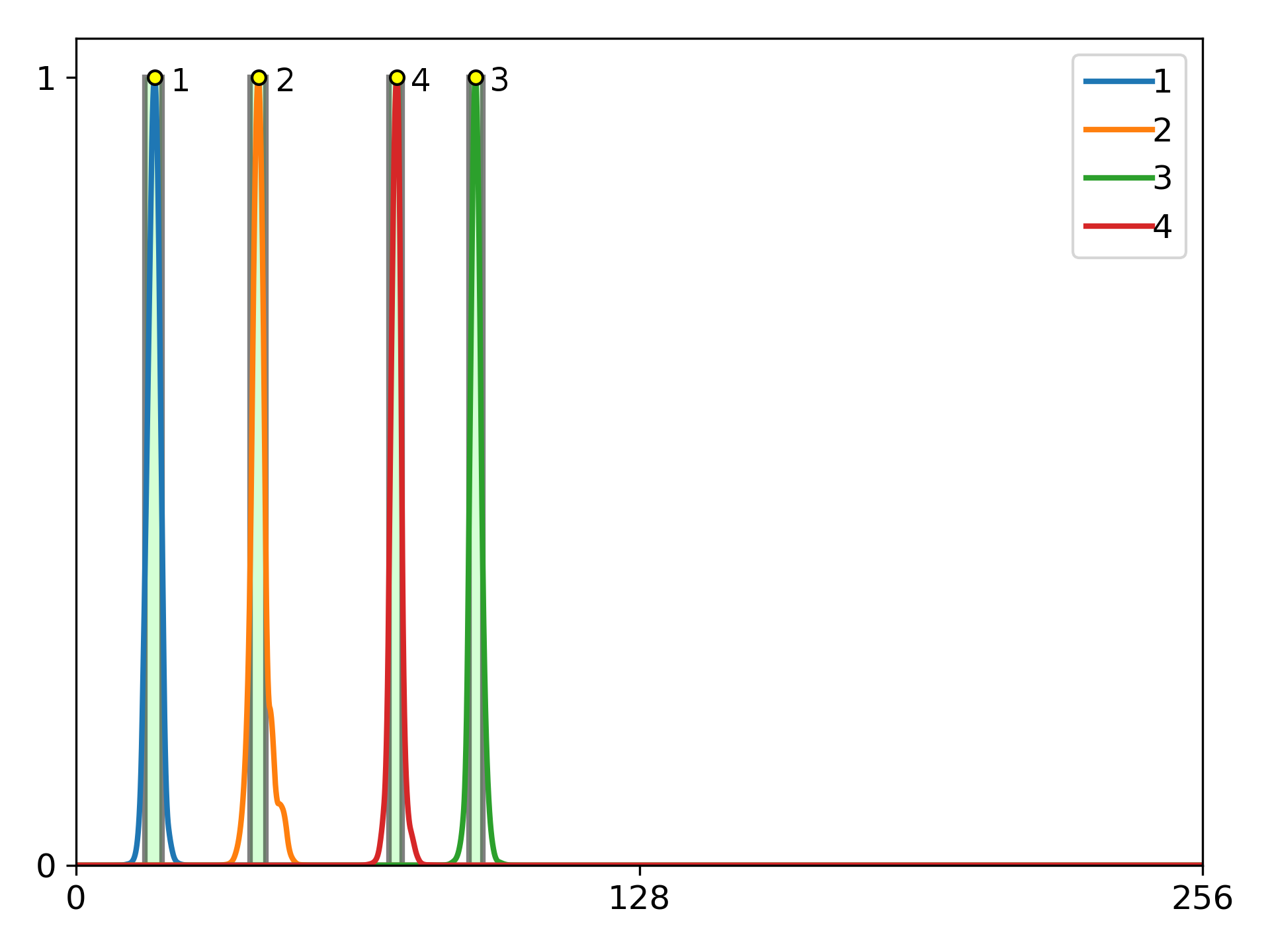}}
 \caption{On the left is the effective potential corresponding to the piecewise constant potential with 256 uniformly iid randomly selected values shown on the left of Figure~\ref{fg:le}. The first, second, third, and fourth deepest local minima are marked and labeled. The small yellow circles signify the positions of these minima. We expect the corresponding eigenfunctions to be centered near the location of the minima, with extent related to the surrounding basin of attraction.  This prediction is plotted in green, and the actual first four eigenfunctions superimposed over the predictions on the right.}
\label{fg:pred1d} 
\end{figure}

Figure~\ref{fg:pred1d} shows the outcome of applying this approach to the 1D Schr\"odinger equation with the potential shown in Figure~\ref{fg:le}. The landscape function was computed using Lagrange cubic finite elements with a uniform mesh of $2,560$ subintervals ($10$ element per constant piece of the potential). The finite element solution was evaluated at $15,360$ equally-spaced points ($6$ per element), with the reciprocals giving the values of the effective potential. The local maxima and local minima of the effective potential were then identified by comparing the value at each point to that of its two immediate neighbors. For comparison, the true eigenvalues and eigenfunctions were computed by using the same finite element discretization and solving the resulting sparse matrix real symmetric generalized eigenvalue problem using a Krylov--Schur solver. The finite element discretizations were implemented using the FEniCS software environment~\cite{Logg2012}, calling the SLEPc library \cite{slepc} for the eigenvalue solves. Note that the locations of the four local minima of the effective potential, indicated by small circles in both plots of Figure~\ref{fg:pred1d}, very nearly coincide with the locations of the maxima of the four corresponding eigenfunctions. Moreover the correspondence respects the ordering of the eigenvalues, in the sense that the $i$th eigenfunction corresponds to the $i$th well for $i=1,2,3,4$. To predict the extent of the localized eigenfunctions, we use as outer boundaries the level curves of the effective potential at an energy level~$E$ set to be $1.875$~times the depth of the wells. This is 150\% of the value we justify in the next section as an approximation of the eigenvalues.  Of course, this
choice is somewhat arbitrary, since the effective support of a localized function is not an absolute notion, but must be defined with respect to some tolerance.  We could as well have chosen a somewhat larger level to get wider regions incorporating more of the tail of the eigenfunctions, or have chosen a somewhat smaller level to get narrower regions.

\begin{figure}[htbp]
  \centerline{%
  \includegraphics[height=2in]{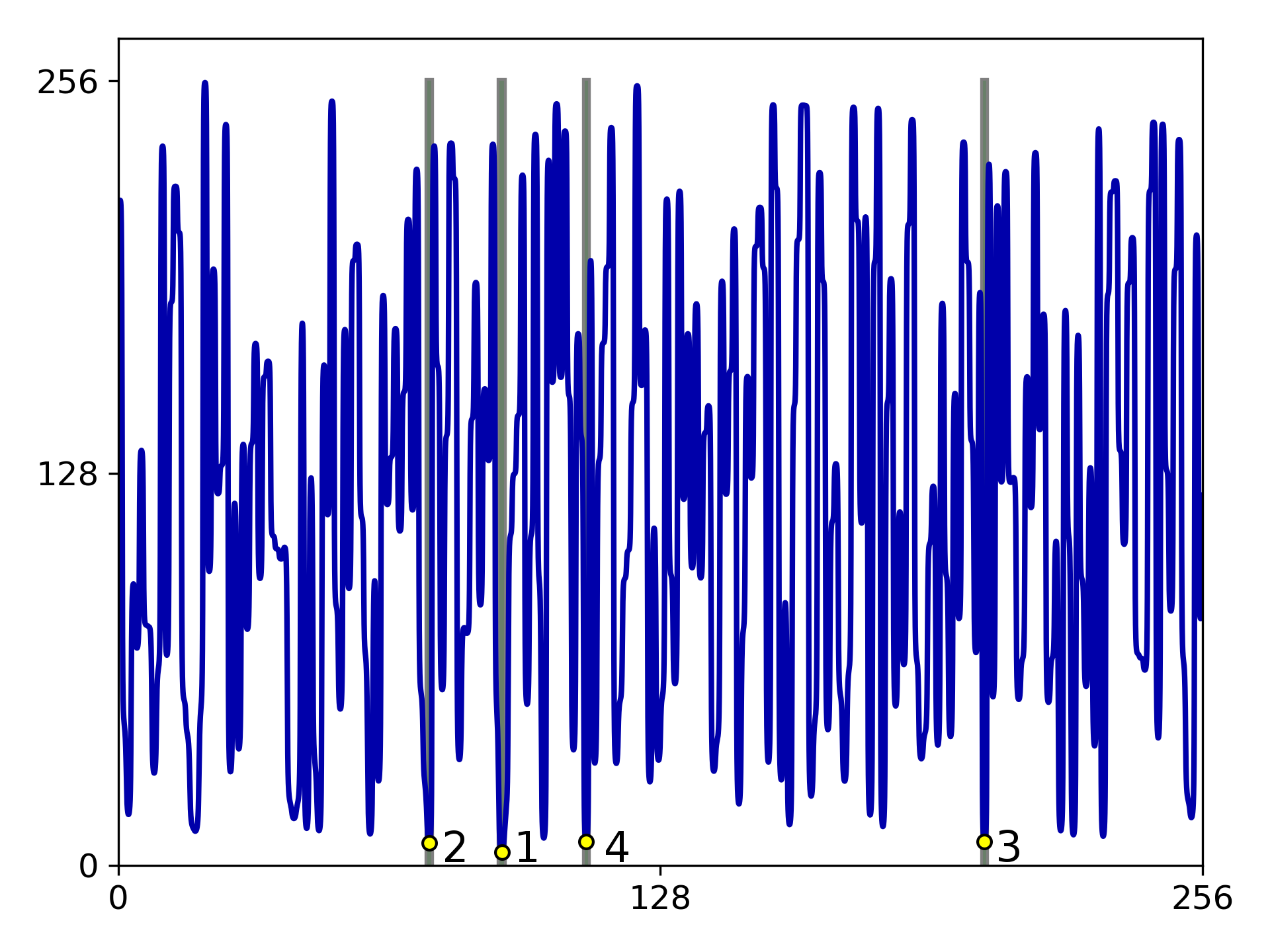}%
  \quad
  \includegraphics[height=2in]{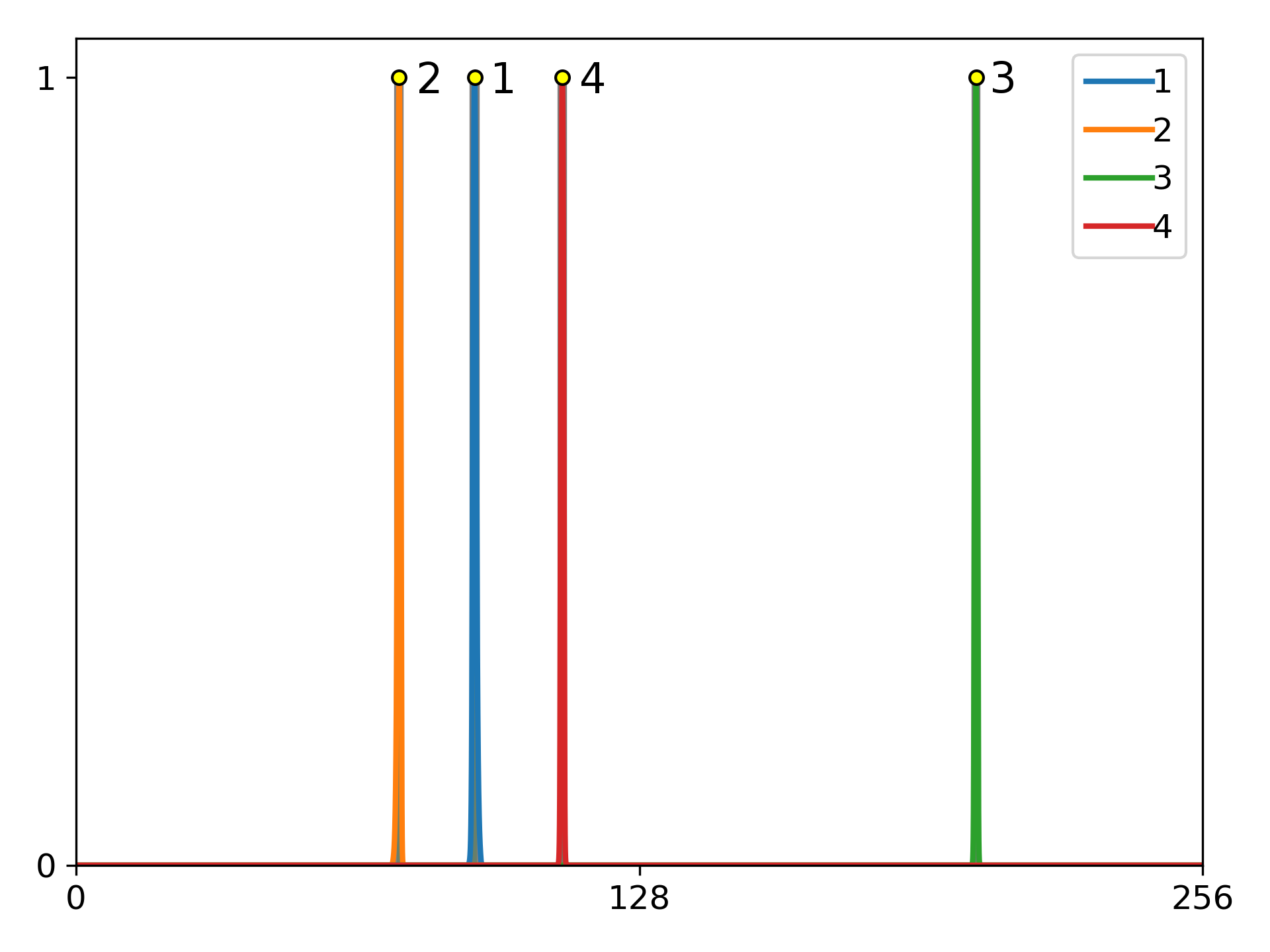}}
 \caption{The effective potential for the potential equal to $64$ times that shown on the left of Figure~\ref{fg:le}, with the four deepest local minima and their wells marked and labeled, and then a comparison with the actual first four eigenfunctions.}
\label{fg:pred1d256} 
\end{figure}

Next, we vary this example by increasing the amplitude of the potential by a factor of $64$, so that it takes values between $0$ and $256$, but is otherwise identical to the potential shown on the left of Figure~\ref{fg:le}. The results analogous to Figure~\ref{fg:pred1d} for this potential are shown in Figure~\ref{fg:pred1d256}. Note that, despite the fact that the potentials are proportional in the two cases, the effective potentials look quite different and the eigenfunctions localize in entirely different places. The eigenvalue maxima occur at $17.85$, $41.47$, $90.75$, and $72.92$ for the smaller potential and at $90.57$, $73.42$, $204.5$, and $110.5$ for the second. These locations again are captured very accurately by the minima of the effective potential, and in the correct order. A crucial difference between the two examples is that the eigenvalues for the problem with the larger potential are much more tightly localized, as predicted by the thinner wells of its effective potential.

\begin{figure}[htbp]
  \centerline{%
  \includegraphics[height=2in]{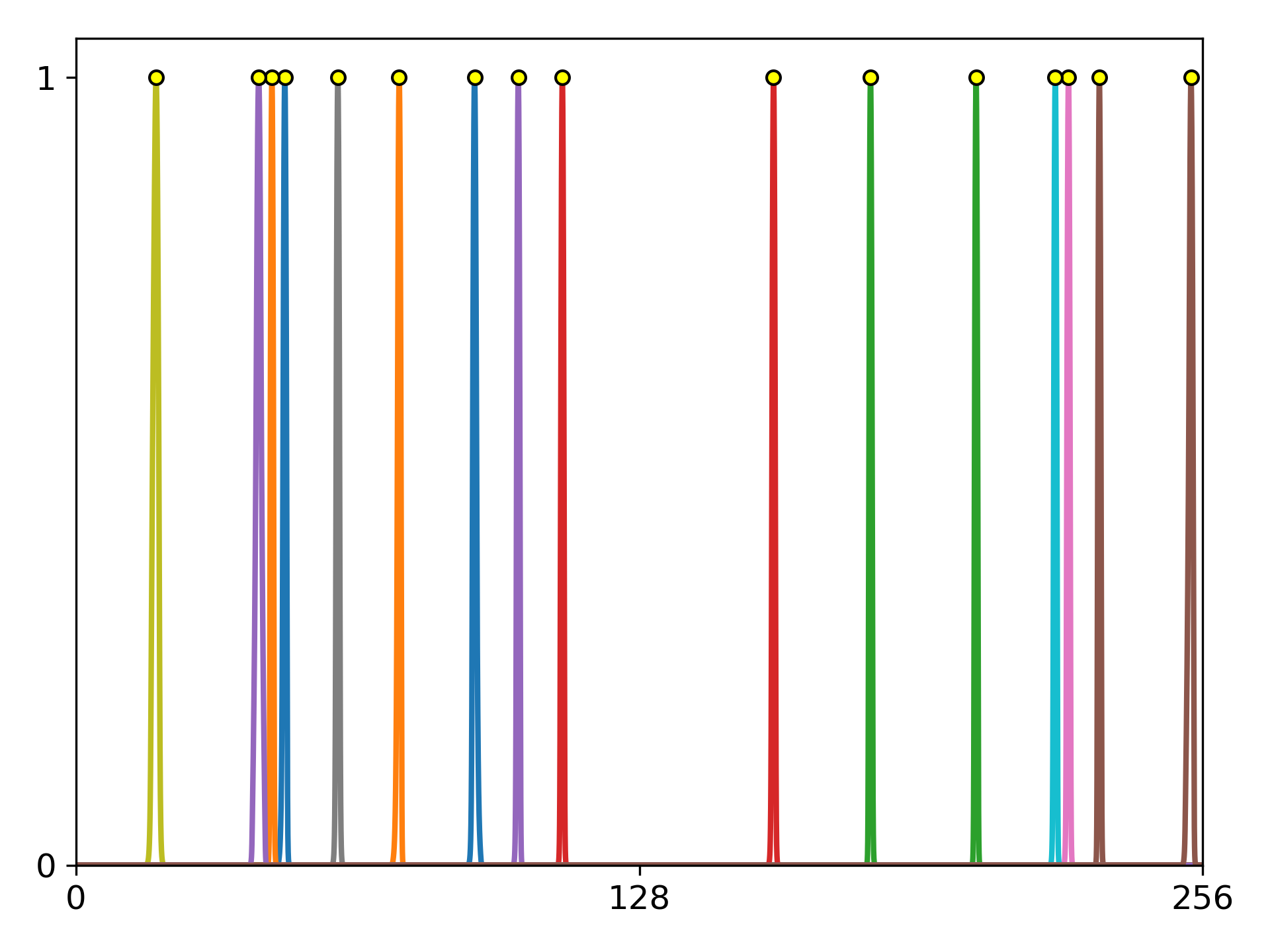}}
 \caption{For the same potential as Figure~\ref{fg:pred1d256} the first 16 local minima
 locations accurately predict the locations of the corresponding eigenfunctions.}
\label{fg:pred1d256-16} 
\end{figure}

\begin{figure}[htbp]
  \centerline{%
  \includegraphics[height=3in]{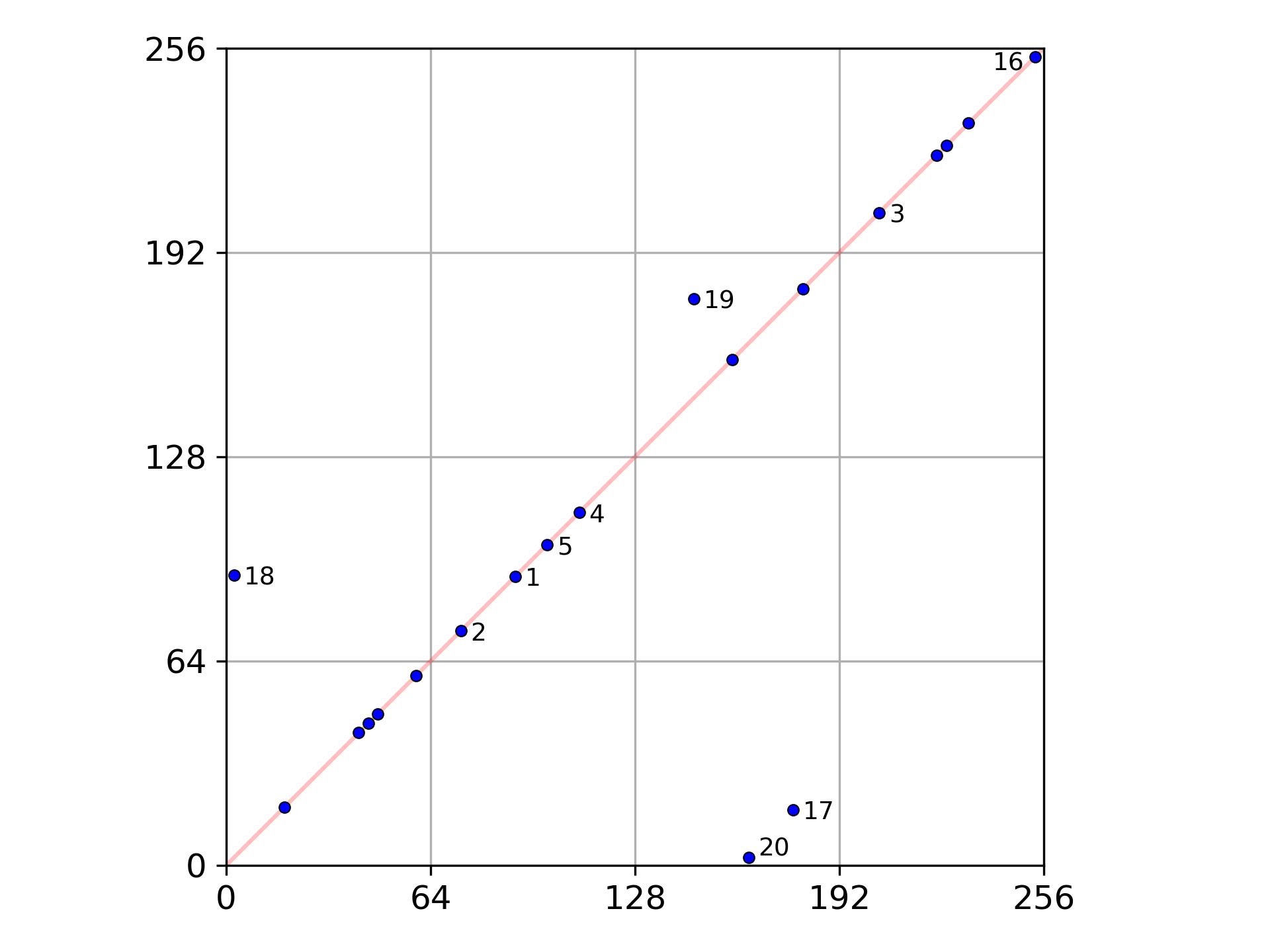}}
 \caption{First 20 local minima on the $x$-axis versus the maximum of the corresponding eigenfunction on the $y$-axis.}
\label{fg:eigsvsmax} 
\end{figure}

In the examples depicted in Figures~\ref{fg:pred1d} and \ref{fg:pred1d256}, we looked at the first four eigenvalues and found that the eigenvalues are very accurately located by the local minima of the effective potential.  We now look at what happens for a larger number of local minima.  Obviously, at some point the eigenfunction locations cannot be predicted by the local minima, since there are infinitely many eigenfunctions and only finitely many local minima.  Figure~\ref{fg:pred1d256-16}, which is similar to the plot on the right-hand side of Figure~\ref{fg:pred1d256}, and, in particular, uses the same potential, shows the locations of the first 16 local minima of $W$ (as yellow dots), plotted over the first 16 eigenfunctions.  We see that for all 16, the location of the local minimum of $W$ predicts very accurately the location of the corresponding eigenfunction.
Figure~\ref{fg:eigsvsmax} and Table~\ref{tb:eigsvsmax} explore the situation further, comparing the location of the $n$th local minimum of $W$, plotted on the $x$-axis, to that of the maximum of the $n$th eigenfunction, plotted on the $y$-axis, for $n=1,\ldots, 20$.  Since these nearly coincide for $n\le 16$, the first 16 points lie very nearly on the line $y=x$.  From then on, however, the points deviate from the line because the ordering of the local minima does not perfectly match the ordering of the most closely associated eigenfunctions. Specifically, as can be seen from Table~\ref{tb:eigsvsmax}, the 17th local minimum of $W$ occurs at the location of the 19th eigenfunction, and the 18th occurs at the location of the 20th eigenfunction.
\begin{table}[htbp]
\centering
\begin{tabular}{rrr@{\hspace{.5in}}rrr}
\toprule
$n$ & $W$ min & eigfn max & $n$ & $W$ min & eigfn max
\\
\midrule
  1 &   90.57 &   90.55 &  11 &   47.43 &   47.43 \\
  2 &   73.42 &   73.43 &  12 &   44.52 &   44.50 \\
  3 &  204.50 &  204.50 &  13 &  180.52 &  180.52 \\
  4 &  110.50 &  110.50 &  14 &  158.48 &  158.48 \\
  5 &  100.48 &  100.48 &  15 &   41.50 &   41.50 \\
  6 &  232.50 &  232.50 &  16 &  253.32 &  253.33 \\
  7 &  225.48 &  225.48 &  17 &  177.50 &   17.43 \\
  8 &   59.48 &   59.48 &  18 &    2.45 &   90.95 \\
  9 &   18.22 &   18.18 &  19 &  146.53 &  177.50 \\
 10 &  222.48 &  222.48 &  20 &  163.58 &    2.47 \\
\bottomrule
\end{tabular}
 \caption{The data plotted in Figure~\ref{fg:eigsvsmax}.}
\label{tb:eigsvsmax} 
\end{table}

\begin{figure}[htbp]
  \centerline{%
  \includegraphics[height=2in]{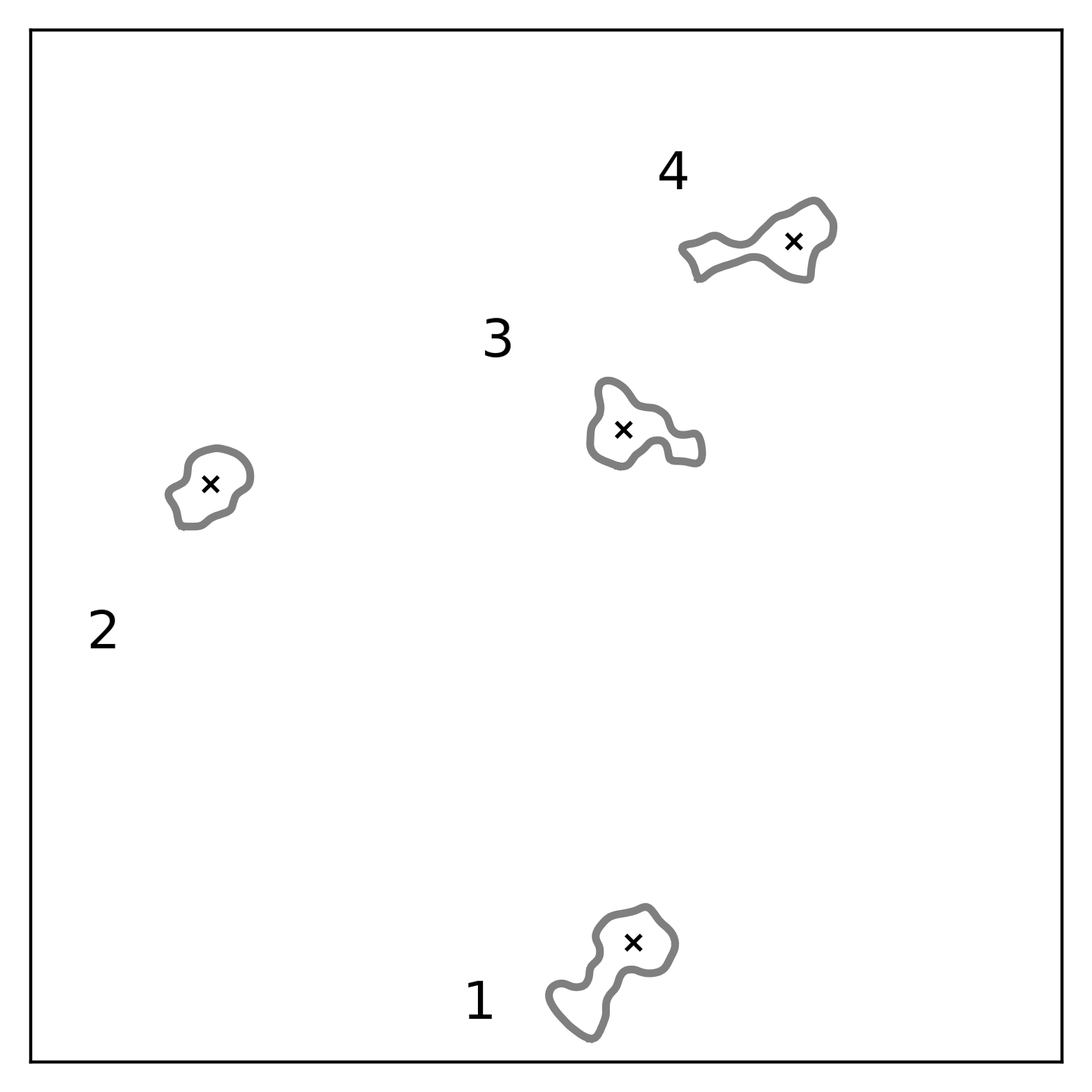}%
  \quad
  \includegraphics[height=2in]{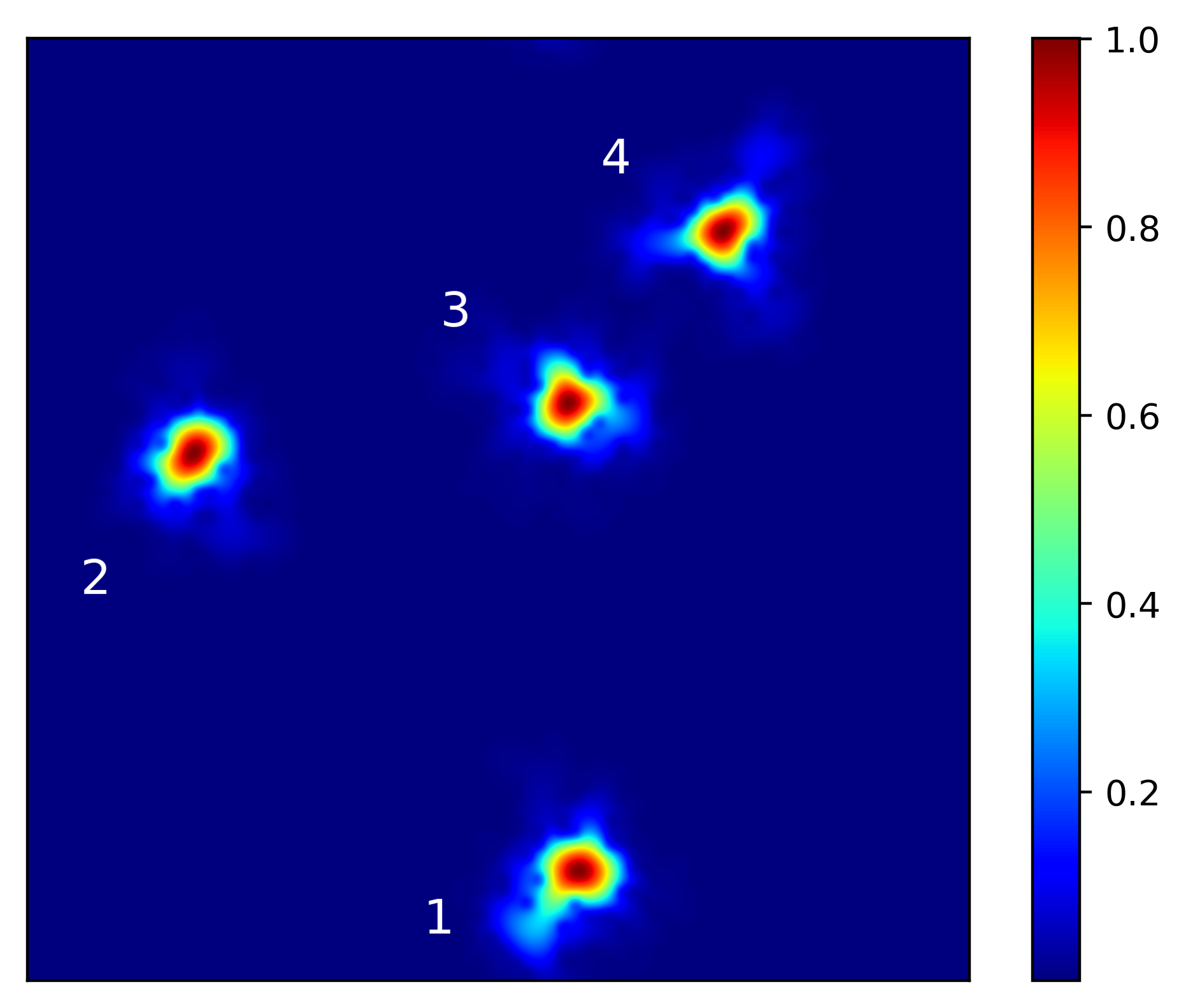}}%
  \medskip
  \centerline{\includegraphics[height=2.5in]{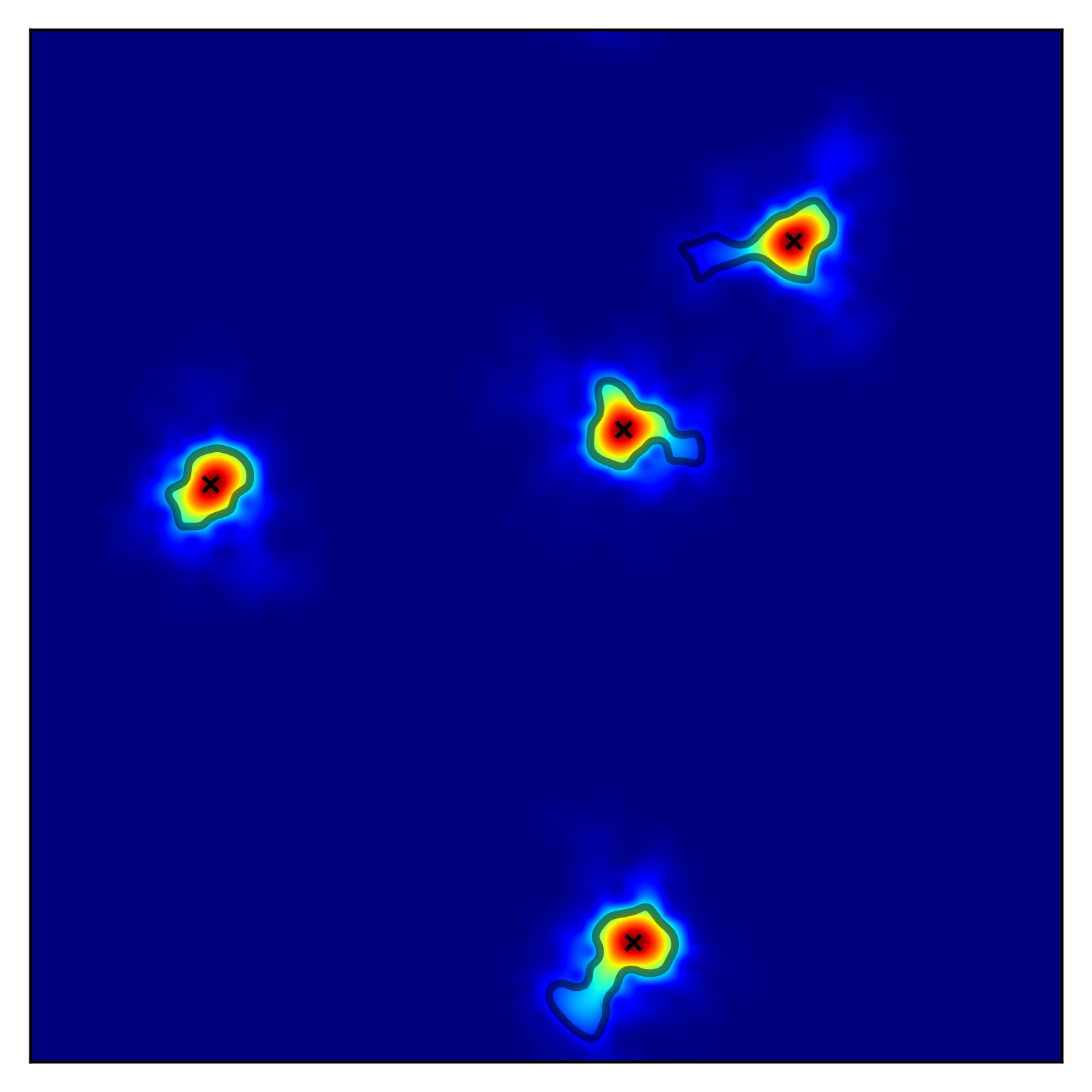}}
 \caption{The first plot shows the prediction for the location of the first four eigenfunctions for the Bernoulli potential of Figure~\ref{fg:bernoulli}. The second plot shows the actual positions of these eigenfunctions, superimposed. The final plot compares the actual positions to the predictions.}
\label{fg:pred2db} 
\end{figure}

We now consider the two-dimensional case where the potential is the random $80\x 80$ Bernoulli potential shown in Figure~\ref{fg:bernoulli}, for which the effective potential is shown in Figure~\ref{fg:berneffpot}. To compute the landscape function we again used the finite element method with Lagrange cubic finite elements on a uniform mesh. The mesh was obtained by dividing each of the unit squares into $10\x 10$ subsquares, each of which was further divided into two triangles, resulting in 1,280,000 triangles altogether. We then evaluated the solution at a uniform grid of $400\x400$ points and found the local minima of the effective potential by comparing each of these values to the values at the eight nearest neighbors (horizontally, vertically, and diagonally). In Figure~\ref{fg:pred2db} the first plot shows the first four local minima of the effective potential. For each, a corresponding sublevel set of the effective potential is shown. The four minima and sublevel sets are our predictors for the locations of the eigenfunctions. The energy level $E$ of the sublevel sets was taken as $1.56$ times the well depth, just slightly larger than the prediction for the eigenvalue, namely $1.5$ times the well depth, which we propose in the next section. Recall that the choice of $E$ is somewhat arbitrary.  This choice gives a good visual match with the apparent support of the eigenfunctions.  The second plot in Figure~\ref{fg:pred2db} is a plot of the sum of four eigenfunctions, each normalized in the $L^{\infty}$ norm. Since they are localized one can easily distinguish the location of each within the sum, which is very close to that predicted. The third plot is a superposition of the first two, to facilitate comparison.

\begin{figure}[htbp]
  \centerline{%
  \includegraphics[height=2.5in]{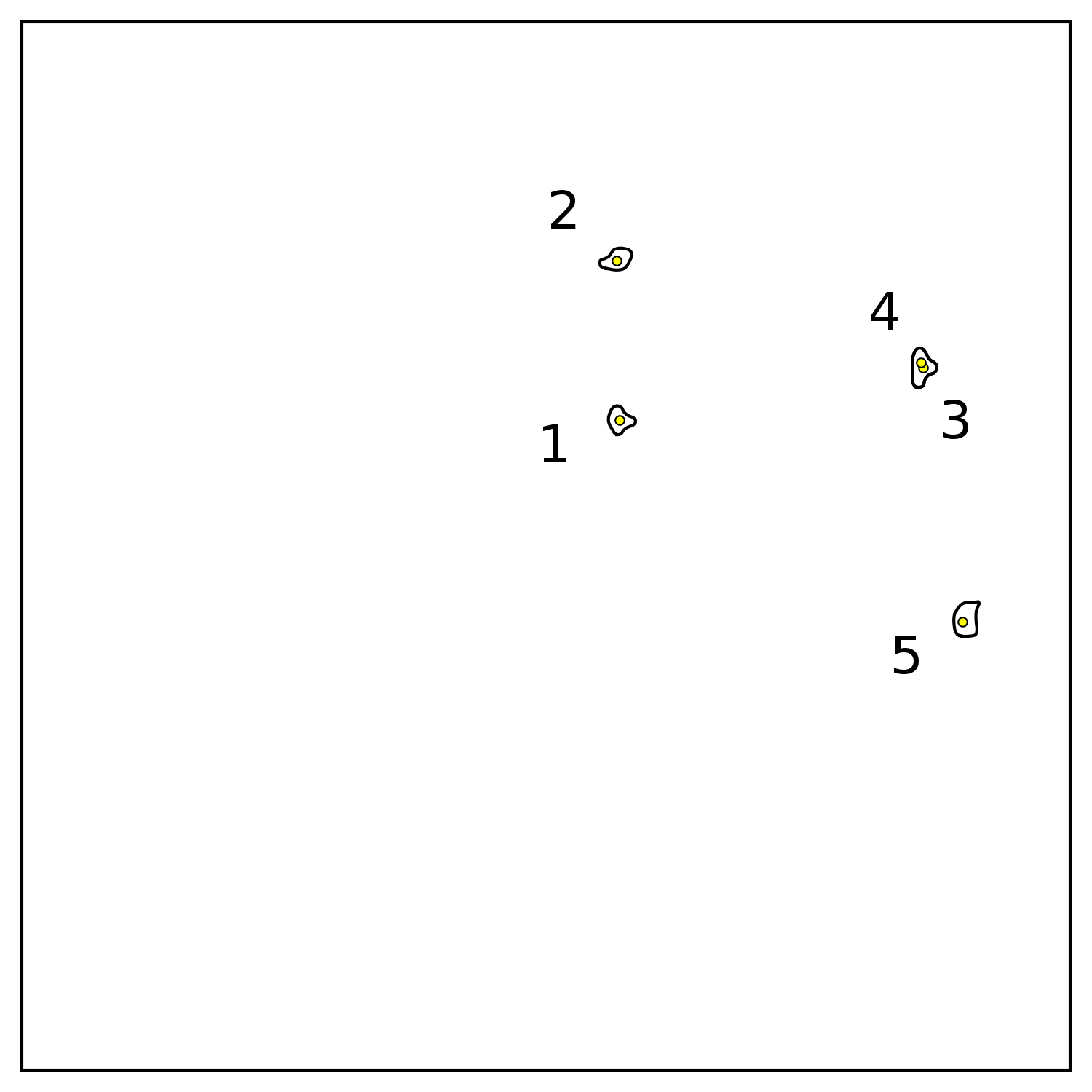}%
  \quad
  \includegraphics[height=2.5in]{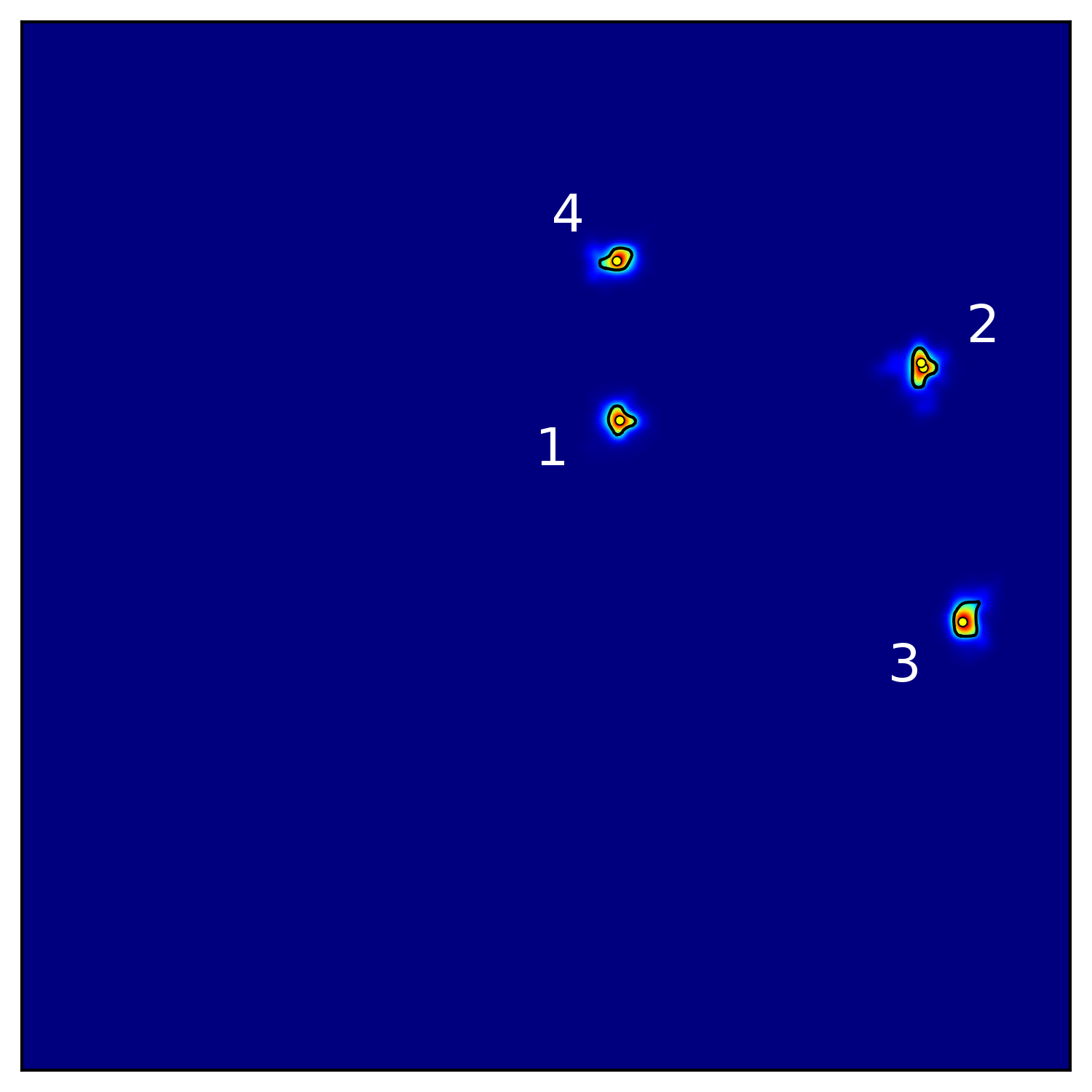}}%
 \caption{For the potential of Figure~\ref{fg:localization}, the correspondence between eigenfunctions and wells of the effective potential does not respect the ordering of the well minima. Moreover, the third and fourth minima effectively define a single well.}
\label{fg:pred2du} 
\end{figure}

In the three cases just considered, there is a clear correspondence between the first four eigenfunctions and the four deepest wells of the effective potential, with each of eigenfunctions, when ordered as usual by increasing eigenvalue, centered  at the corresponding local minimum of the effective potential, ordered by depth of the minimum. However, this ideal situation does not always pertain. When two eigenvalues, or two of the minima, are nearly equal, their ordering may not be respected by the correspondence. Another situation which may arise is that the basin surrounding one of the minima may include another. In that case the second minima does not lead to a separate eigenfunction. Both of these issues arise in the the case of the uniformly random potential of Figure~\ref{fg:localization}. Figure~\ref{fg:pred2du} shows the first five local minima of the effective potential, and their basins. The numbers provided show the ordering by the depth of the wells. Note that the third and fourth minima nearly coincide in location, and they only contribute one well, even though they are, technically, two distinct minima. The actual minimum values and eigenvalues are given in Table~\ref{tb:vals}. It reveals that the third and fourth minima are not only close in location, but nearly coincide in value as well.
Moreover the difference between their value and values of the preceding and following minima is rather small.  These close values account for the fact that the correspondence between the eigenfunctions and minima clearly visible in Figure~\ref{fg:pred2du} does not respect the precise ordering. Nonetheless, the structure of the effective potential clearly provides a lot of information on the localization structure of the eigenfunctions.
To account for such near coincidences we could seek to develop an algorithm to identify clusters of minima with nearly equal values and relate them to clusters of nearly equal eigenvalues. However we shall not pursue this direction here.

\begin{table}[htb]
\centering
\renewcommand{\arraystretch}{1.2}
\begin{tabular}{@{}llllll@{}}
\toprule
& \multicolumn{1}{c}{1} & \multicolumn{1}{c}{2} & \multicolumn{1}{c}{3} & \multicolumn{1}{c}{4} & \multicolumn{1}{c}{5}\\
\cline{2-6}
minima: & 2.3061 & 2.4246 & 2.4763 & 2.4796 & 2.5370 \\
eigenvalue:  & 3.6112 & 3.6618 & 3.7075 & 3.7717 & 3.9190\\
\bottomrule
\end{tabular}
 \caption{The values of the effective potential at its first five minima, and the first five eigenvalues for the uniformly random potential of Figure~\ref{fg:localization}. Cf.~Figure~\ref{fg:pred2du}.}
\label{tb:vals} 
\end{table}

\begin{figure}[htbp]
  \centerline{%
  \includegraphics[height=2in]{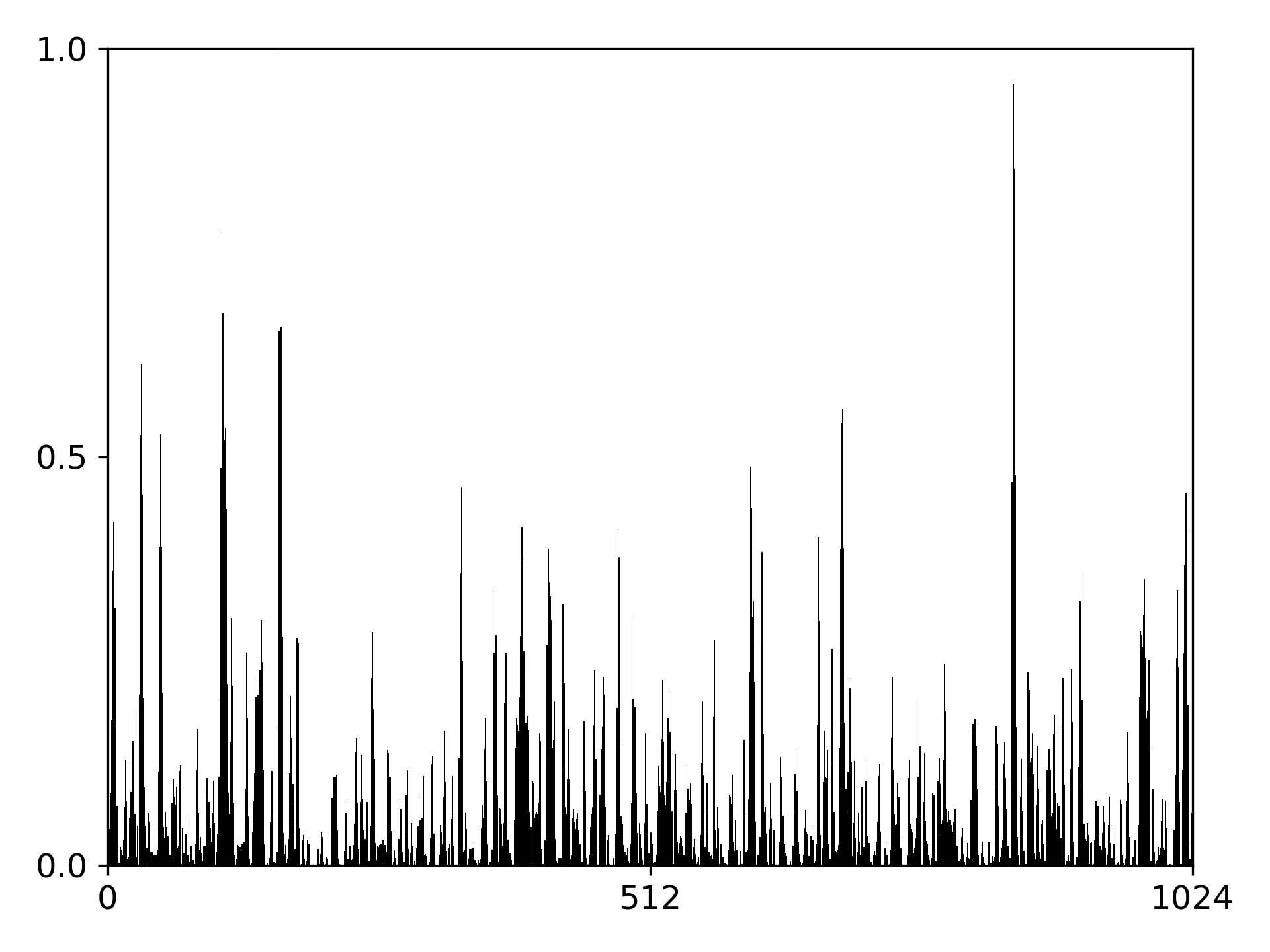}%
  \quad
  \includegraphics[height=2in]{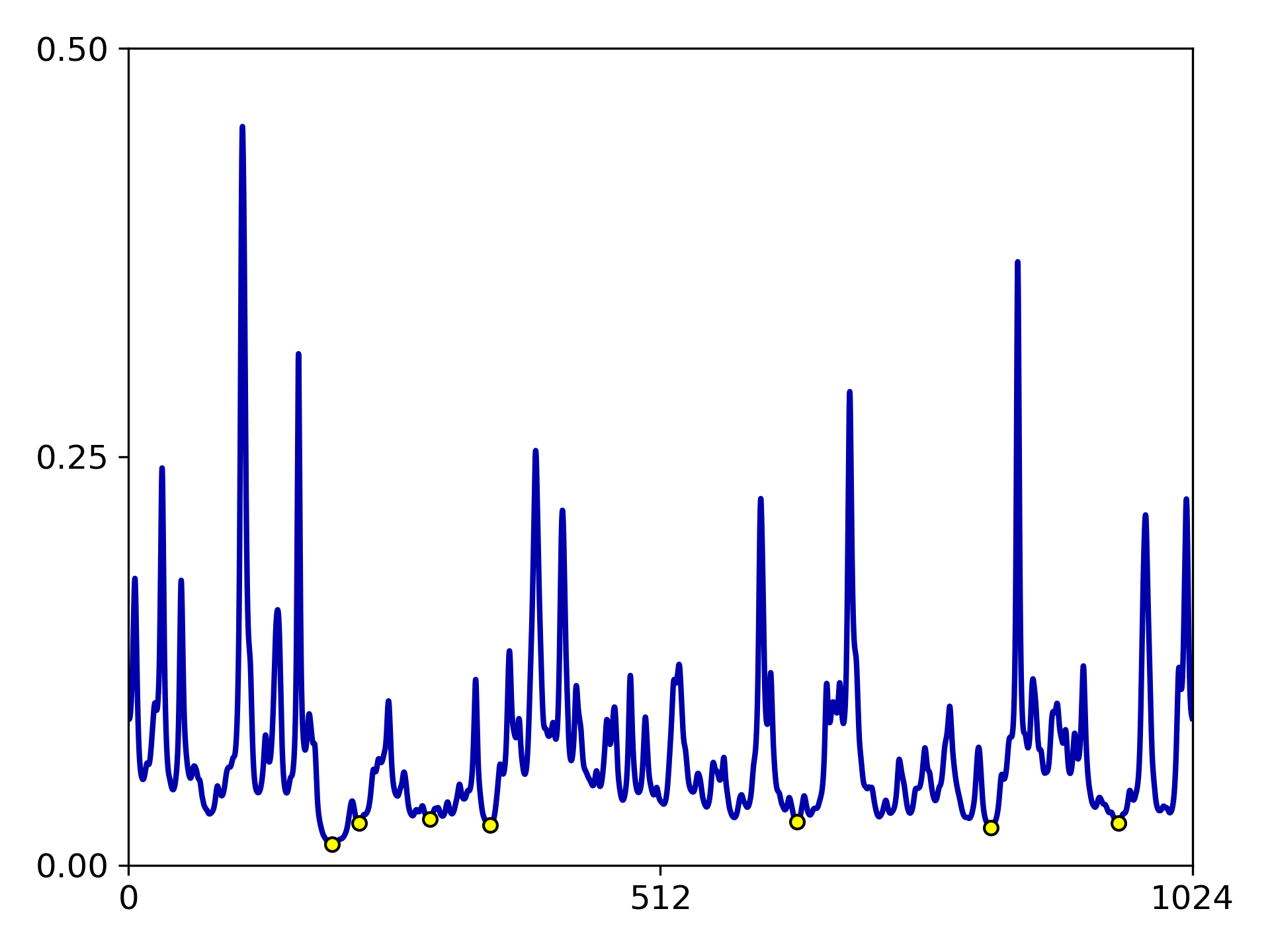}}
  \medskip
  \centerline{%
  \includegraphics[height=2.5in]{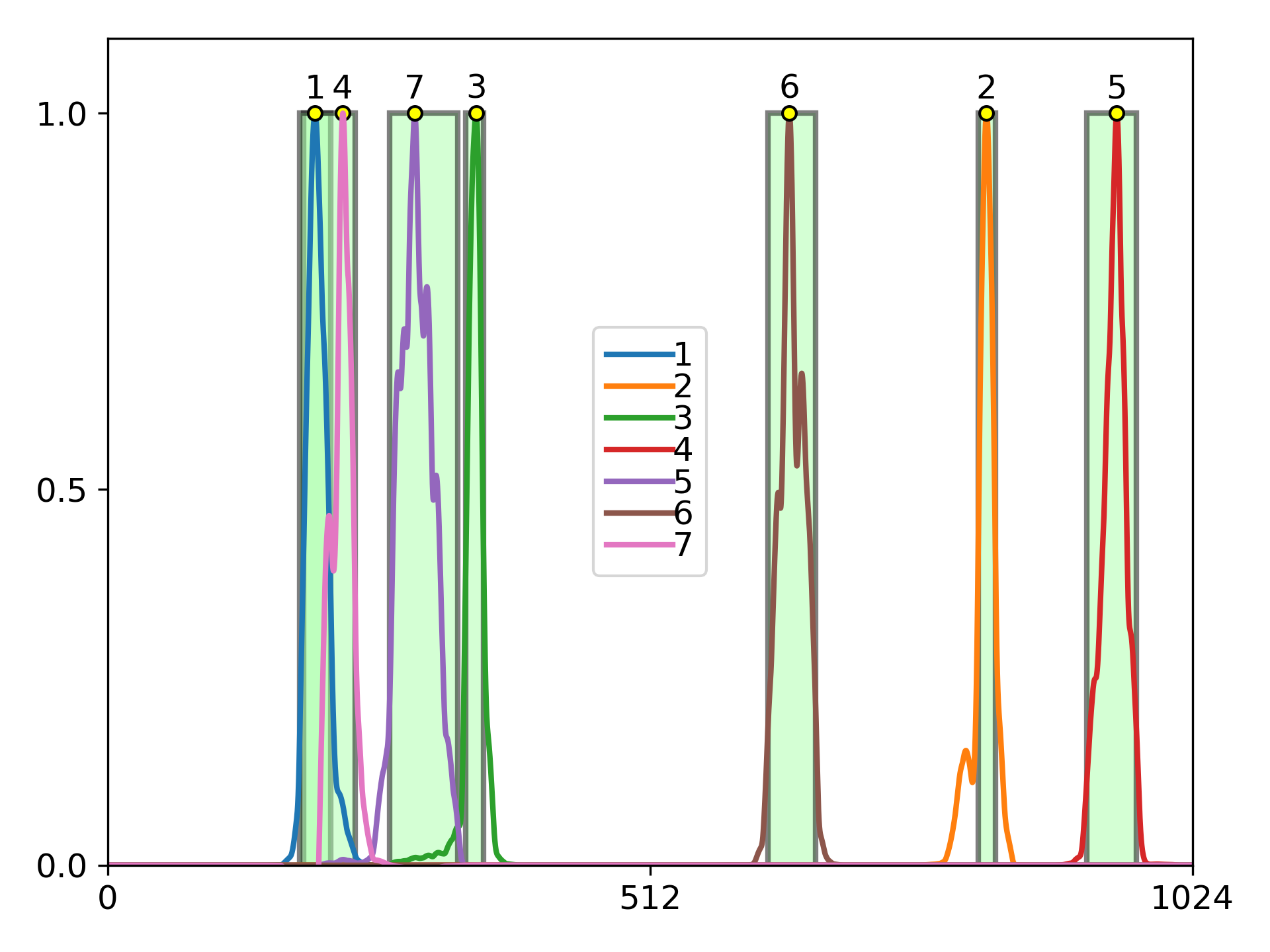}}%
 \caption{A correlated potential and the corresponding effective potential with its seven lowest minima marked. On the top left is the original correlated potential, on the top right the corresponding effective potential. On the bottom are the first seven eigenfunctions with the color giving the order as indicated by the inset. A rectangle is superimposed on each eigenfunction. The heights of the rectangles and eigenfunctions are normalized to unity, while each rectangle's width indicates the extent of the localized eigenfunction as predicted by the effective potential (using a sublevel set at an energy level equal to $1.875$ times the depth of the wells as was done for Figure~\ref{fg:pred1d}). The small yellow circles give the position of the local minima and
 the numbers over them, the order of the local minima.  Notice that the order differs from the order of the corresponding eigenfunctions after the first three.}
\label{fg:corr1d} 
\end{figure}

Thus far we have examined random piecewise constant potentials with  values taken independently and identically distributed according to some probability distribution (uniform or Bernoulli). In the final example of this section, we consider a potential for which the values are \emph{correlated} rather than independent. To generate values for the potential, we use circulant embedding to convert uncorrelated Gaussian $N(0,1)$ random vector samples to correlated Gaussians~\cite{Kroese2015}. We take a 1-dimensional example, in which the potential is piecewise constant with $n$ unit length pieces, with $n$ even. Define $q_i=q_{n-i}=\sigma \exp(-d i)$, $0\le i\le n/2$, where $d$ is a positive constant, and let $Q$ be the diagonal matrix with entries $q_0, q_1, \cdots q_{n-1}$. A sample vector for the values of $V$  is obtained by squaring the components of the vector $F^{-1}Q Fz$ where $z$ is a vector of length $n$ with components sampled independently from a normalized Gaussian distribution, and $F$ is the discrete Fourier transform on $\mathbb C^n$. This type of a random potential is typically created by optical speckles in a Bose-Einstein condensate. See, e.g., \cite{Modugno}, \cite{Falco}. It is quite challenging to derive rigorous probabilistic results when the correlation is not negligibly small, especially in higher dimensions. The landscape theory, however, continues to apply. We consider an example 
with $n=1,024$, $\sigma=1.0$, $d=0.01$,
shown in Figure~\ref{fg:corr1d} along with the corresponding effective potential. Note that, although the potential and effective potential look quite different from the previous examples, the effective potential still has clearly defined wells which allow us to apply our theory. In the final plot in Figure~\ref{fg:corr1d} we use the effective potential as before to predict the location and extent of the first seven eigenfunctions. As in the uncorrelated cases, the well minima very nearly coincide with the peak of the eigenfunctions. The first three minima, in order, correspond to the first three eigenfunctions, but after that the order is not the same. The discrepancy in ordering is not very significant, however, since both the minimum values and the eigenvalues are very close to one another, as indicated in Table~\ref{tb:corr1d}.

\begin{table}[htb]
\centering
\renewcommand{\arraystretch}{1.2}
\begin{tabular}{@{}llllllll@{}}
\toprule
& \multicolumn{1}{c}{1} & \multicolumn{1}{c}{2} & \multicolumn{1}{c}{3} & \multicolumn{1}{c}{4} & \multicolumn{1}{c}{5} & \multicolumn{1}{c}{6} & \multicolumn{1}{c}{7}\\
\cline{2-8}
min:  & 0.012752 & 0.022966 & 0.024642 & 0.025647 & 0.025673 & 0.026626 & 0.028347\\
eig: & 0.016534 & 0.030069 & 0.031620 & 0.032953 & 0.033822 & 0.034479 & 0.034735\\
\bottomrule
\end{tabular}
 \caption{The values of the effective potential at its first seven minima, and the first seven eigenvalues for the correlated potential of Figure~\ref{fg:corr1d}.}
\label{tb:corr1d}
\end{table}

In this section we have shown how we can deduce the approximate locations and
approximate supports of eigenfunctions just by processing the effective potential, without
solving eigenvalue problems.  We remark that this information could be refined
to give an approximation of the precise shape of the eigenfunction.  To do so,
one could solve for the eigenfunction with a standard PDE eigensolver, but with
the domain taken as a regular domain just slightly larger than
the approximate support, and with Dirichlet boundary conditions. Because
of the localization, this computational domain will be much smaller than the original
domain and this computation much less expensive than a global eigenvalue solve.
The development and study of such algorithms, however, goes beyond the scope of
this paper, and is left for future work.

\section{Eigenvalue prediction}\label{sec:eval}

We now turn to the question of predicting eigenvalues of the Schr\"odinger operator $H$ from the effective potential. As a simple illustration of the utility of the effective potential for eigenvalue estimation, we start by recalling the basic lower bound on the fundamental eigenvalue in terms of the potential $V$, and show how it can be improved by using the effective potential.

\begin{figure}[htbp]
   \centerline{\includegraphics[width=2.5in]{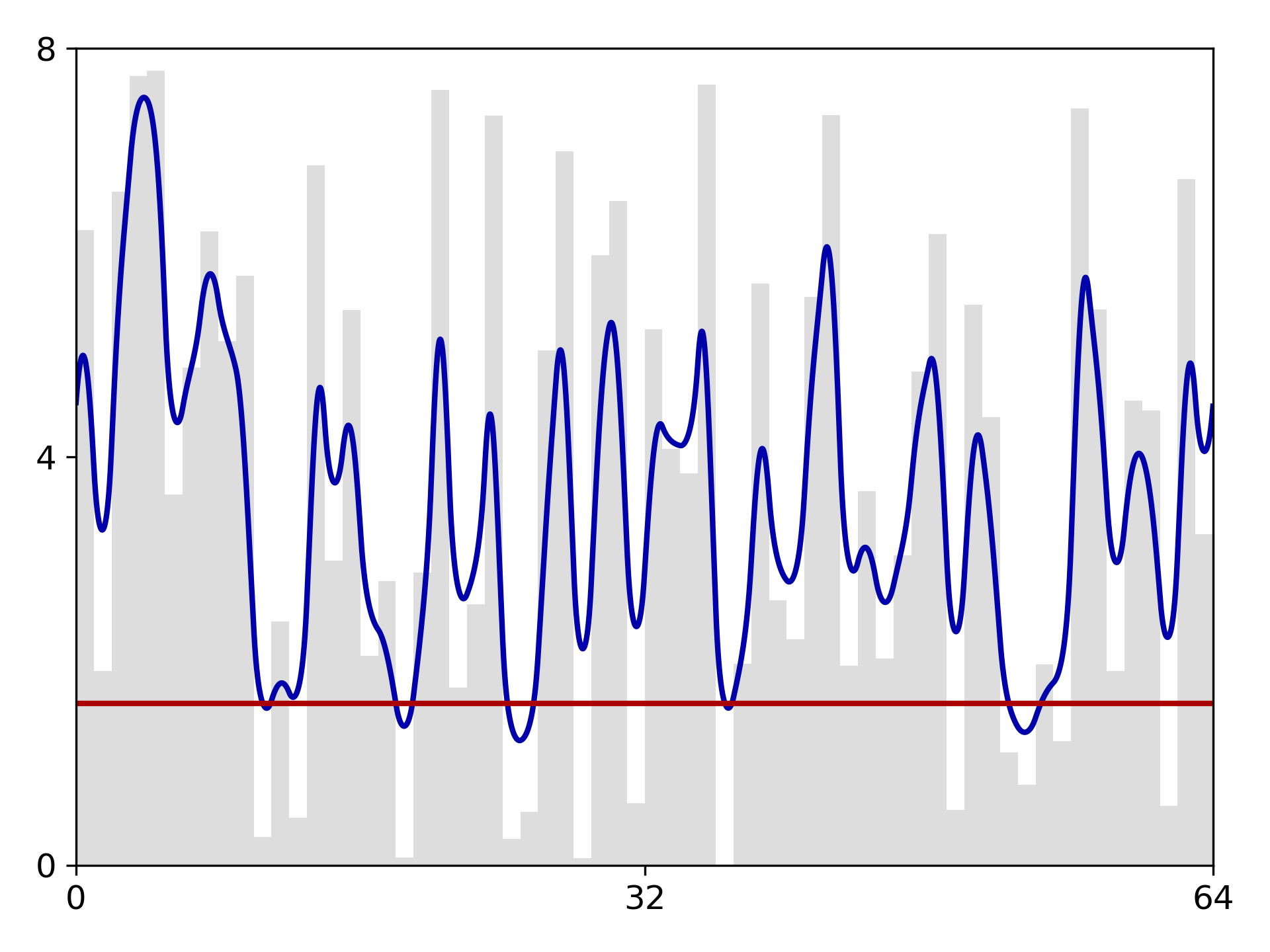}}
 \caption{A random potential (in shaded gray), the corresponding effective potential (the solid blue line), and the fundamental eigenvalue (the horizontal red line).}
\label{fg:comp} 
\end{figure}

For an eigenfunction $\psi$ of $H$ with eigenvalue $\lambda$, normalized to have $L^2$ norm $1$, we have
\begin{equation}\label{ekp}
\lambda = (H\psi,\psi) = \|\grad \psi\|^2 + (V\psi,\psi)
\end{equation}
which represents the decomposition into kinetic and potential energy. Dropping the kinetic energy term and replacing $V$ by its infimum gives a lower bound on the eigenvalues:
\begin{equation}\label{vbd}
\lambda \ge \inf V.
\end{equation}
Now we use the fundamental identity of Proposition~\ref{th:ep} to decompose the eigenvalue in terms of the effective potential:
$$
(H\psi,\psi) = (u^2 L\phi,\phi) + (W\psi,\psi),
$$
where $\phi=\psi/u$. In view of the form of $L$ \eqref{defL}, the first term on the right hand side is positive, so dropping it and replacing $W$ by its infimum gives another lower bound:
\begin{equation}\label{wbd}
\lambda \ge \inf W.
\end{equation}
Figure~\ref{fg:comp} allows one to
compare the two bounds for a random potential with $64$ values chosen uniformly iid in the range $[0,8]$.  The fundamental eigenvalue for this realization is $1.58$, indicated on the plot in red. The infimum of $V$ is, however, very near zero: $0.00009$, and so the bound \eqref{vbd} is nearly worthless.  (In this realization $\inf V$ happens to be particularly small, but the expected value of $1/65=0.015$ is again of little use.) By contrast, $\inf W = 1.22$, which is a useful lower bound.  In fact, the fundamental eigenvalue is equal to about $1.3~\inf W$. We shall see below that this factor of roughly $1.25$ or $1.3$ applies for a wide range of random potentials in one dimension.

In the remainder of this section, we shall consider two approaches to eigenvalue prediction, one based solely on the local minima values of the effective potential, and the other based on a variant of Weyl's law utilizing the effective potential. Although we shall explain the thinking behind these approaches, it has to be noted that neither has yet been justified rigorously.

\subsection{Eigenvalues from minima of the effective potential}\label{ssec:eval1}

In this section we discuss the approximation of the eigenvalues of the Schr\"odinger operator $H$ using the effective potential $W$. We shall be interested in both the approximation of individual eigenvalues, and in the distribution of the eigenvalues.  The latter is captured by the \emph{density of states} (DOS). Defined precisely, the DOS is the distribution on $\R$ obtained by summing the delta functions centered at each eigenvalue. For visualization, it is often converted to a piecewise constant function with respect to a partition of the real line into intervals of some length $\epsilon>0$, with the value over any interval being the integral of the DOS over the interval (so a plot of this function is a histogram of the eigenvalues using bins of width $\epsilon$).  The integrated density of states (IDOS) is the integral of the DOS from $-\infty$ to a real value $E$.  The resulting function $N:\R\to \mathbb N$ is simply the eigenvalue counting function, with $N(E)$ defined as the number of eigenvalues $\le E$.

\begin{figure}[htbp]
  \centerline{%
  \includegraphics[height=1.25in]{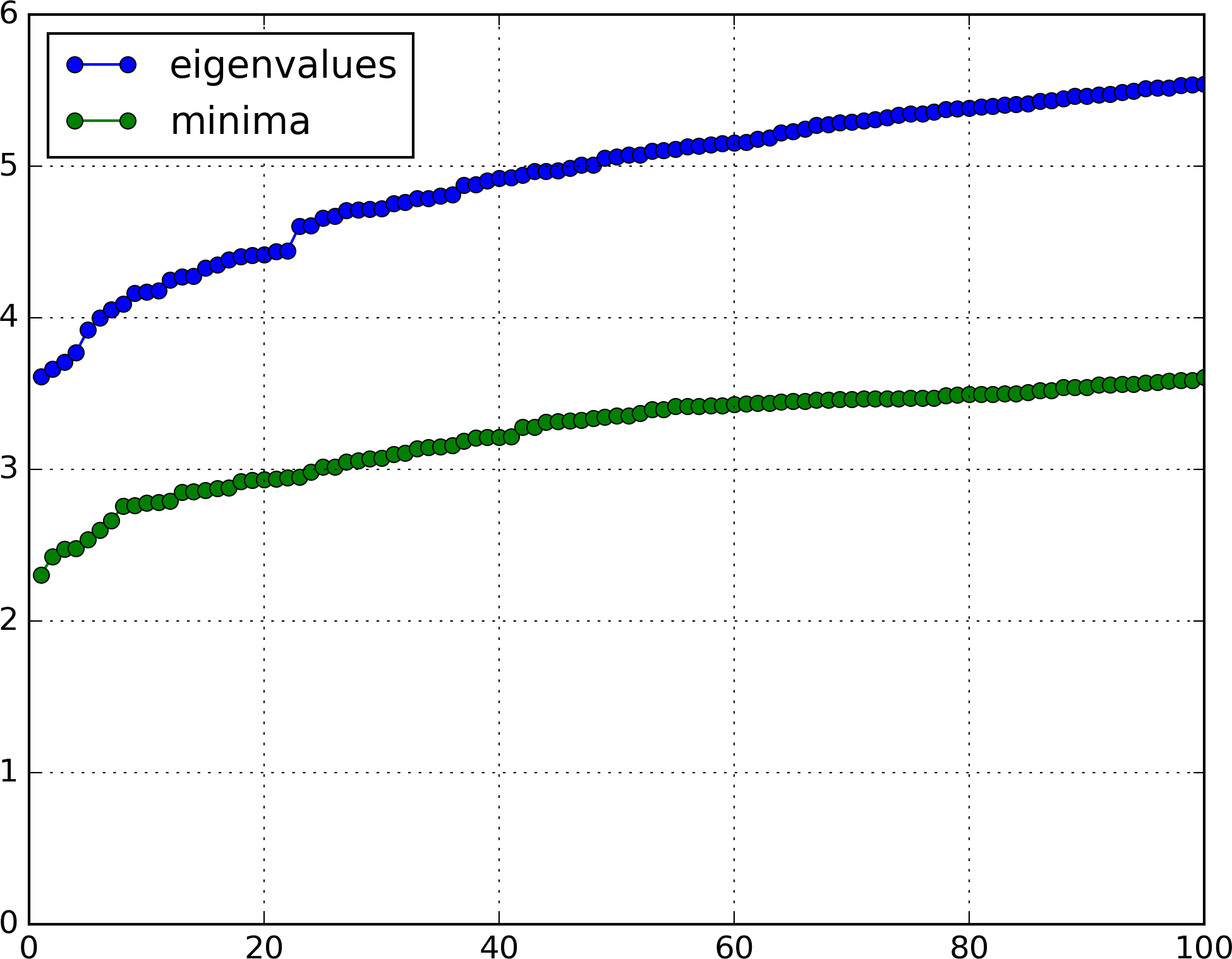}
  \includegraphics[height=1.25in]{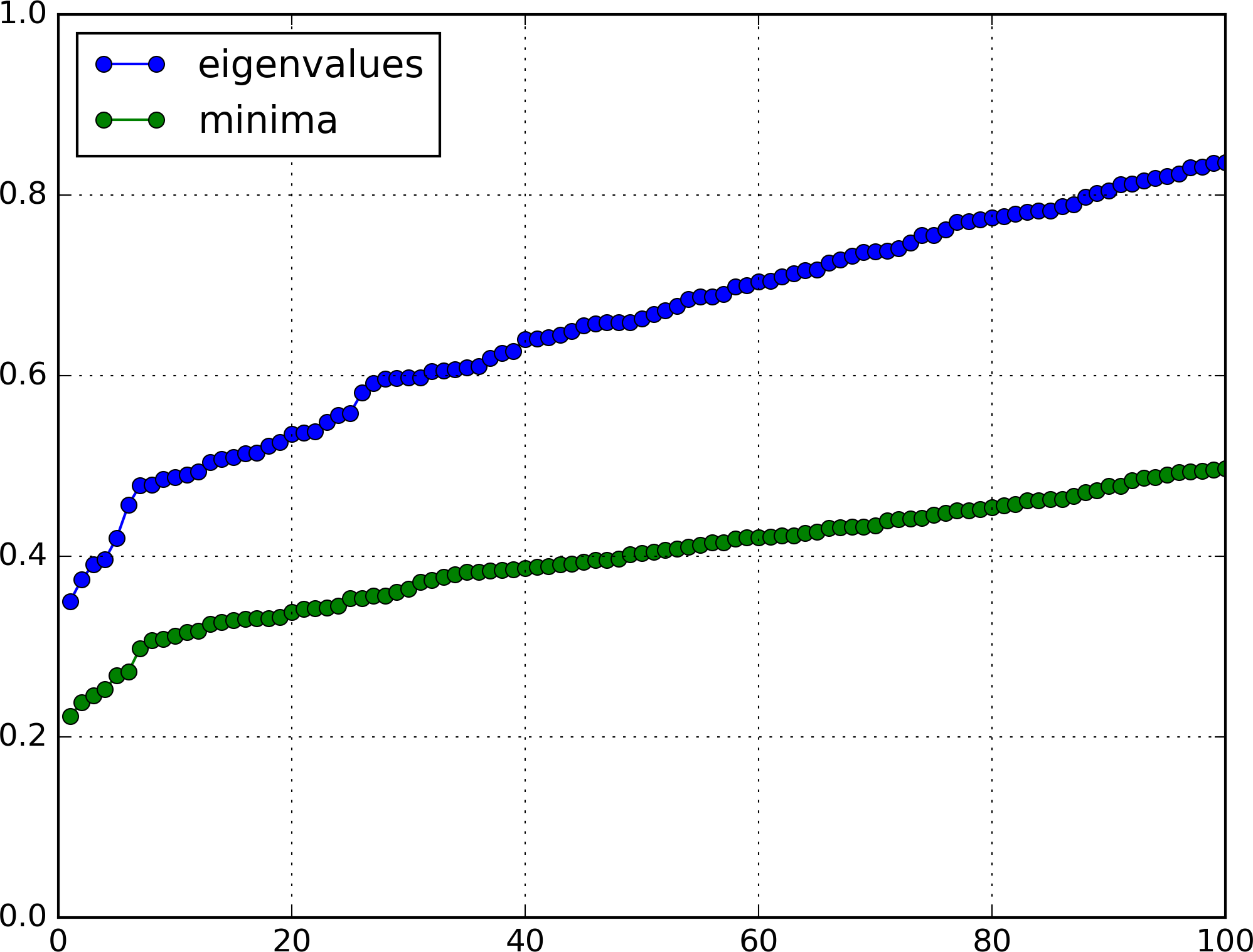}
  \includegraphics[height=1.25in]{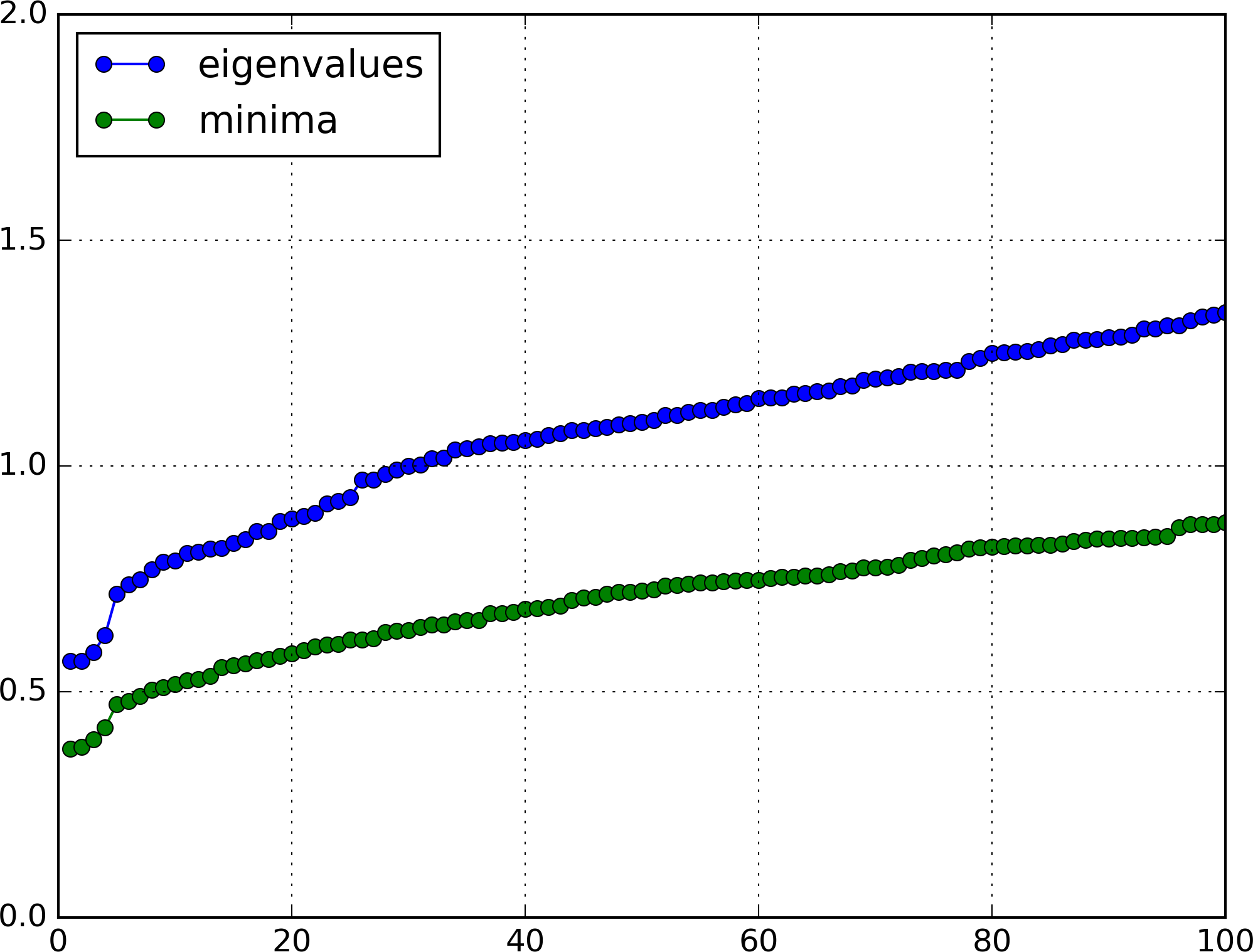}
  }
  \smallskip
  \centerline{%
  \includegraphics[height=1.25in]{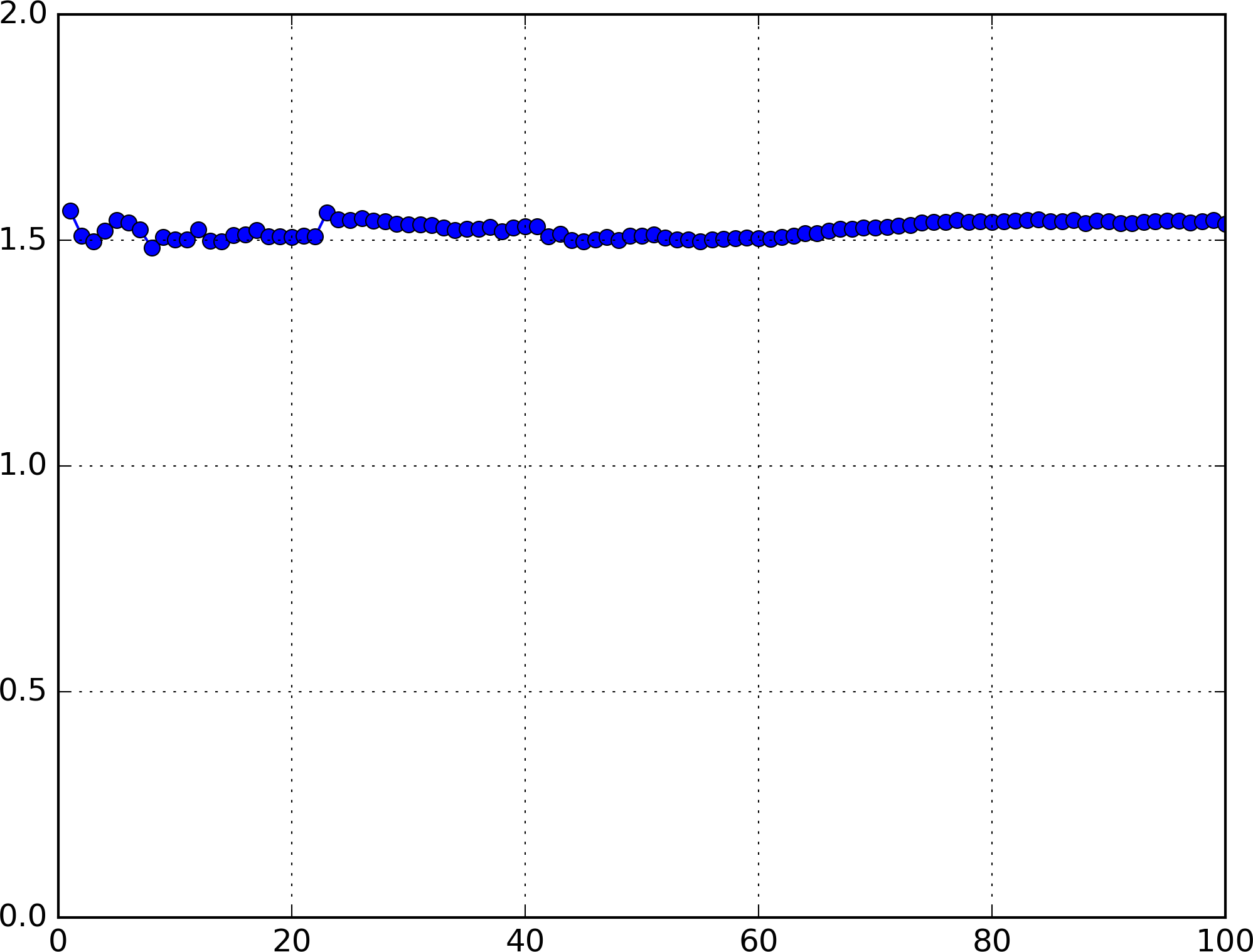}
  \includegraphics[height=1.25in]{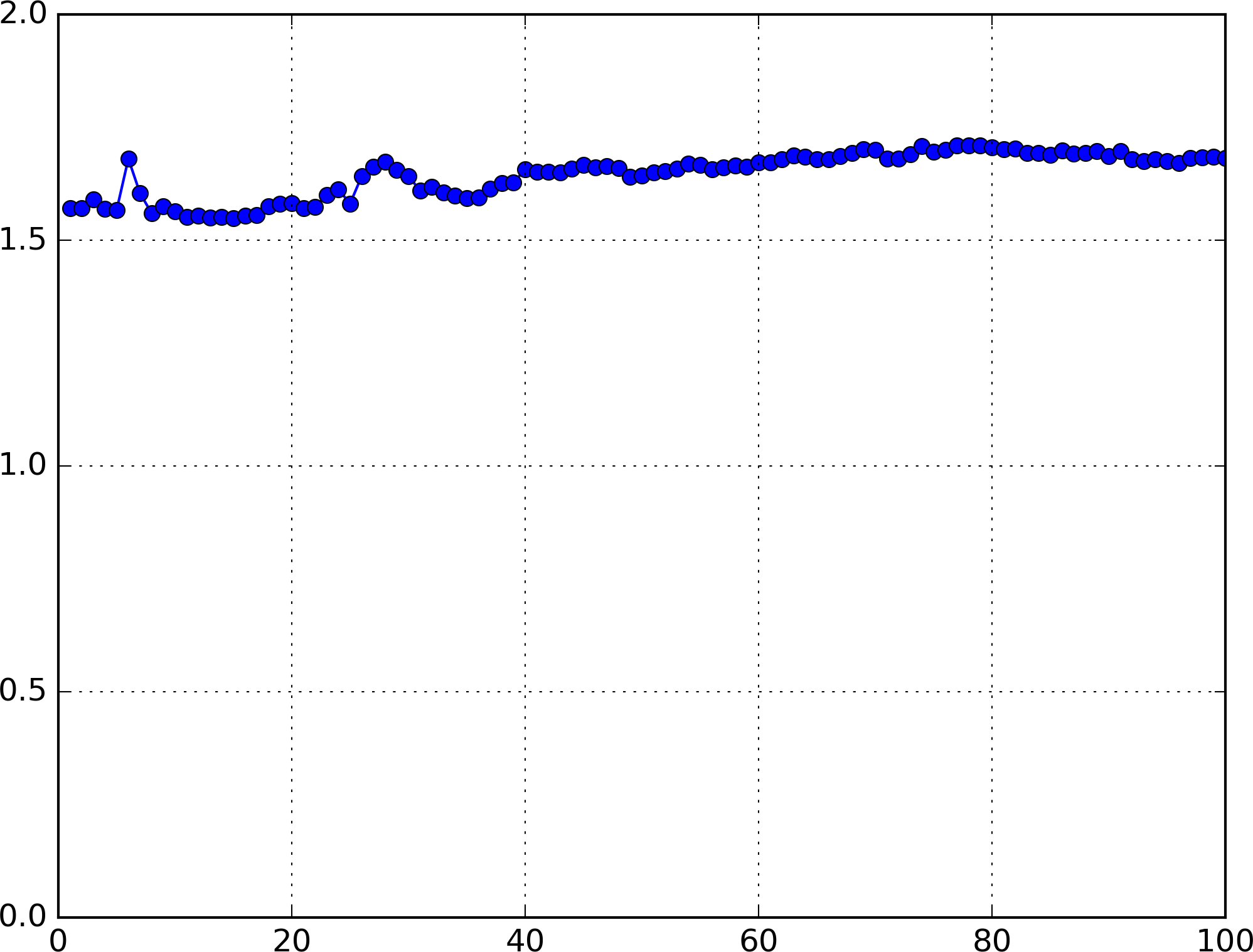}
  \includegraphics[height=1.25in]{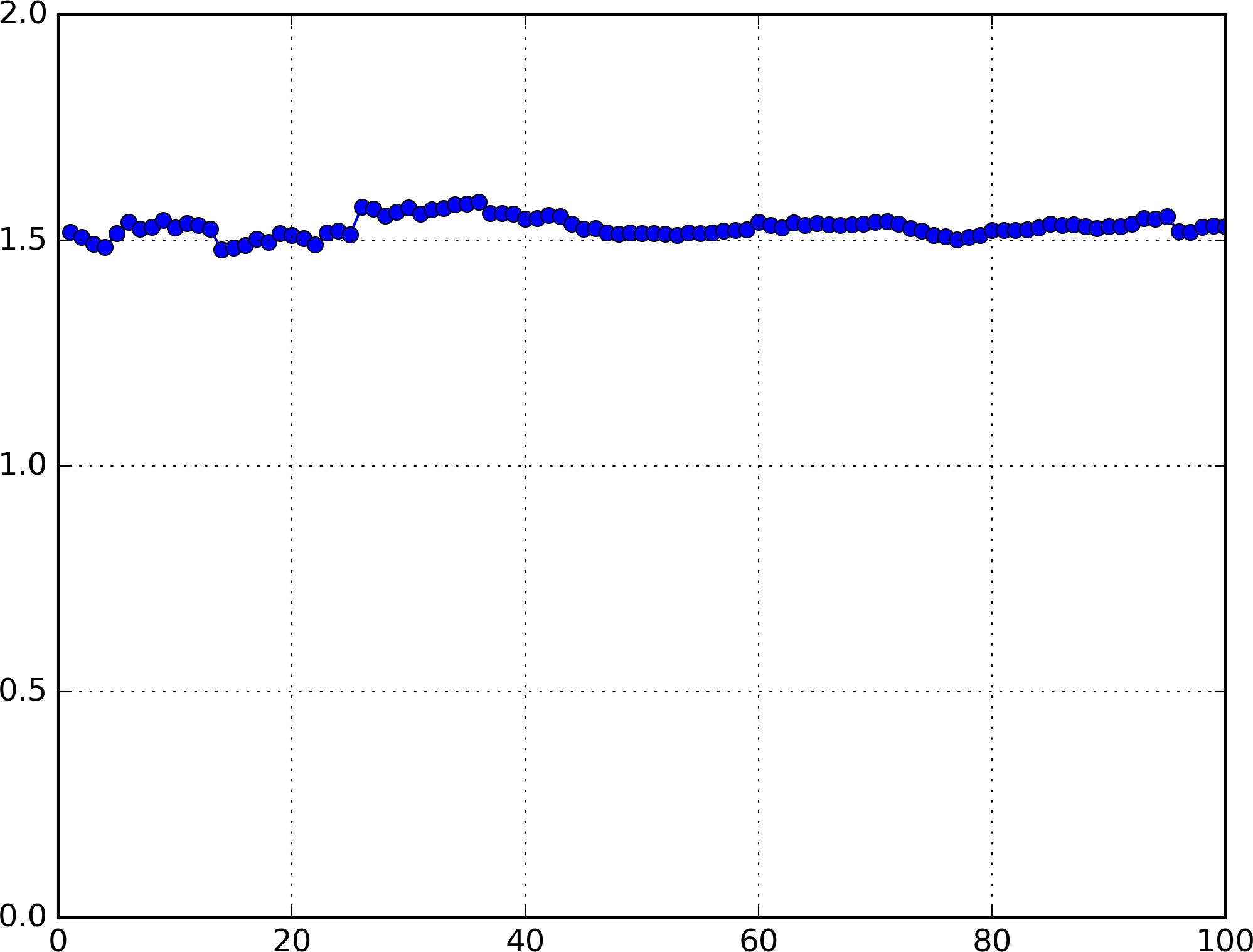}
  }
  \bigskip
  \centerline{%
  \includegraphics[height=1.25in]{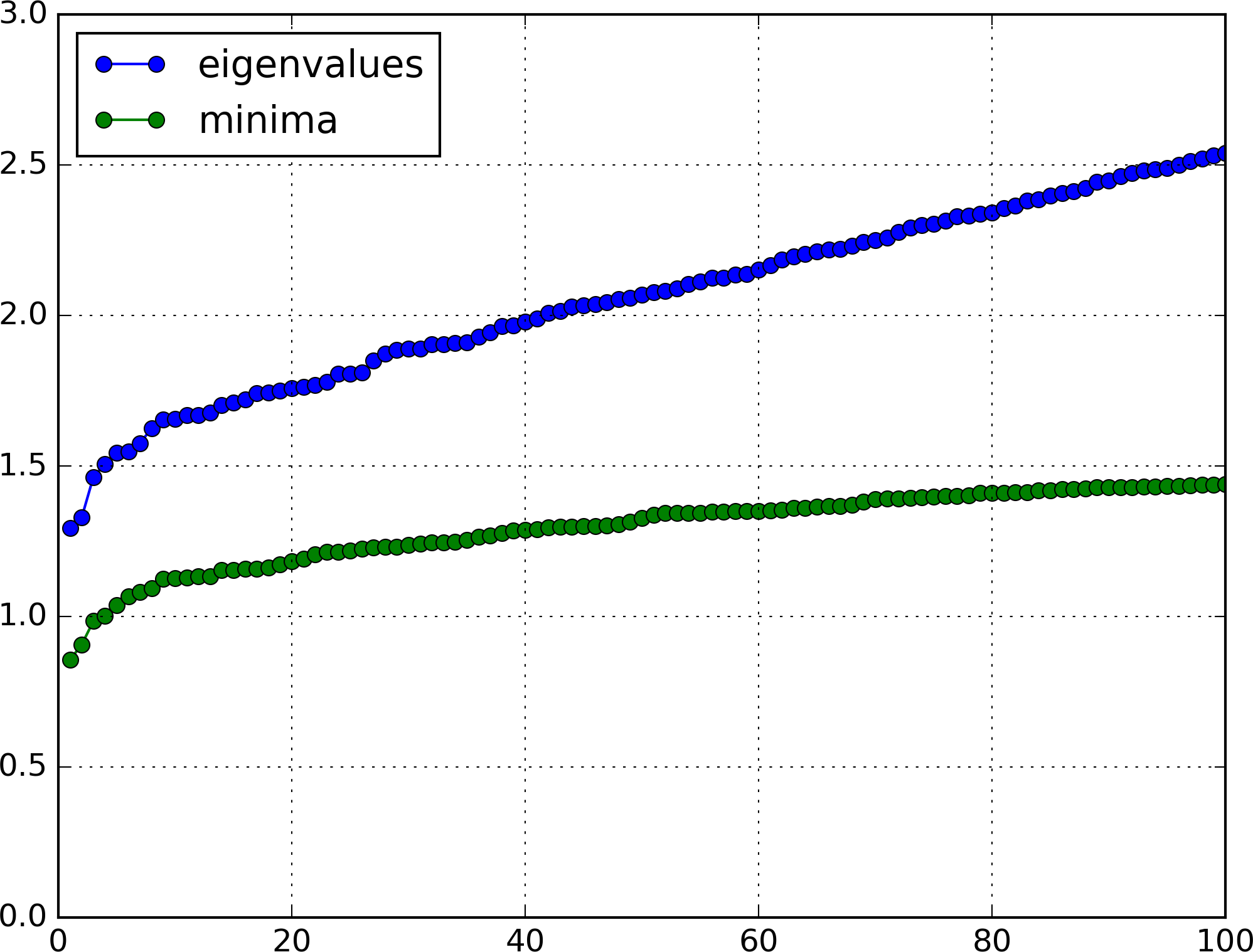}
  \includegraphics[height=1.25in]{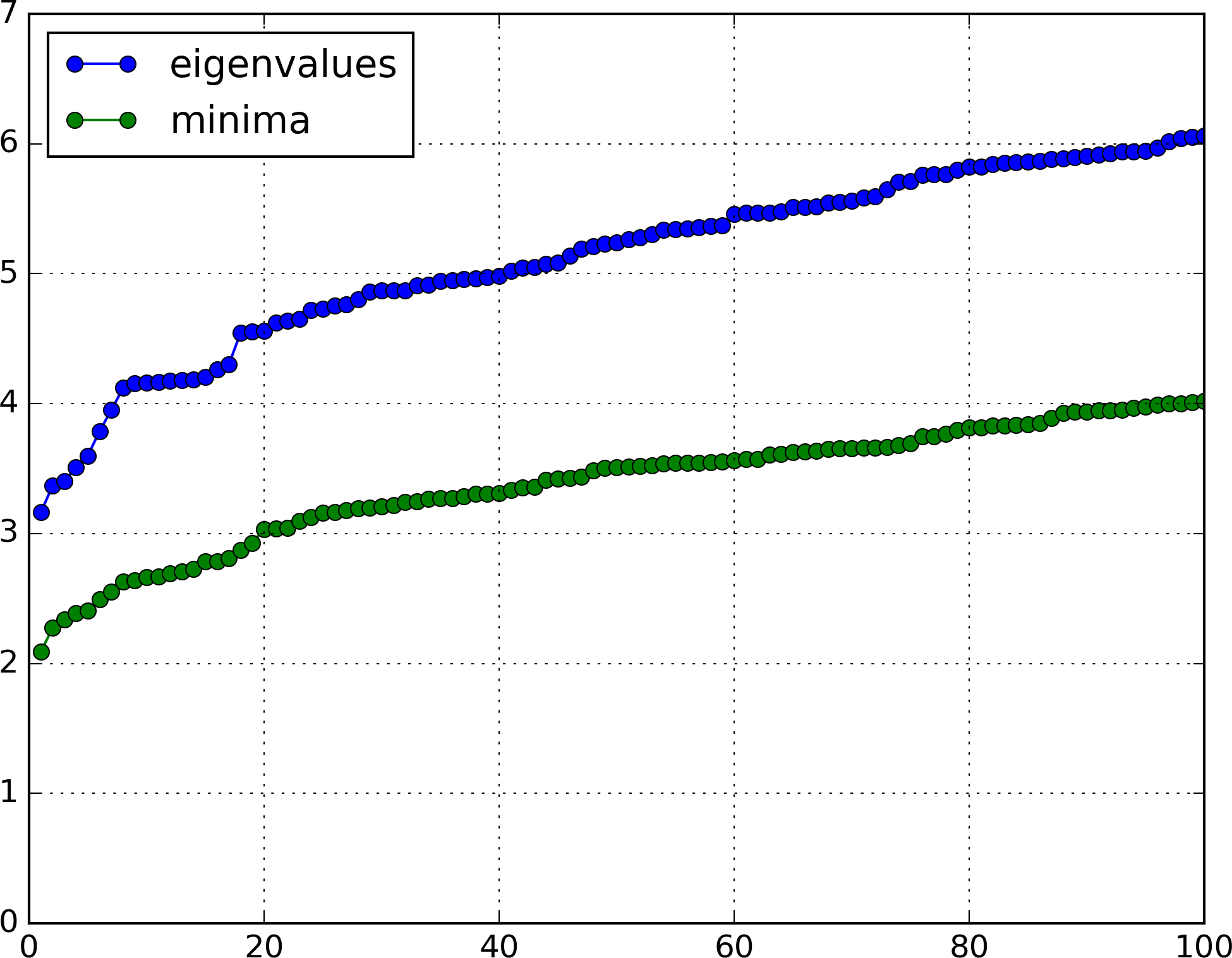}
  \includegraphics[height=1.25in]{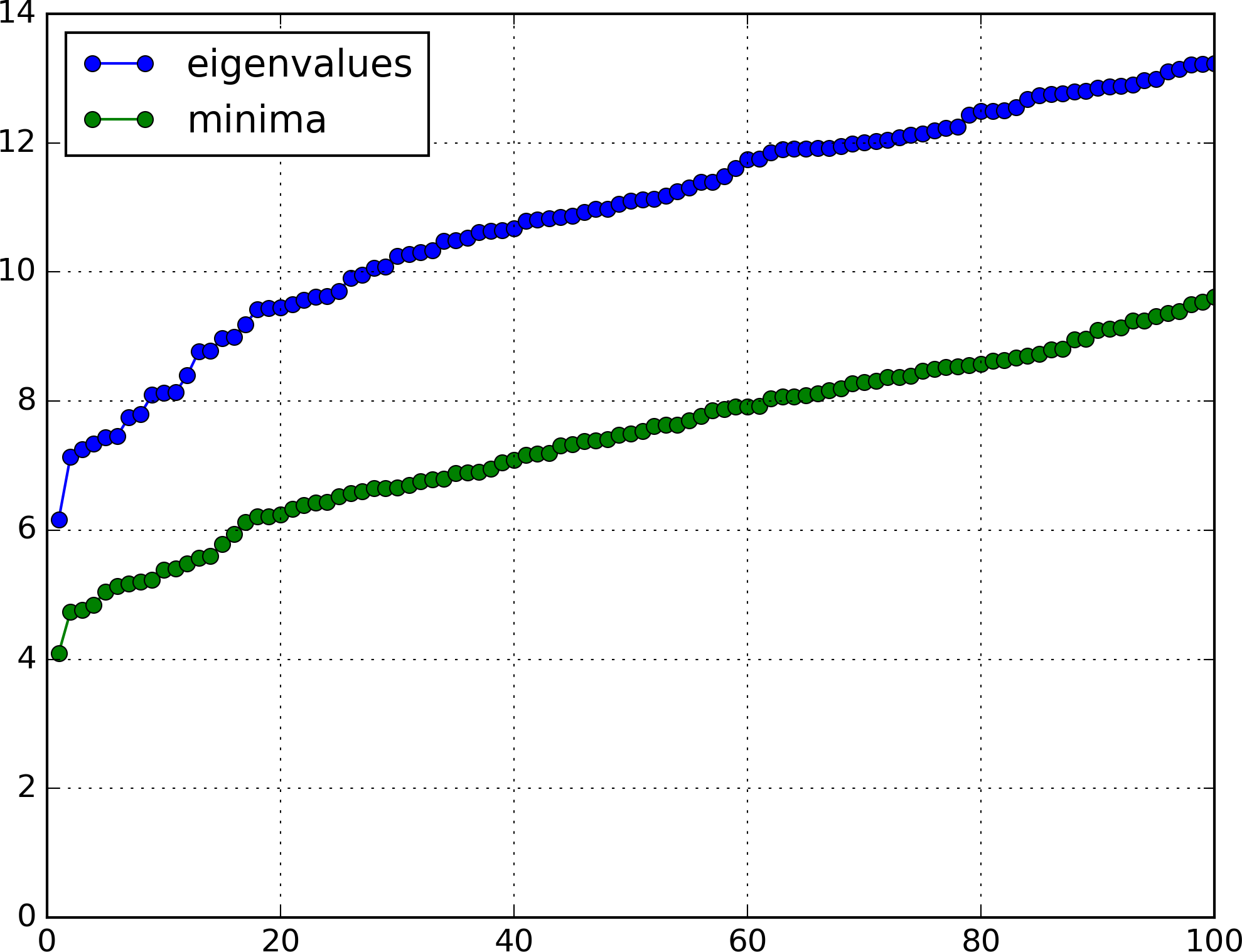}
  }
  \smallskip
  \centerline{%
  \includegraphics[height=1.25in]{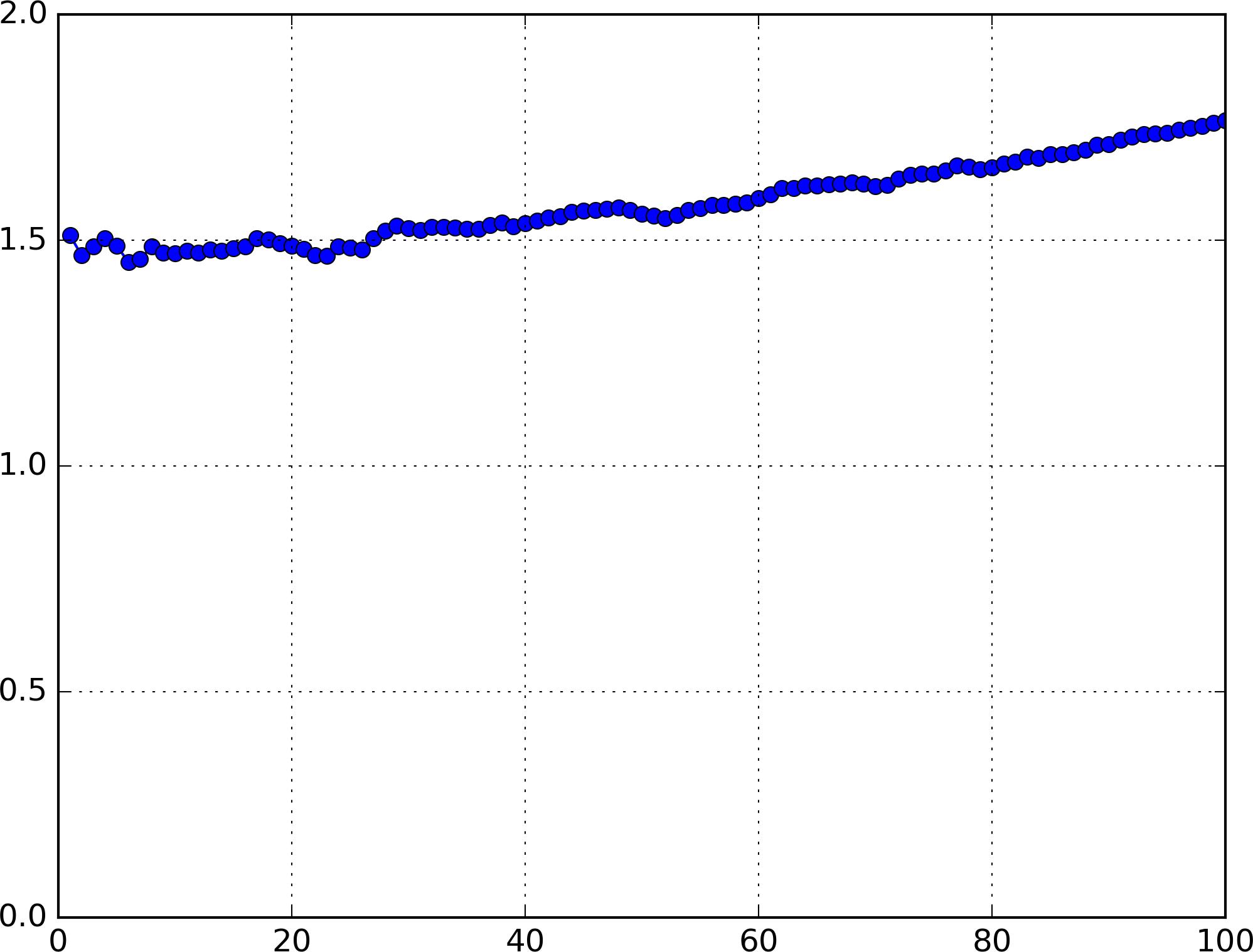}
  \includegraphics[height=1.25in]{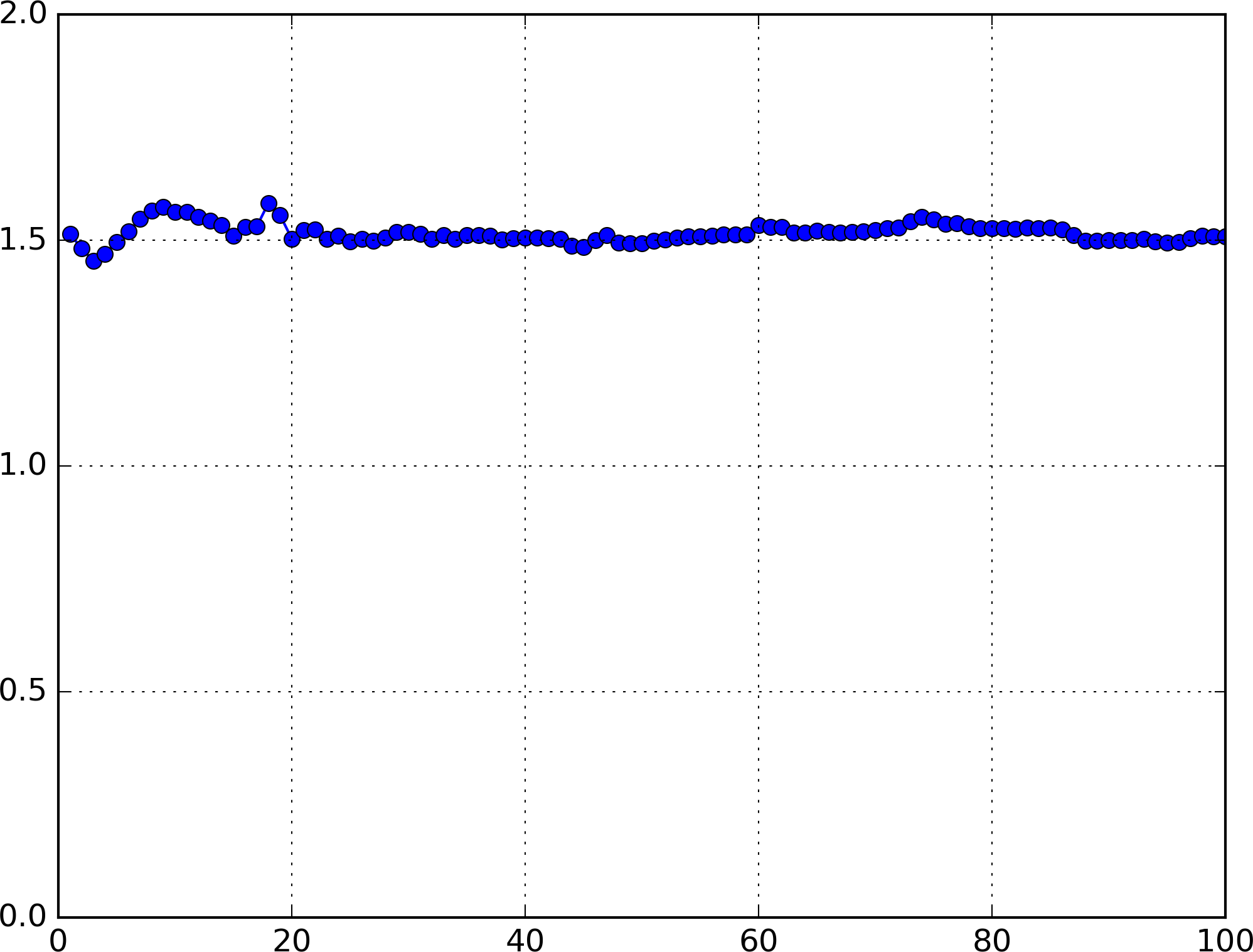}
  \includegraphics[height=1.25in]{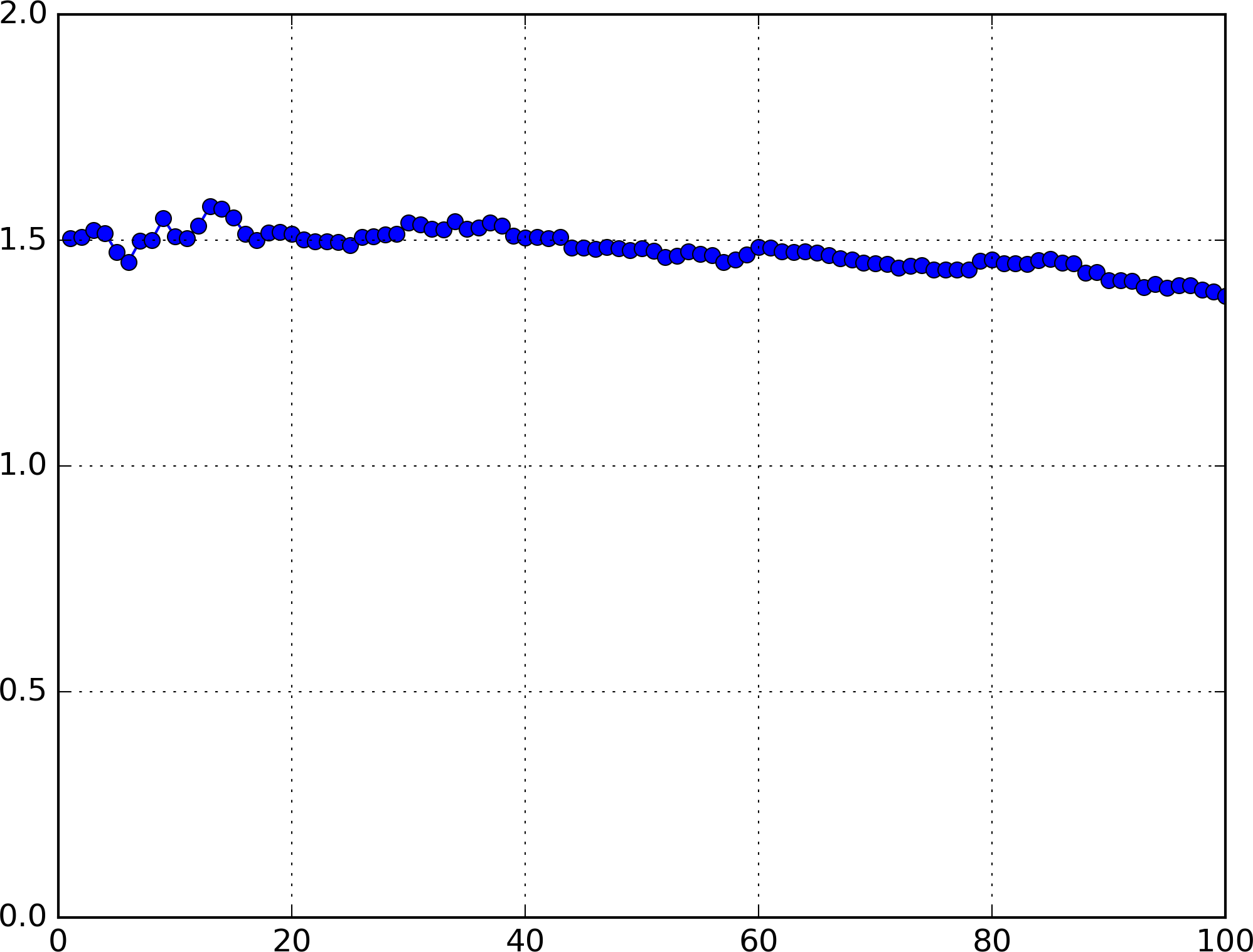}
  }
 \caption{A comparison of $100$ eigenvalues with the corresponding minimum values of $W$ for six different potentials with the ratio of the two shown beneath each one. Top two rows: (left) random $80\x 80$ piecewise potential with values chosen uniformly in $[0, 20]$; (center) Bernoulli potential of Figure~\ref{fg:bernoulli}; (right) correlated potential of Figure~\ref{fg:correlated}. Bottom two rows: random $40\x 40$ piecewise potential with values chosen uniformly in (left) $[0, 4]$, (center) $[0, 16]$, (right) $[0, 64]$. In all these cases, one can notice that the ratio of the minimum values of $W$ to the corresponding eigenvalues remains remarkably close the value 1.5.}
\label{fg:ratio} 
\end{figure}

\begin{figure}[htbp]
   \centerline{\includegraphics[width=58mm]{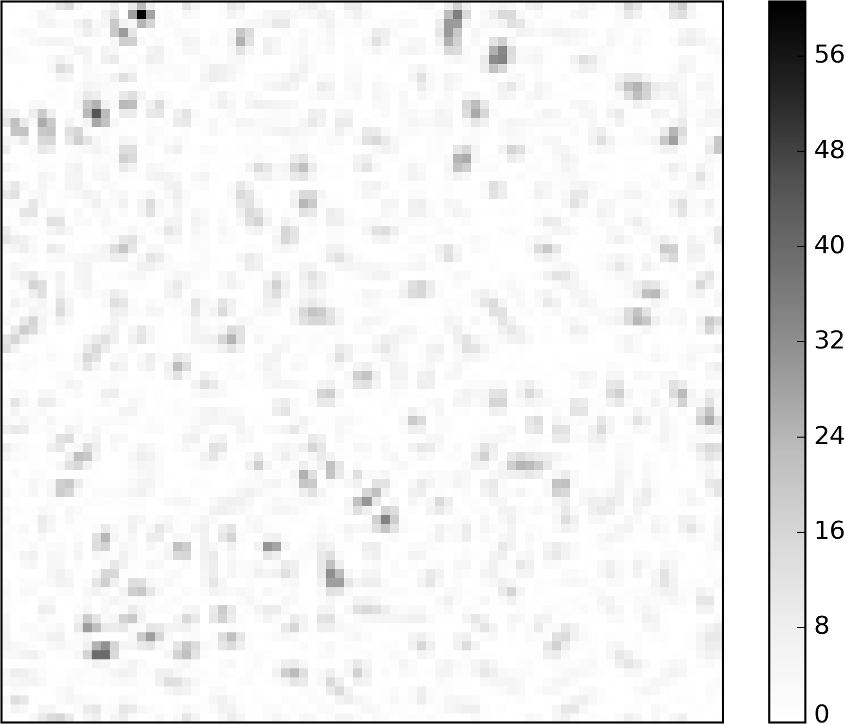}\quad
   \includegraphics[width=65mm]{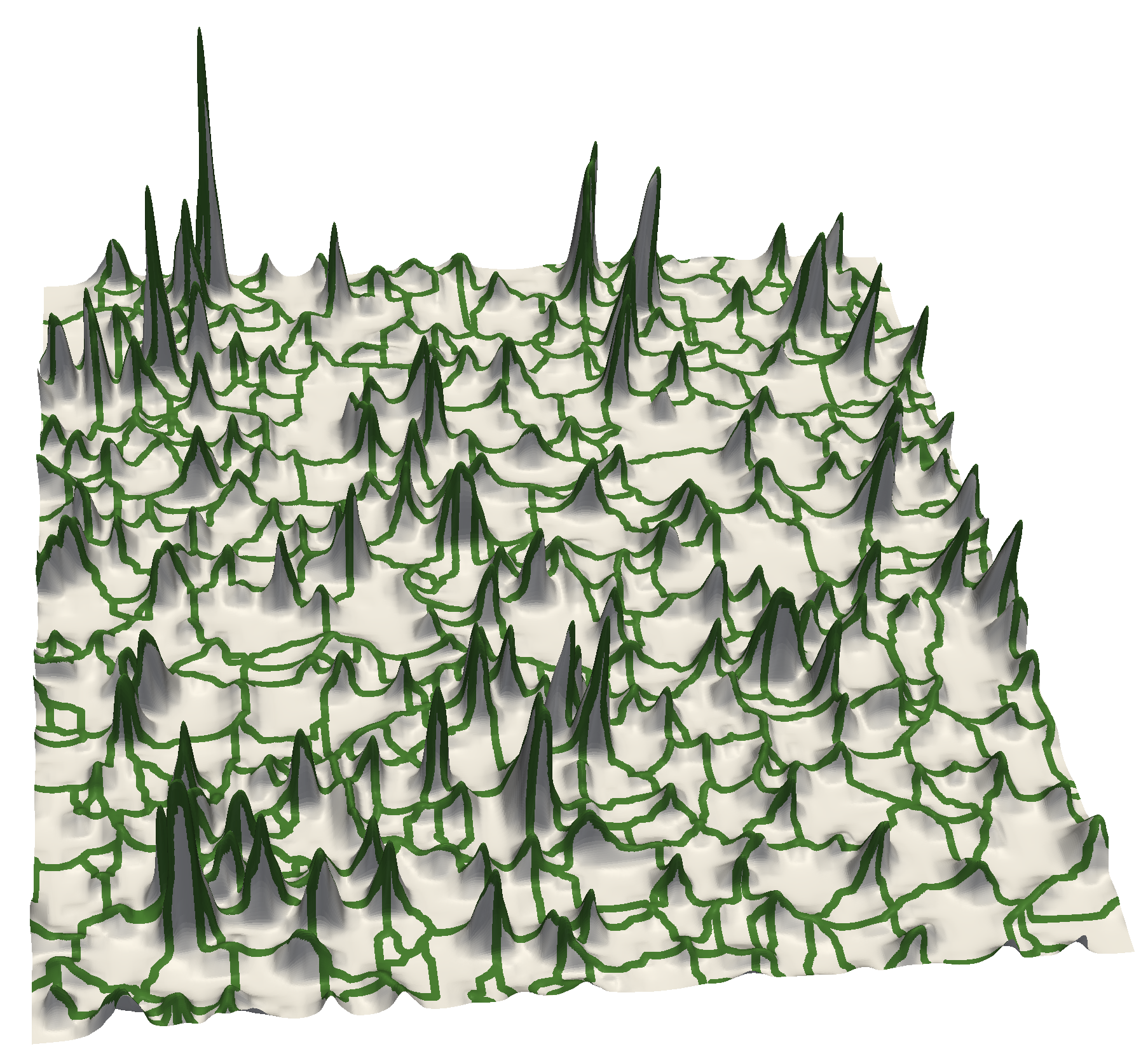}}%
 \caption{A correlated potential with $80\x 80$ pieces and the corresponding effective potential.  The potential takes values ranging from $2.364\x 10^{-9}$ to $60.56$, while the effective potential's values range from $0.374$ to $57.17$.}
\label{fg:correlated} 
\end{figure}

To motivate our first approach to eigenvalue prediction, consider again the potential shown in Figure~\ref{fg:localization}, a piecewise constant function with $80\x 80$ pieces and constant values chosen randomly and independently from $[0, 20]$. In the first plot of Figure~\ref{fg:ratio} we have computed the values of the effective potential~$W$ at its local minima and compared the first 100 of these, in increasing order, to the first 100 eigenvalues of $W$. Just below we plot the quotients of each of these eigenvalues divided by the corresponding local minimum value of $W$. Observe that the quotient is quite constant, taking on a value of roughly $1.5$. We shall endeavor to explain this value below, but first we observe that this ratio of roughly $1.5$ between the $m$th eigenvalue and the $m$th minimum value of the effective potential holds over a wide range of random potentials in two dimensions. In the remainder of the first two rows of Figure~\ref{fg:ratio} we show the same results also for the Bernoulli potential of Figure~\ref{fg:bernoulli} and the correlated potential shown in Figure~\ref{fg:correlated}. (The correlated potential was constructed with circulant embedding as discussed above for one dimension, except that we of course used the two-dimensional discrete Fourier transform, and we took as aperture function $\sigma \chi(d |t|)$ with $\sigma=4$ and $d=0.05$, where $\chi$ is the characteristic function of the unit interval, instead of $\sigma\exp(-d|t|)$ as we took previously.) The three pairs of plots in the final two rows of the figure show the case of a uniformly random $40\x40$ potential with values taken between $0$ and $4$, $16$, and $64$, respectively. The quotient is quite close to being constant with value $1.5$, particularly in the last two cases. In the first case, with the lowest disorder, the ratio drifts away from $1.5$ to about $1.75$ after approximately 50 eigenvalues.

\begin{figure}[htbp]
   \centerline{\includegraphics[width=60mm]{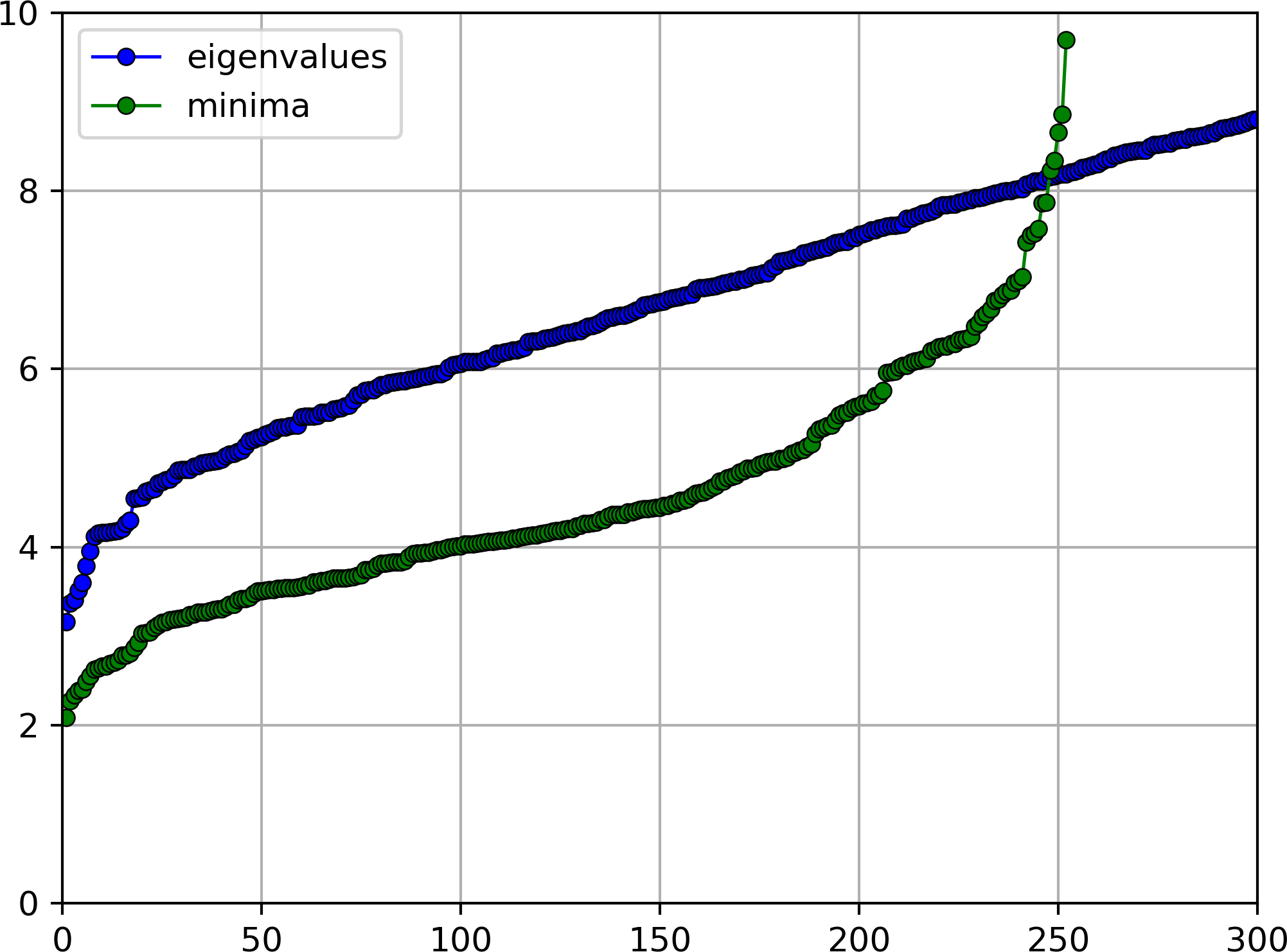}
   \includegraphics[width=60mm]{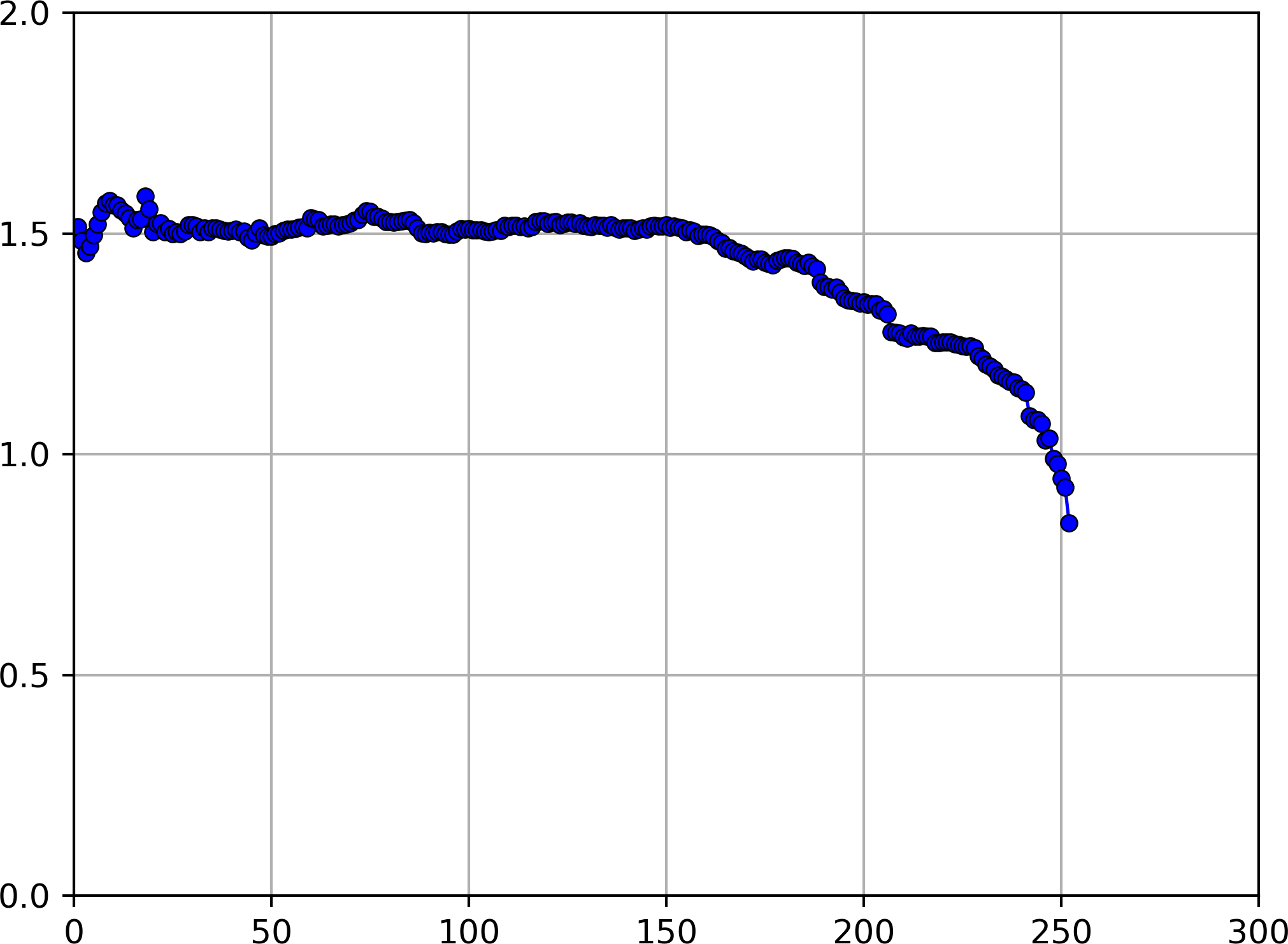}}%
 \caption{A plot of \emph{all} the minima of $W$ versus the eigenvalues, and, on right, their ratio, for the same potential realization as the fifth case shown in Figure~\ref{fg:ratio}.}
\label{fg:ratioall} 
\end{figure}

Of course, the minima of the effective potential can only predict a limited number of eigenvalues.  Indeed, $W$ has only a finite number of local minima, while there are infinitely many eigenvalues.  In Figure~\ref{fg:ratioall} we revisit the fifth case shown in Figure~\ref{fg:ratio},
with a uniformly random $40\x 40$ potential with values in $[0,16]$. For this realization of the potential, $W$
has exactly 252 local minima.  These are plotted alongside the first 300 eigenvalues in the left-hand side of Figure~\ref{fg:ratioall}, and their ratios with the corresponding eigenvalues are plotted on the left.  We see that, in this case, the ratio of $1.5$ remains quite accurate for more than 150 of the 252 minima.

In Table~\ref{tb:ratio} we display the mean and standard deviation of these ratios computed for the first 10, 50, and 100 eigenvalues for each of the six potentials of Figure~\ref{fg:ratio}. The main observation is that, across all this range, the ratio stays quite close to $1.5$.

\begin{table}[htb]
\centering
\renewcommand{\arraystretch}{1.2}
\begin{tabular}{@{}lrrcl@{\hskip-1pt}lcl@{\hskip-1pt}lcl@{\hskip-1pt}l@{}}
\toprule
& & & \multicolumn{2}{c}{10 eigs} & & \multicolumn{2}{c}{50 eigs} & & \multicolumn{2}{c}{100 eigs}\\
\cmidrule{4-5} \cmidrule{7-8} \cmidrule{10-11}
potential & \ct{$n_{\text{c}}$} & & \ct{mean} & \ct{SD} & &  \ct{mean} & \ct{SD} & & \ct{mean} & \ct{SD} \\
\midrule
uniform $[0,20]$ & 80 & & 1.519 & 0.024 & & 1.520 & 0.018 & & 1.524 & 0.017 \\
Bernoulli & 80 & & 1.585 & 0.034 & & 1.607 & 0.041 & & 1.646 & 0.049 \\
correlated & 80 & & 1.519 & 0.018 & & 1.532 & 0.028 & & 1.530 & 0.022 \\
uniform $[0,4]$ & 40 & & 1.479 & 0.018 & & 1.510 & 0.034 & & 1.582 & 0.087 \\
uniform $[0,16]$ & 40 & & 1.518 & 0.041 & & 1.515 & 0.026 & & 1.516 & 0.021 \\
uniform $[0,64]$ & 40 & & 1.503 & 0.025 & & 1.511 & 0.024 & & 1.477 & 0.043 \\
\bottomrule
\end{tabular}
 \caption{The mean and standard deviation for the ratio of first 10, 50, and 100 eigenvalues to the corresponding minima values of $W$.  For a wide range of potentials in two dimensions the ratio is roughly $1.5$ across many eigenvalues.}
\label{tb:ratio}
\end{table}

Table~\ref{tb:ratio1d} is similar, but shows the results for a variety of potentials in one dimension. Again, we see that the ratio of the eigenvalue to the corresponding minimum values of the effective potential is roughly constant, as it was in two dimensions. However the constant value we find in one dimension is about $1.25$ or $1.3$ rather than the value of $1.5$ we saw in two dimensions.

\begin{table}[htb]
\centering
\renewcommand{\arraystretch}{1.2}
\begin{tabular}{@{}lrrcl@{\hskip-1pt}lcl@{\hskip-1pt}lcl@{\hskip-1pt}l@{}}
\toprule
& & & \multicolumn{2}{c}{10 eigs} & & \multicolumn{2}{c}{25 eigs} & & \multicolumn{2}{c}{50 eigs}\\
\cmidrule{4-5} \cmidrule{7-8} \cmidrule{10-11}
potential & \ct{$n_{\text{c}}$} & & \ct{mean} & \ct{SD} & &  \ct{mean} & \ct{SD} & & \ct{mean} & \ct{SD} \\
\midrule
uniform $[0,4]$ & 256 & & 1.303 & 0.026 & & 1.321 & 0.029 & & 1.300 & 0.067 \\
 & 1024 & & 1.302 & 0.019 & & 1.301 & 0.020 & & 1.304 & 0.022 \\
uniform $[0,16]$ & 256 & & 1.322 & 0.023 & & 1.308 & 0.031 & & 1.240 & 0.099 \\
 & 1024 & & 1.274 & 0.018 & & 1.296 & 0.031 & & 1.294 & 0.027 \\
Bernoulli & 256 & & 1.301 & 0.033 & & 1.316 & 0.050 & & 1.296 & 0.130 \\
 & 1024 & & 1.262 & 0.014 & & 1.272 & 0.026 & & 1.266 & 0.073 \\
correlated & 256 & & 1.335 & 0.055 & & 1.404 & 0.090 & & 1.310 & 0.171 \\
 & 1024 & & 1.280 & 0.018 & & 1.286 & 0.017 & & 1.303 & 0.042 \\
\bottomrule
\end{tabular}
 \caption{The mean and standard deviation for the ratio of first 10, 25, and 50 eigenvalues to the corresponding
 minima values of $W$, tabulated here for 8 different types of random potentials in one dimension.}
\label{tb:ratio1d}
\end{table}

\begin{figure}[htbp]
   \centerline{\includegraphics[width=73mm]{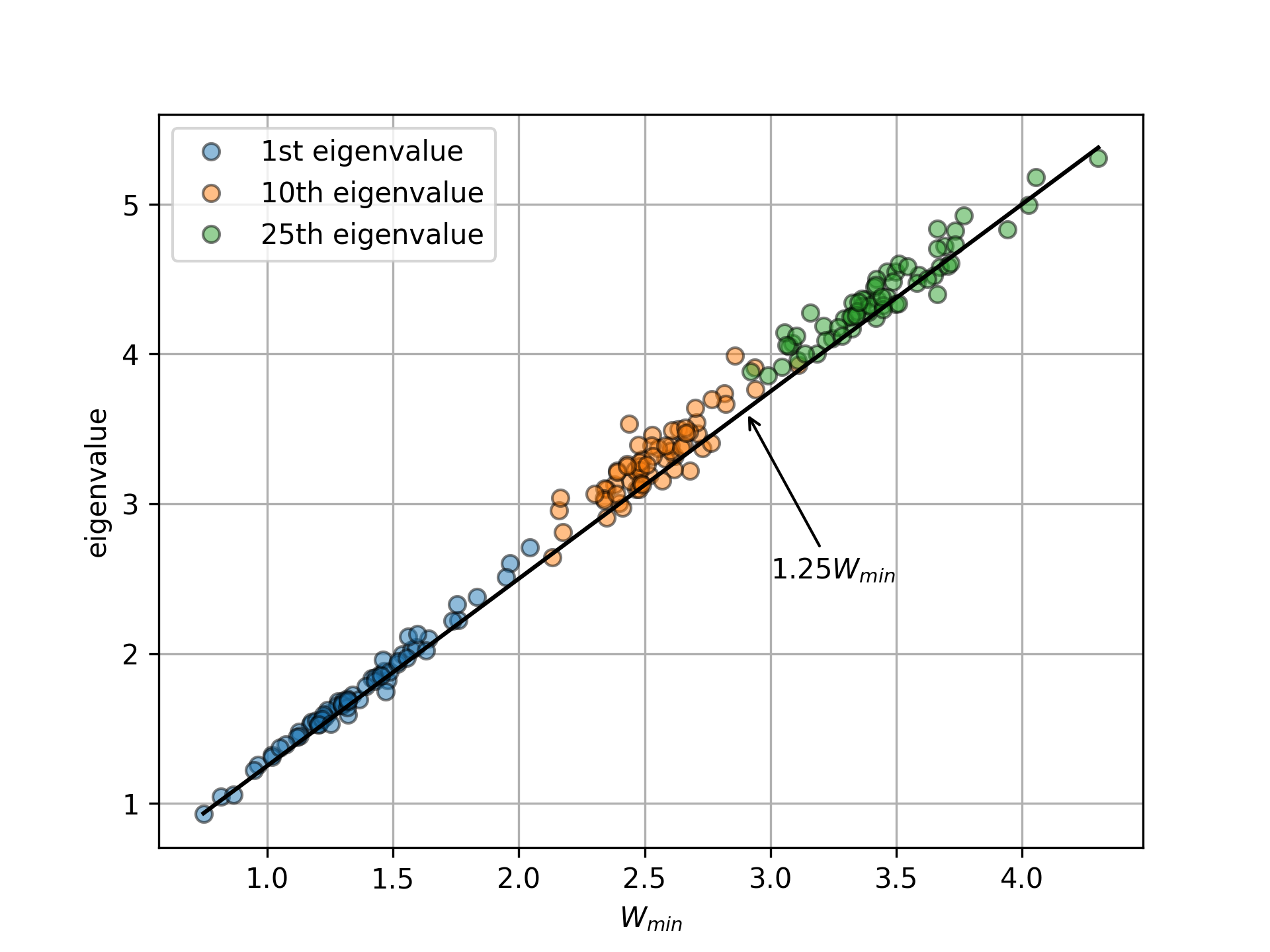}\hspace{-8mm}
   \includegraphics[width=73mm]{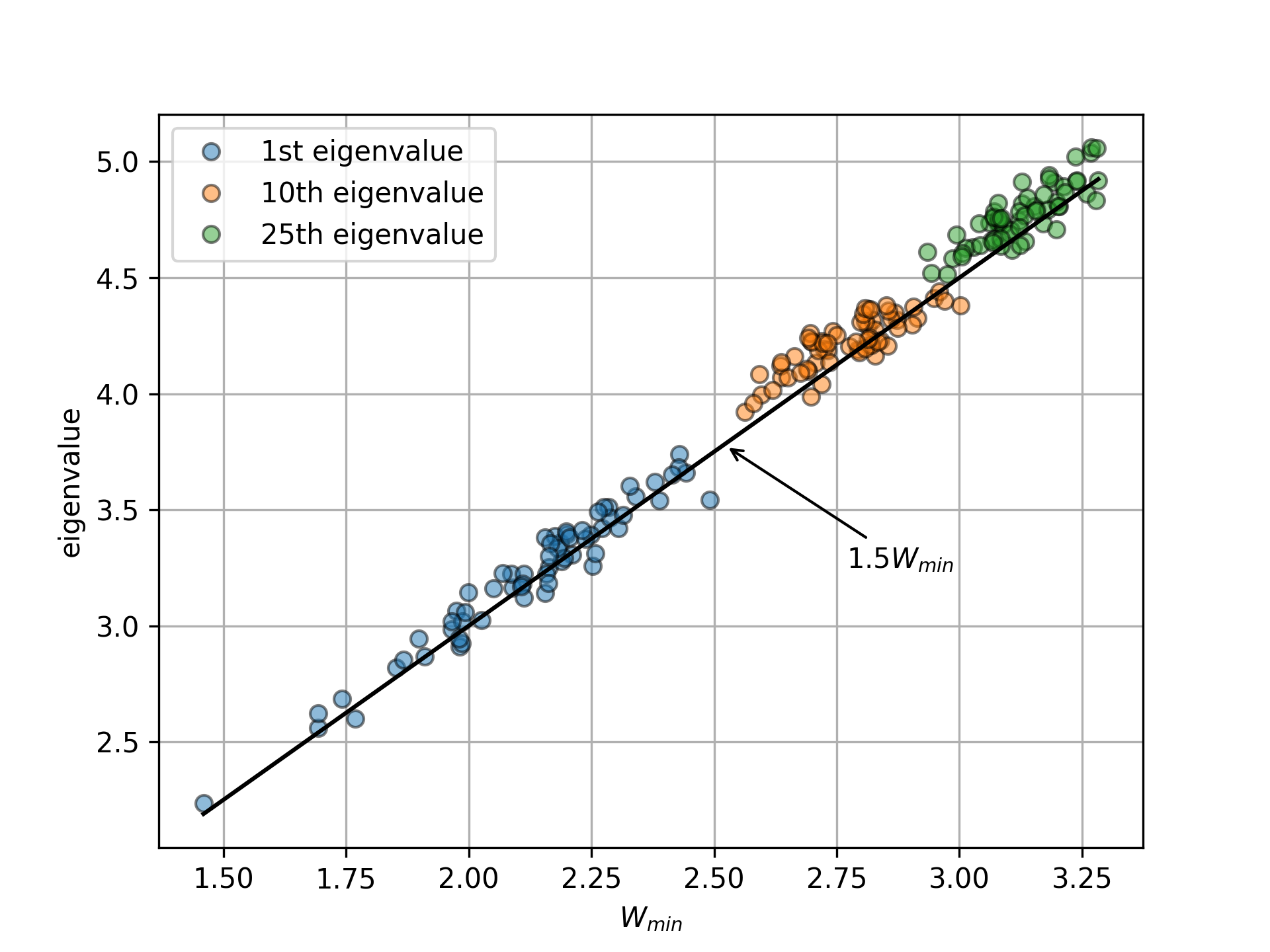}}%
 \caption{The 1st, 10th, and 25th eigenvalues versus the corresponding minima values of the effective potential,
 for $64$ independent realizations of a random potential.  Left in 1D where the black line shown
 is $\lambda=1.25 W_{\text{min}}$, right in 2D with $\lambda=1.5 W_{\text{min}}$.}
\label{fg:scatter} 
\end{figure}

Finally, in Figures~\ref{fg:scatter} we plot the 1st, 10th, and 25th eigenvalues
versus the corresponding minima values of $W$ for numerous different realizations of a random potential.
The first figure displays 64 realizations of a 1D potential on $[0, 256]$ with $256$ values selected uniformly iid from $[0, 16]$, while the second figure displays $64$ realizations of a 2D potential on $[0, 40]\x[0,40]$ with
$1,600$ random values again chosen uniformly iid from $[0, 16]$.  We see that in the first case the points line up
well along the line $\lambda=1.25 W_{\text{min}}$, and in the second along the line $\lambda=1.5 W_{\text{min}}$.

From this and other evidence, we conclude that in many cases
\begin{equation}\label{evalap}
\lambda\approx (1 + \frac{n}{4})\,W_{\text{min}}.
\end{equation}
Here $\lambda$ is one of the lower eigenvalues of $H$, for which the corresponding eigenfunction is localized to a subdomain $\Omega_0$, $W_{\text{min}}$ is the minimum value of the effective potential on that subdomain, and $n$ is the number of dimensions (thus far $1$ or $2$).  The constant $1+n/4$, i.e.,
$1.25$ in 1D and $1.5$ in 2D, is a rough approximation in accord with our observations and which we now further justify heuristically.  It is remarkable that this constant, though dimension-dependent, is independent of the specific realization of the potential and of the parameters of its probability distribution, the size of the domain, etc.

We now give some heuristic support of the eigenvalue approximation \eqref{evalap}.  Our argument will be rather crude, and we do not claim it fully explains the numerical evidence presented above. Let $\psi$ denote one of the eigenfunctions associated to a smaller eigenvalue. We assume $\psi$ to be localized, i.e., essentially supported in a small subdomain $\Omega_0$, for which it is the fundamental Dirichlet eigenfunction. We also assume that on the subdomain $\Omega_0$ the landscape function $u$ is well approximated by a constant multiple of the fundamental eigenfunction $\psi$ of the subdomain.  This is roughly supported by experimental results such as shown in Figure~{\ref{fg:le}}.  Another supporting argument comes from the expansion of the constant $1$ on $\Omega_0$ in terms of the Dirichlet eigenfunctions of the domain, retaining
only the first term $c_0\psi$, and dropping the terms coming from the eigenfunctions which change sign.  Then $u\approx (c_0/\lambda)\psi$, indeed a multiple of $\psi$.  We may use these two assumptions, together with the definition $Hu=1$ of the landscape function, to approximate the Rayleigh quotient:
$$
\lambda= \frac{\int_{\Omega} \psi H\psi\,dx}{\int_{\Omega}\psi^2\,dx}
\approx \frac{\int_{\Omega_0} \psi H\psi\,dx}{\int_{\Omega_0}\psi^2\,dx}
\approx \frac{\int_{\Omega_0} u Hu\,dx}{\int_{\Omega_0}u^2\,dx} = \frac{\int_{\Omega_0} u\,dx}{\int_{\Omega_0}u^2\,dx}.
$$
Next we assume that on $\Omega_0$ the landscape function $u$ (or $\psi$ which we are supposing is a constant multiple of $u$ there) can be approximated by the simplest sort of positive bump-like function, the positive part of a concave quadratic function.  After rotating and translating the coordinate system, this means that
$$
u\approx u_{\text{max}}[1-\sum(x_i/a_i)^2 ] \text{ on $\Omega_0\approx\{\,x\in\R^n\,|\, \sum(x_i/a_i)^2 \le 1\,\}$},
$$
for some positive constants $a_i$. Thus our approximation to the eigenvalue is
$$
\lambda\approx \frac{\int_{\Omega_0} u\,dx}{\int_{\Omega_0} u^2\,dx} =
\frac c{u_{\text{max}}},\quad c = \frac{\int_{\Omega_0} [1-\sum(x_i/a_i)^2 ]\,dx}
{\int_{\Omega_0}[1-\sum(x_i/a_i)^2 ]^2\,dx}.
$$
Finally, we compute $c$ using the change of variables $\hat x_i = x_i/a_i$ to convert the integrals in the numerator and denominator into integrals over the unit ball $B$ which can be computed with polar coordinates.  This gives $c= 1+n/4$ (independent of the values of the $a_i$ and so of the size and shape of the ellipsoid).  Since $1/u_{\text{max}}= W_{\text{min}}$, this indeed gives the approximation \eqref{evalap}.

\begin{figure}[htbp]
  \centerline{%
  \includegraphics[height=2.in]{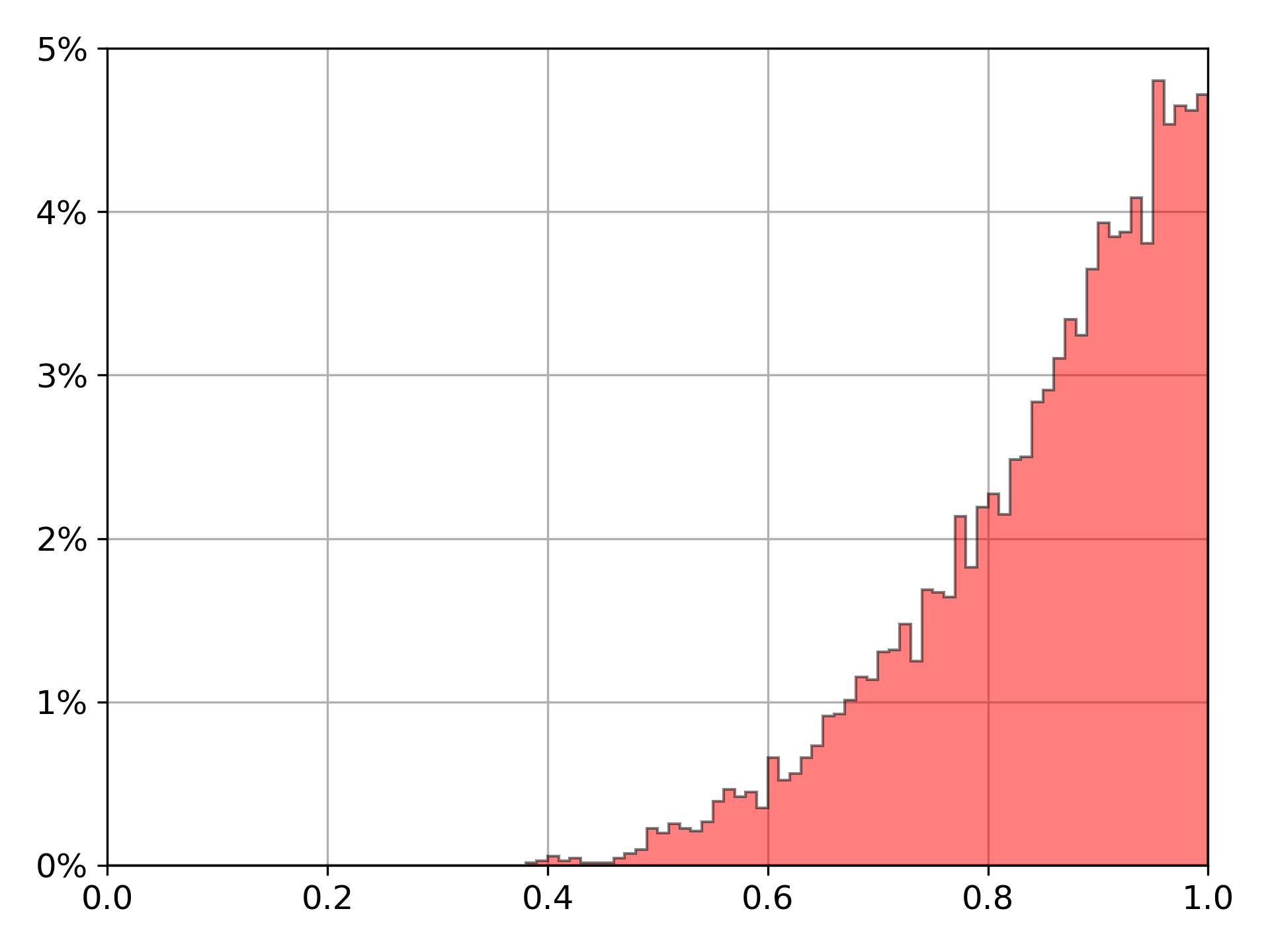}%
  \quad
  \includegraphics[height=2.in]{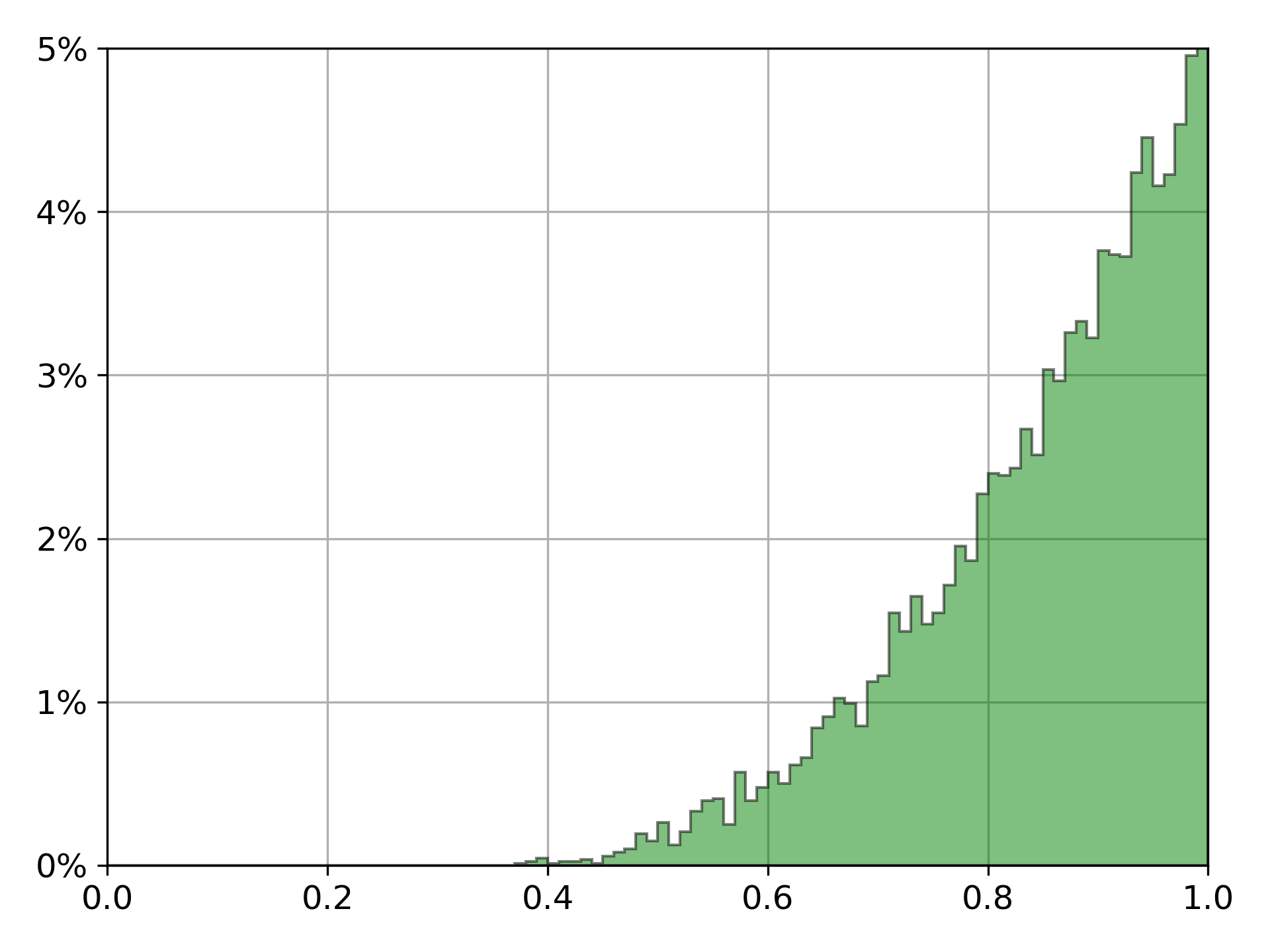}}
  \centerline{\includegraphics[height=2.in]{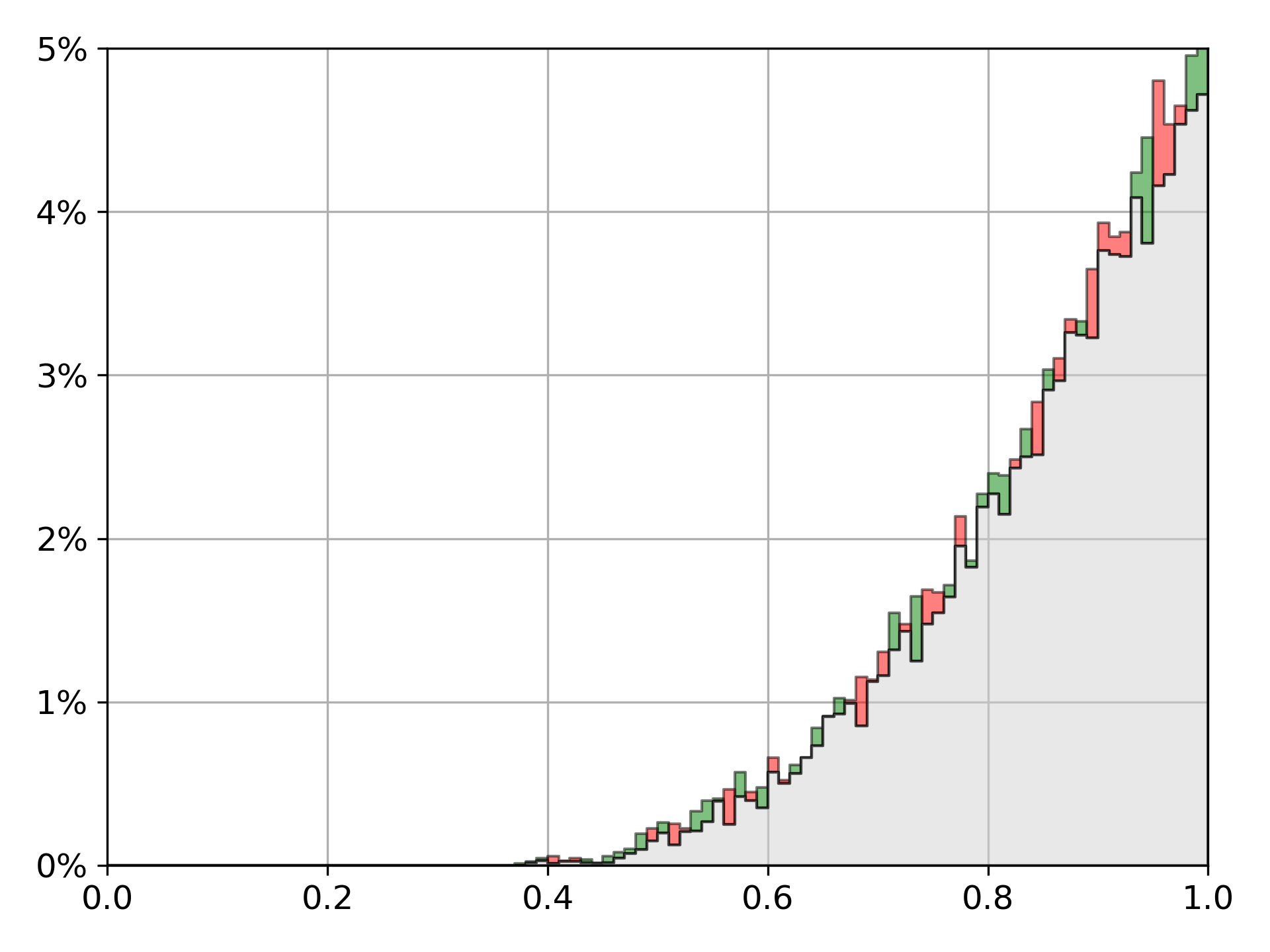}}
 \caption{Density of states on $[0,1]$ for the case of a uniformly random potential with $2^{19}$ pieces, displayed with 100 bins. On top left, histogram of the actual eigenvalues. On top right, histogram of $1.25$ times the local minima of $W$. The bottom highlights the differences between the two.}
\label{fg:dos} 
\end{figure}

The results we have shown in Figure~\ref{fg:ratio} and Tables~\ref{tb:ratio} and \ref{tb:ratio1d} demonstrate that the approximation \eqref{evalap} can be used to estimate 100 eigenvalues with errors of a few percent. Note that it is nonetheless very cheap to apply \eqref{evalap}, the cost being that of solving a single source problem and the extraction of some maxima, much less than the cost of computing many eigenvalues. As another example of the utility of \eqref{evalap} we now use it to approximate the density of states in the interval $[0,1]$ of the 1D Schr\"odinger operator for which the piecewise constant potential has $2^{19}=524,288$ pieces. (Specifically, we compute on the interval $[0, 2^{19}]$ and assign the random values to unit subintervals uniformly iid in $[0,4]$.) We display the DOS as a histogram with 100 bins. The top left plot in Figure~\ref{fg:dos} shows the actual density of states, requiring the computation of all 7,122 of the eigenvalues which belong to the interval $[0,1]$.  The finite element mesh we use to approximate the Schr\"odinger operator in this case has $10\x2^{19}$ elements, and we use piecewise cubic finite elements, so that the problem has about 15.7 million degrees of freedom.  The calculation of 7,000 eigenvalues is thus a very large computation. It required about 40 CPU hours on a workstation with an Intel Core i7-4930K processor, using sprectral slicing and the Krylov-Schur method of SLEPc. However, an accurate approximation of the density of states can be obtained quickly using the effective potential without resorting to the computation of any eigenvalues.  This approximation is shown in the top right plot of Figure~\ref{fg:dos}. It is 
a histogram of all those values of $1.25 W_{\text{min}}$ which belong to the same interval $[0, 1]$ (8,800 in all). The computation of these values is much less demanding.  It required slightly over 5 minutes of CPU time of the same workstation, i.e., was about 480 times faster. The two histograms are compared in the bottom plot.

\subsection{Eigenvalues from a variant of Weyl's law}\label{ssec:eval2}
\begin{figure}[htbp]
\centerline{%
  \includegraphics[width=2.5in]{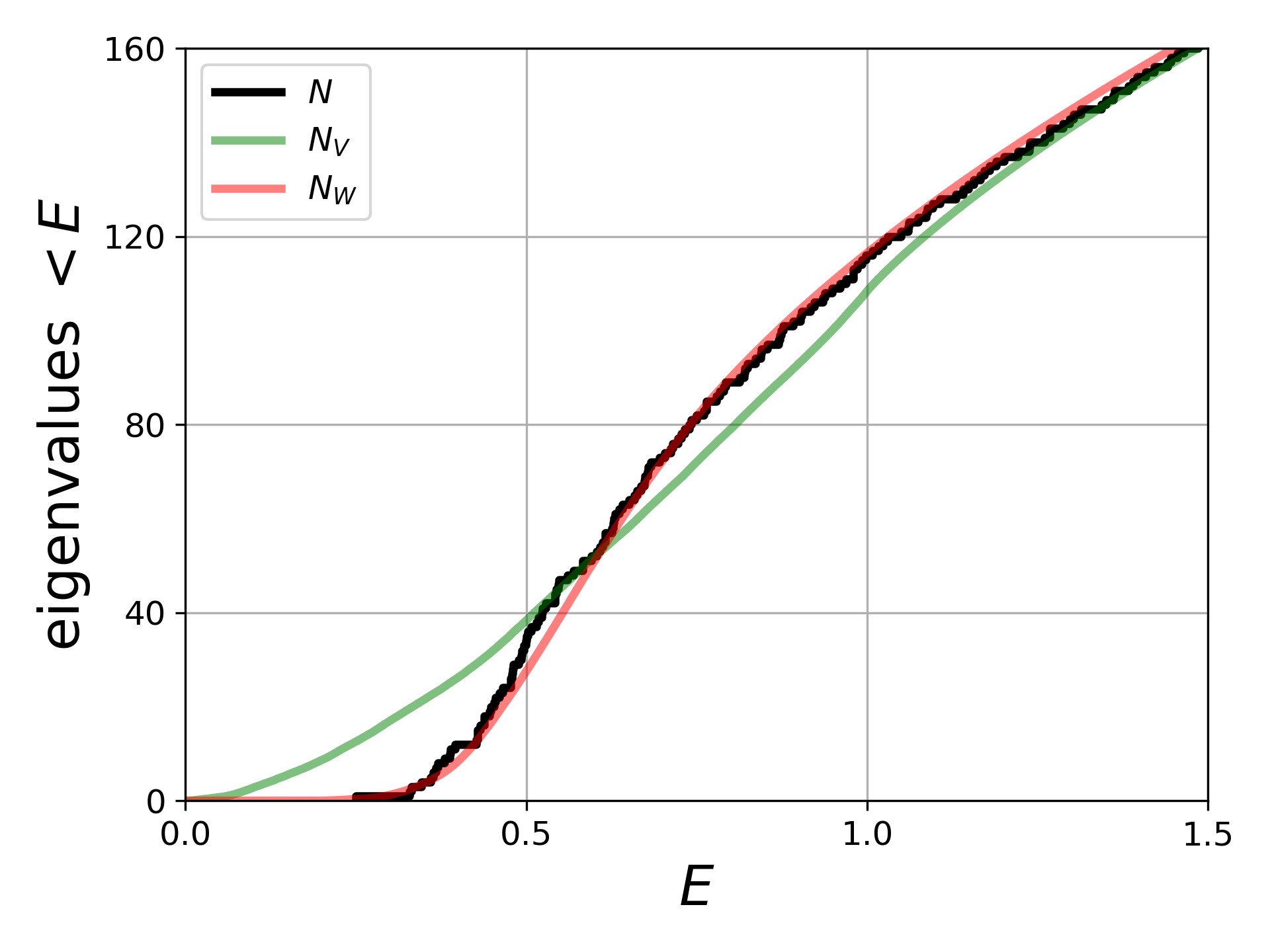} \quad
  \includegraphics[width=2.5in]{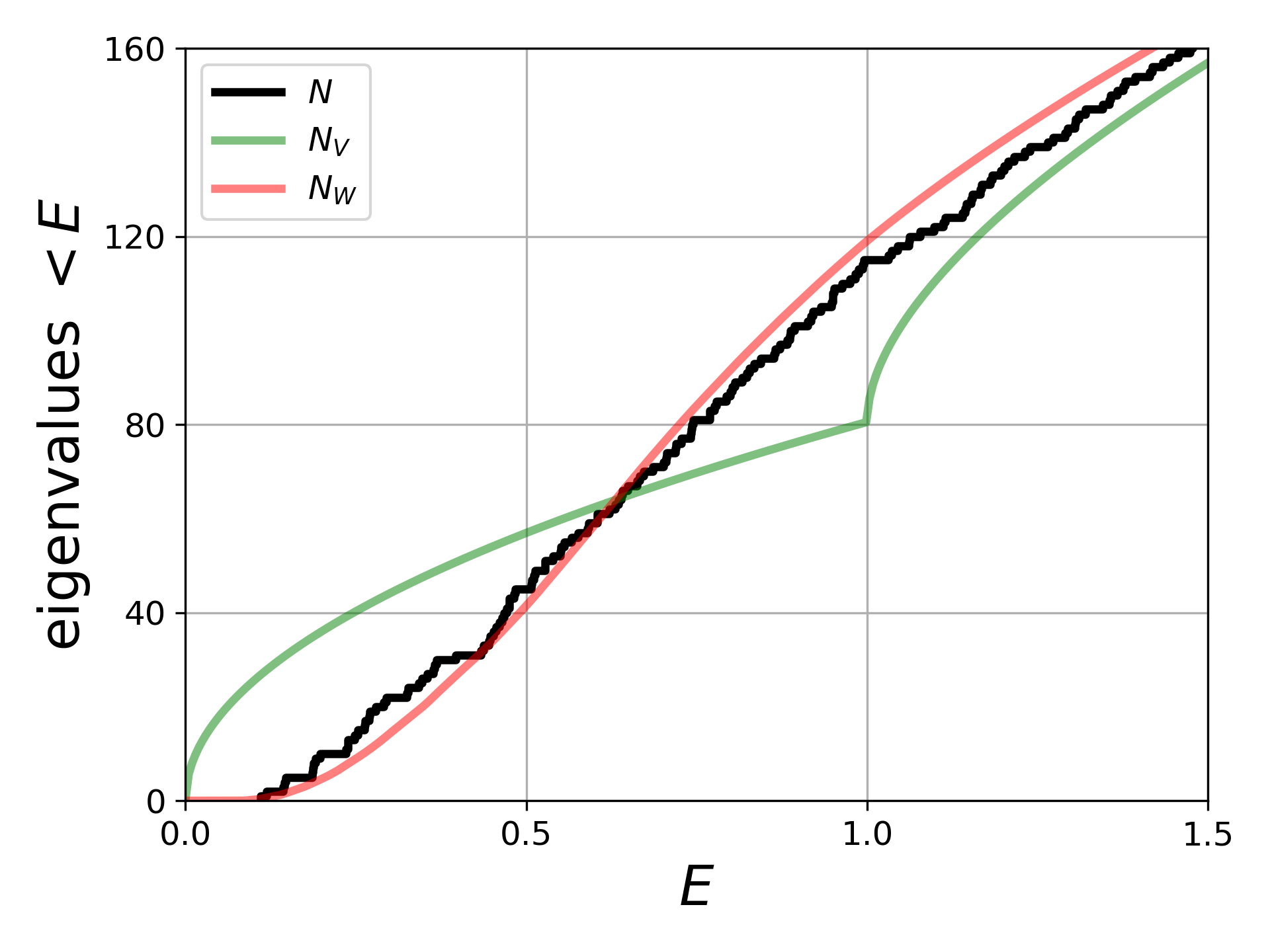}}
\centerline{%
  \includegraphics[width=2.5in]{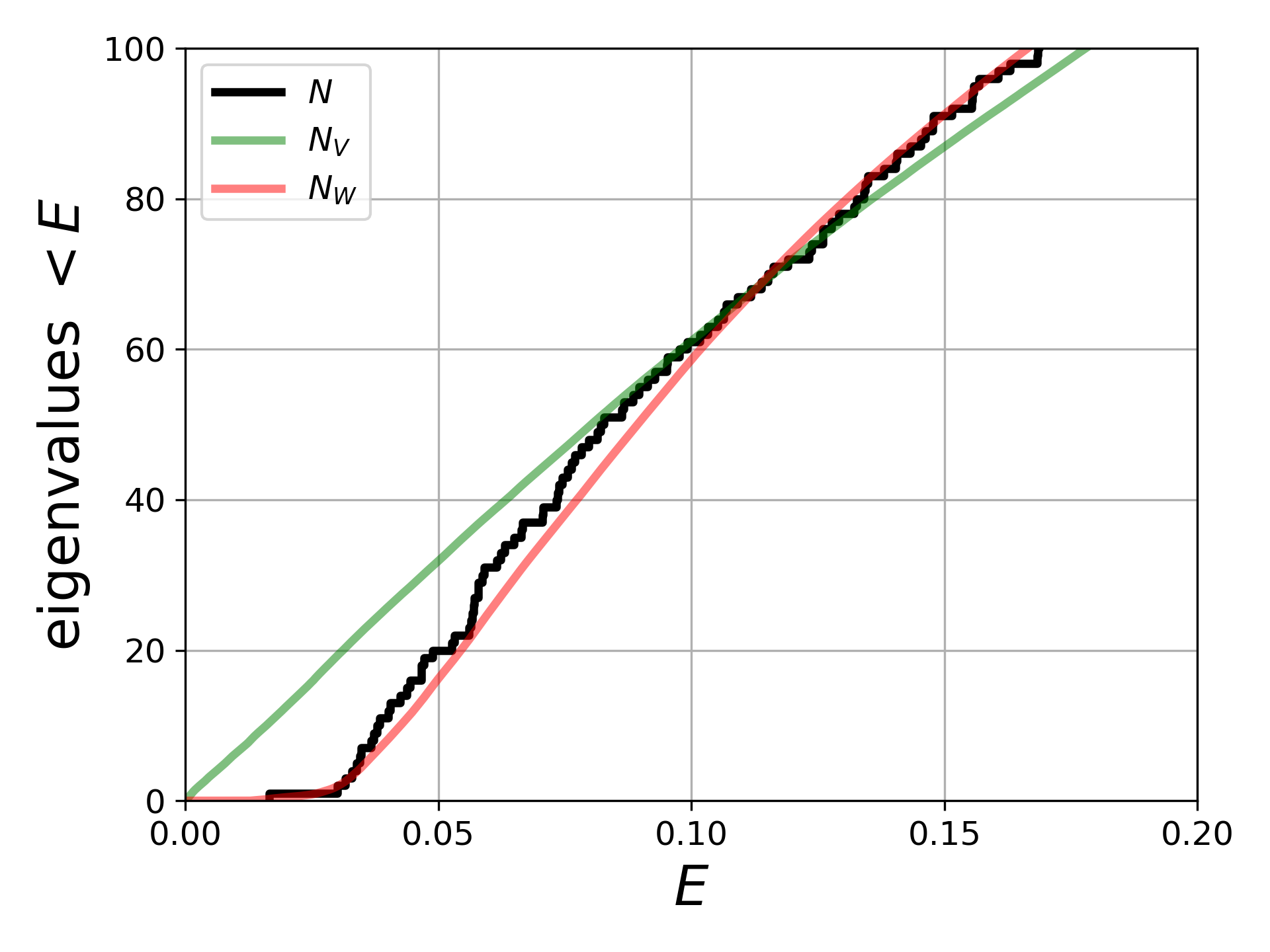} \quad
  \includegraphics[width=2.5in]{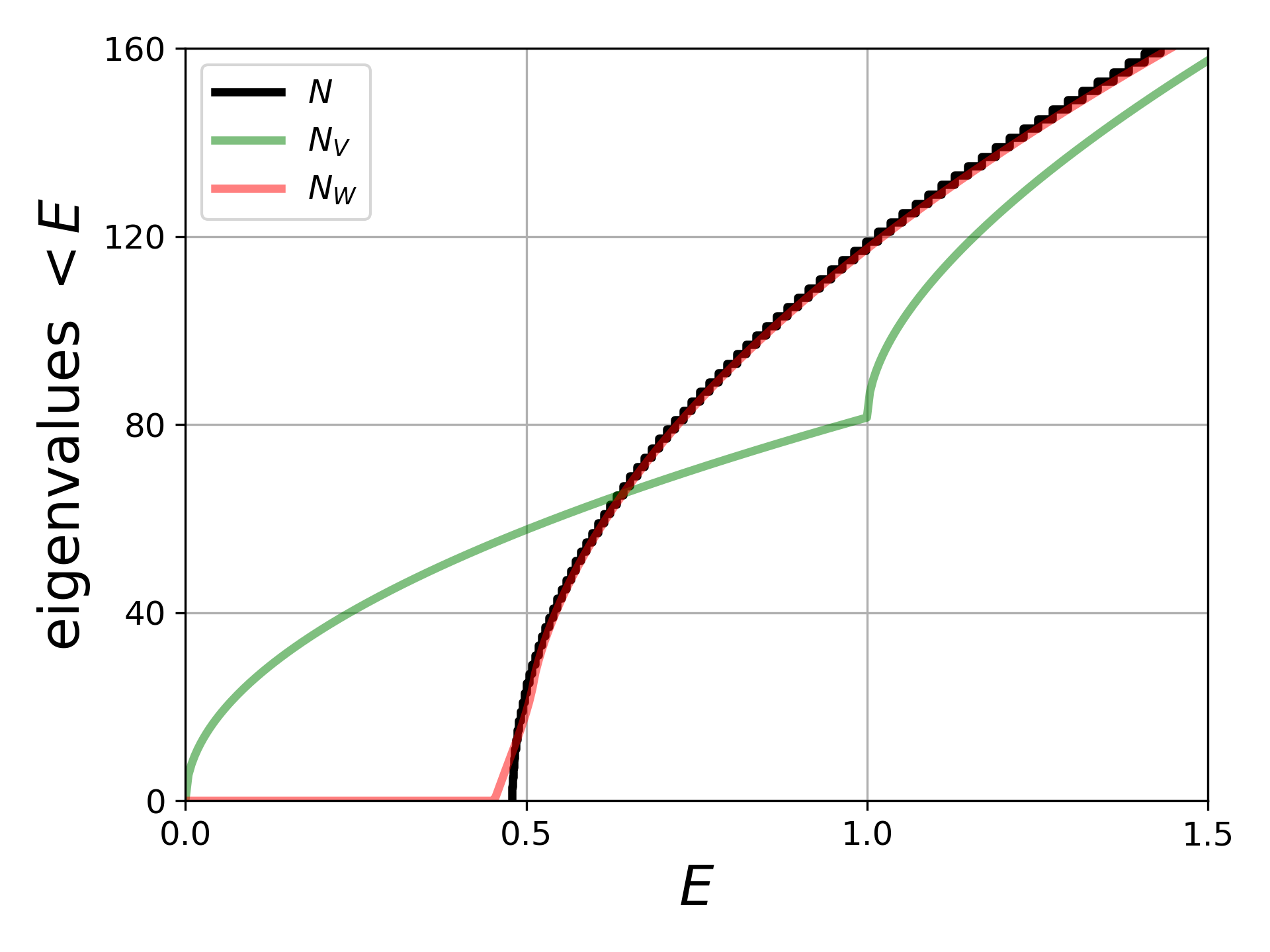}}
\caption{The eigenvalue counting function $N$, the Weyl's law approximation $N_V$, and the effective Weyl's law approximation $N_W$ for some potentials in one dimension. Top row: uniform and Boolean random potentials. Bottom row: correlated and periodic Boolean potential.}
\label{fg:weyl1d} 
\end{figure}

The approximation \eqref{evalap} can be used to predict, at most, one eigenvalue per local minimum of $W$.  Intuitively, it approximates a localized eigenfunction by the fundamental mode of the well around the local minimum. Now we present an alternative approach, which again relies on the effective potential $W$, but which gives some sort of prediction for all the eigenvalues.  For large eigenvalues, it provides information similar to Weyl's law. Recall that Weyl's law for the Schr\"odinger equation is an asymptotic formula for the eigenvalue counting function:
\begin{equation}\label{weyl}
N(E)\sim N_V(E):= (2\pi)^{-n}\operatorname{vol}\{\,(x,\zeta)\in\Omega\x\R^n\,|\,
V(x) + |\zeta|^2\le E\,\}\text{ as $E\to\infty$},
\end{equation}
where the volume term is the $2n$-dimensional measure of the indicated subset of the
phase space $\Omega\x\R^n$. Assuming smoothness and growth conditions for the potential, Weyl's law holds asymptotically as $E$ tends to $+\infty$ and so for a large number of eigenvalues~\cite[Theorem~6.8]{zworski}.  Inverting the counting function, we can thus view Weyl's law as furnishing an approximation of the $n$th eigenvalue, which is asymptotically valid for $n$ large.   Weyl's law is generally not expected to be accurate for a small number of eigenvalues.  However, experimentally we have found that, for the sorts of random potentials considered in this paper, a variant of Weyl's law invoking the effective potential $W$ gives very good results right down to the first few eigenvalues, while remaining asymptotically correct. The variant, which we shall refer to as the \emph{effective Weyl's law}, is obtained by simply replacing the potential $V$ in \eqref{weyl} by the effective potential $W$:
$$
N_W(E) = (2\pi)^{-n}\operatorname{vol}\{\,(x,\zeta)\in\Omega\x\R^n\,|\,
W(x) + |\zeta|^2\le E\,\}.
$$
Figure~\ref{fg:weyl1d} compares the true eigenvalue counting function $N$ (shown in black), Weyl's law (green), and the effective Weyl's law (red), for four different types of potential. The first is uniformly random iid with 512 pieces and values in $[0,1]$. The second is Bernoulli potential where the 512 random values are either $0$ or $1$, each with probability $1/2$. The third is a correlated Gaussian squared potential like that of Figure~\ref{fg:corr1d}. The fourth is quite different: the 512 Boolean values $0$ and $1$ are assigned alternately.  For this potential, there is no localization. Nonetheless, we see that, in each case, Weyl's law becomes a good approximation only after 100 or so eigenvalues and, in every case, it incorrectly predicts many eigenvalues in the interval from $0$ to the least eigenvalue. By contrast, the effective Weyl's law provides a very good approximation of the counting function for many eigenvalues, starting from the first. It is revealing to compare the two potentials which take on only Boolean values (the second and the fourth). Because the classical Weyl's law is unaffected by rearrangement of the potential, it gives the same prediction for the counting function in both cases. But the actual counting functions differ very significantly, a fact which is well captured by the effective Weyl's law. (Similar results were published in \cite{ADJMF2016}.)

\begin{figure}[htbp]
\centerline{%
  \includegraphics[width=2.5in]{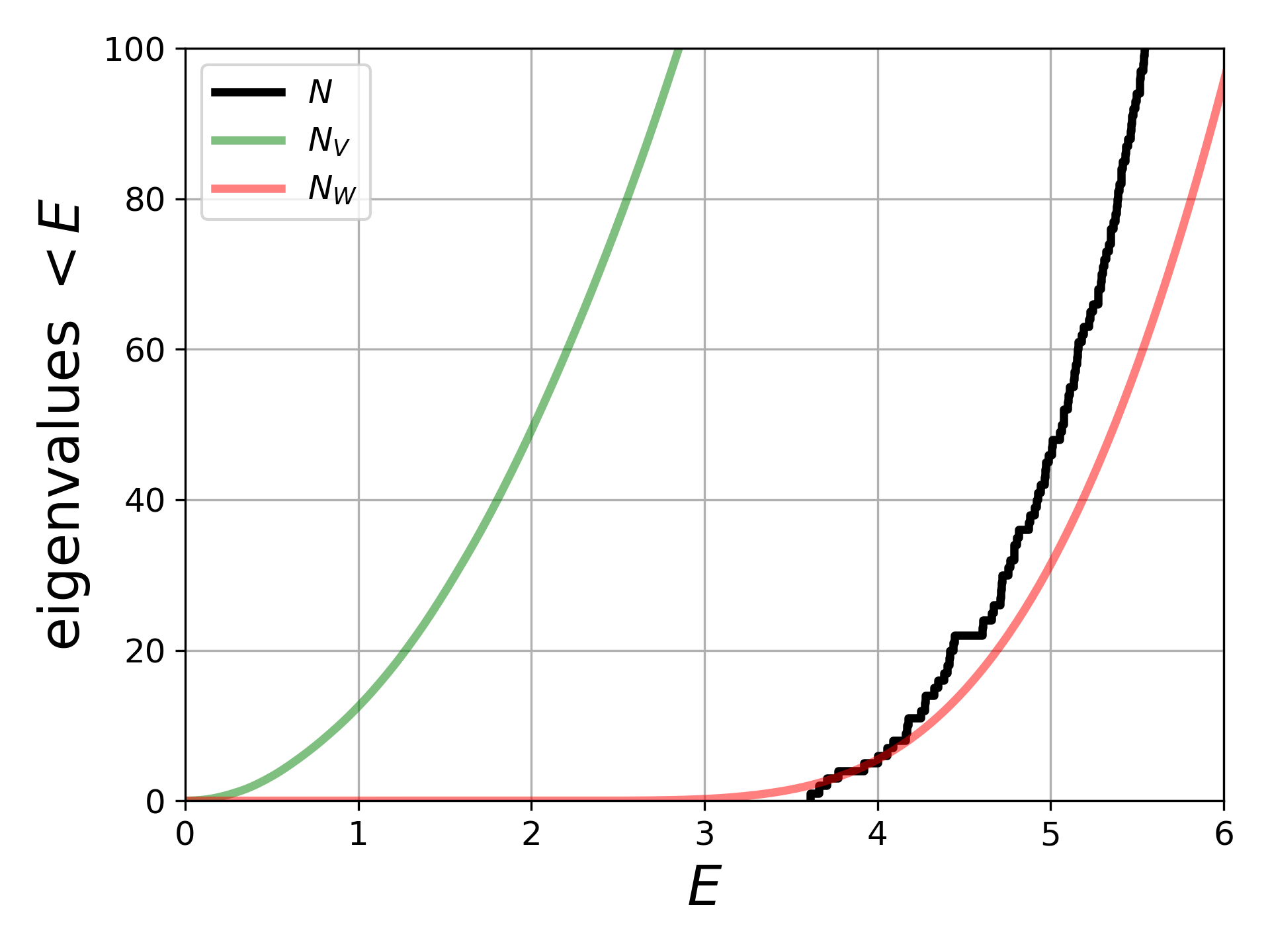} \quad
  \includegraphics[width=2.5in]{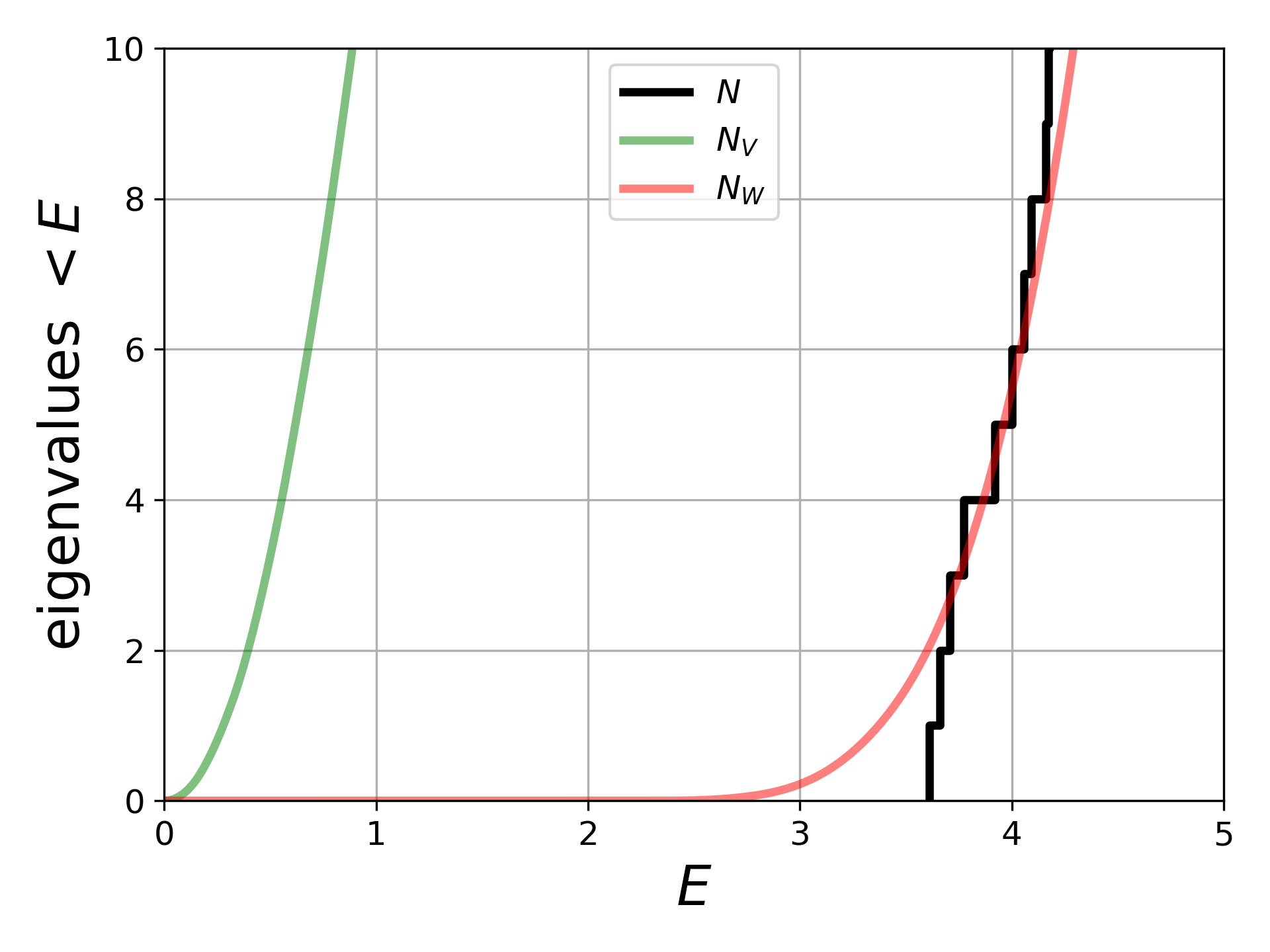}}
\caption{The eigenvalue counting function $N$, the Weyl's law approximation $N_V$,  and the effective Weyl's law approximation $N_W$ for the 2D potential of Figure~\ref{fg:localization}, showing the first 100 eigenvalues on left and restricting to the first 10 on right.}
\label{fg:weyl2d} 
\end{figure}

Finally, in Figure~\ref{fg:weyl2d} we show similar results for a single 2D potential, namely the uniformly random $80\x80$ potential of Figure~\ref{fg:localization}. The first plot shows the first 100 eigenvalues, while the second zooms in on the first 10 eigenvalues.
The predictive power of the effective Weyl's law does not seem to be as great as in one dimension, but again it displays a great
improvement over the classical Weyl's law for small eigenvalues.

\section{Conclusion}\label{sec:conc}
We have demonstrated numerically that the effective potential, defined as the reciprocal of the localization landscape function,
accurately captures a great deal of information about the localization properties of a 
random potential, and shown how to employ it to predict eigenvalues and eigenfunctions.
These predictions are attained by a solving a single PDE source problem,
without the direct solution of any eigenvalue
problems, and so at a very low computational price.  The wells of the effective potential reveal the main localization subdomains, and the values of its minima are found to be very good predictors of the corresponding fundamental eigenvalues. We have tested this approach on piecewise constant
potentials with several types of random distributions, uniformly random and Bernoulli, 
and with a certain correlated distribution, both in one and two dimensions. We have further used the effective potential
to predict the density of states and obtain good precision even for small eigenvalues, something
which is not attained by the classical Weyl law asymptotics.
In highly demanding computations where the Schr\"odinger equation has to be solved for a large number eigenfunctions and eigenvalues (as for instance in semiconductor physics), the resulting computational efficiency makes it now possible to reproduce numerically, to analyze, and to understand the behavior of quantum disordered materials.
\bibliographystyle{siamplain}
\bibliography{specpred}

\end{document}